\documentclass{elsarticle}
  
\usepackage[margin=1in]{geometry}

\usepackage{bbm}
\usepackage{graphicx}
\usepackage{pstool}
\usepackage{wrapfig}
\usepackage{caption}
\usepackage{subcaption}
\usepackage[section]{placeins} 

\usepackage{fancyhdr} 
\usepackage{lastpage} 
\usepackage{extramarks} 

\usepackage{xcolor} 
\usepackage{enumerate} 
\usepackage{paralist} 
\usepackage{amsmath, amsthm, amssymb, mathtools}
\usepackage{mathabx, pifont, stmaryrd} 
\usepackage[explicit]{titlesec} 
\usepackage{etoolbox} 
\usepackage{bibentry} 
\makeatletter\let\saved@bibitem\@bibitem\makeatother 
\usepackage[colorlinks, bookmarksopen, bookmarksnumbered,
            citecolor=red,urlcolor=red]{hyperref} 
\makeatletter\let\@bibitem\saved@bibitem\makeatother 

\usepackage{algorithm}
\usepackage{algorithmic}

\newtheorem{remark}{Remark}

\theoremstyle{definition}

\usepackage{bm} 
\usepackage{amsfonts} 

\DeclareMathOperator*{\argmax}{arg\,max}
\DeclareMathOperator*{\argmin}{arg\,min}

\newcommand{\func}[3]{\ensuremath{#1 : #2 \rightarrow #3}}

\newcommand{\norm}[1]{\ensuremath{\left\| #1 \right\|}}

\newcommand{\suchthat}{\mathrel{}\middle|\mathrel{}}

\newcommand{\pder}[2]{\ensuremath{\frac{\partial #1}{\partial #2}}}

\newcommand{\Ccal}{\ensuremath{\mathcal{C}}}
\newcommand{\Dcal}{\ensuremath{\mathcal{D}}}

\newcommand{\Gcal}{\ensuremath{\mathcal{G}}}

\newcommand{\Ical}{\ensuremath{\mathcal{I}}}

\newcommand{\Ncal}{\ensuremath{\mathcal{N}}}
\newcommand{\Ocal}{\ensuremath{\mathcal{O}}}

\newcommand{\Gbb}{\ensuremath{\mathbb{G}}}

\newcommand{\Rbb}{\ensuremath{\mathbb{R} }}

\newcommand\Fbm{{\ensuremath{\bm{F}}}}

\newcommand\Hbm{{\ensuremath{\bm{H}}}}
\newcommand\Ibm{{\ensuremath{\bm{I}}}}
\newcommand\Jbm{{\ensuremath{\bm{J}}}}

\newcommand\Nbm{{\ensuremath{\bm{N}}}}

\newcommand\Pbm{{\ensuremath{\bm{P}}}}
\newcommand\Qbm{{\ensuremath{\bm{Q}}}}
\newcommand\Rbm{{\ensuremath{\bm{R}}}}
\newcommand\Sbm{{\ensuremath{\bm{S}}}}
\newcommand\Tbm{{\ensuremath{\bm{T}}}}
\newcommand\Ubm{{\ensuremath{\bm{U}}}}

\newcommand\Xbm{{\ensuremath{\bm{X}}}}

\newcommand\ubm{{\ensuremath{\bm{u}}}}
\newcommand\vbm{{\ensuremath{\bm{v}}}}
\newcommand\wbm{{\ensuremath{\bm{w}}}}
\newcommand\xbm{{\ensuremath{\bm{x}}}}

\newcommand\zbm{{\ensuremath{\bm{z}}}}

\newcommand\mubold{{\ensuremath{\boldsymbol{\mu}}}}

\newcommand\rhobold{{\ensuremath{\boldsymbol{\rho}}}}
\newcommand\taubold{{\ensuremath{\boldsymbol{\tau}}}}

\newcommand\Phibold{{\ensuremath{\boldsymbol{\Phi}}}}

\newcommand\Psibold{{\ensuremath{\boldsymbol{\Psi}}}}
\newcommand\Xibold{{\ensuremath{\boldsymbol{\Xi}}}}

\newcommand\zerobold{\ensuremath{\mathbf{0}}}
\newcommand\onebold{\ensuremath{\mathbf{1}}}

\usepackage{tikz}
\usepackage{pgfplots}
\usepackage{pgfplotstable, filecontents, booktabs}
\pgfplotsset{compat=1.9}

\usetikzlibrary{pgfplots.groupplots}
\usepgfplotslibrary{fillbetween}
\usetikzlibrary{calc,fit,matrix,arrows,automata,positioning,shapes}
\usetikzlibrary{arrows.meta}

\pgfplotsset{select coords between index/.style 2 args={
    x filter/.code={
        \ifnum\coordindex<#1\fi
        \ifnum\coordindex>#2\fi
    }
}}

\tikzset{
 invisible/.style={opacity=0},
 visible on/.style={alt={#1{}{invisible}}},
 alt/.code args={<#1>#2#3}{%
   \alt<#1>{\pgfkeysalso{#2}}{\pgfkeysalso{#3}}
 },
}

\usepackage{setspace}

\newbool{fastcompile}
\setbool{fastcompile}{false}

\newcommand{\pgrad}{\nabla}

\begin{document}
\title{Accelerated solutions of convection-dominated partial differential equations
 using implicit feature tracking and empirical quadrature}

\author[rvt1]{Marzieh Alireza Mirhoseini\fnref{fn1}}
\ead{malireza@nd.edu}

\author[rvt1]{Matthew J. Zahr\fnref{fn2}\corref{cor1}}
\ead{mzahr@nd.edu}

\address[rvt1]{Department of Aerospace and Mechanical Engineering, University
               of Notre Dame, Notre Dame, IN 46556, United States}
\cortext[cor1]{Corresponding author}

\fntext[fn1]{Graduate Student, Department of Aerospace and Mechanical
             Engineering, University of Notre Dame}
\fntext[fn2]{Assistant Professor, Department of Aerospace and Mechanical
             Engineering, University of Notre Dame}


\begin{abstract}	
This work introduces an empirical quadrature-based hyperreduction procedure and
greedy training algorithm to effectively reduce the computational cost of solving
convection-dominated problems with limited training. The proposed approach
circumvents the slowly decaying $n$-width limitation of linear model reduction
techniques applied to convection-dominated problems by using a nonlinear
approximation manifold systematically defined by composing a low-dimensional
affine space with bijections of the underlying domain. The reduced-order model is
defined as the solution of a residual minimization problem over the nonlinear
manifold. An online-efficient method is obtained by using empirical quadrature
to approximate the optimality system such that it can be solved with
mesh-independent operations. The proposed reduced-order model
is trained using a greedy procedure to systematically sample the
parameter domain. The effectiveness of the proposed approach is
demonstrated on two shock-dominated computational fluid dynamics benchmarks.
\end{abstract}

\maketitle

\textbf{Data sharing}: Data sharing is not applicable to this article as no new data were created or analyzed in this study.

\section{Introduction}
\label{sec:intro}
Projection-based model reduction methods search for an approximate
solution to a partial differential equation (PDE) in a low-dimensional
(affine) subspace. The reduced subspace is constructed in a time-intensive
offline phase, where the PDE solution is queried at a number of sampled
parameter configurations and compressed to form an orthonormal basis. From
the reduced solution space, the reduced-order model (ROM) is defined via a
projection of the governing equations and serves as a surrogate for the expensive
PDE in the online stage. This approach has been shown to accelerate
PDE solutions by orders of magnitude with minimal loss in error, particularly
for diffusion-dominated problems \cite{veroy2003posteriori,grepl2005posteriori}.
However, the approach is fundamentally limited---requires extensive training
and large reduced subspaces, and has poor generalizability beyond training---for
convection-dominated problems owing to the slowly decaying Kolmogorov $n$-width
of the solution manifold \cite{ohlberger_reduced_2016}.

This fundamental limitation of projection-based model reduction has been
widely documented and a number of solutions have been proposed. One such
solution was the development of localized ROMs where the affine approximation
is replaced by a collection of affine approximations
\cite{haasdonk_training_2011,dihlmann_model_2011,amsallem_nonlinear_2012}.
While these methods can offset the Kolmogorov $n$-width issue by using
a collection of smaller bases, their accuracy and prediction capability
is bounded by that of the union of all localized affine spaces
\cite{ohlberger_reduced_2016} and generalizability remains an issue.
Another class of methods that aim to overcome the Kolmogorov $n$-width
limitation of traditional model reduction are based on online adaptivity
of the reduced subspace, which breaks the offline-online decomposition.
Such approaches have been used to accelerate optimization problems
\cite{arian_trust-region_2000, agarwal_trust-region_2013, yue_accelerating_2013, zahr_progressive_2015}, correct reduced models with the PDE solution near
discontinuities \cite{constantine_reduced_2012,lucia_reduced_2003}, or fully
adapt the reduced space online using rank-one updates \cite{peherstorfer_model_2020},
partial PDE solutions \cite{huang2023predictive,zucatti2023adaptive},
and $h$-refinement \cite{carlberg_adaptive_2015}. While these methods
have been shown to dramatically improve the predictive capability
of ROMs, their speedup potential is much more limited than methods
that leverage an offline-online decomposition. Finally, the methods
that aim to directly resolve the issues posed by the slowly
decaying Kolmogorov $n$-width of convection-dominated problems leverage
nonlinear trial manifolds. These approximation manifolds are usually
constructed as a nonlinear parametrization of an affine space (e.g.,
a transformation of the underlying spatial domain)
\cite{ohlberger_nonlinear_2013, welper_interpolation_2017, reiss_shifted_2018, rim_transport_2018, rim_displacement_2018, cagniart_model_2019, nair_transported_2019, black_projection-based_2020, taddei_registration_2020, welper_transformed_2020, torlo_model_2020, mojgani_physics-aware_2020, rim2023manifold, bansal2021model, ferrero2022registration, taddei2021registration,mirhoseini2023model} or directly using machine learning techniques
\cite{kashima_nonlinear_2016, hartman_deep_2017, omata_novel_2019,
      lee_model_2020, xu_multi-level_2020, maulik_reduced-order_2020,
      kim2020efficient}.
Such methods have been shown to substantially improve the predictive capability of
ROMs; however, they present their own difficulties, e.g., large training costs, code
intrusivity, and general and efficient construction of invertible domain mappings.

In our previous work \cite{mirhoseini2023model}, we introduced a model reduction method,
\textit{implicit feature tracking}, that constructs a nonlinear approximation manifold
by composing a low-dimensional affine space with a space of bijections of the underlying
domain. The reduced-order model with implicit feature tracking (ROM-IFT) is defined as
a residual minimization problem over the reduced nonlinear manifold that simultaneously
determines the reduced coordinates in the affine space and the domain mapping that
minimizes the residual of the unreduced PDE discretization, which is analogous to
standard minimum-residual reduced-order models,
except our method expands the optimization space to include admissible domain mappings.
The nonlinear trial manifold is constructed by using the proposed residual minimization
formulation to determine domain mappings that cause parametrized features (e.g.,
discontinuities, extrema, vortical structures) to align in a reference domain for a
set of training parameters. Because the feature is stationary in the reference domain,
i.e., the convective nature of solution removed, the snapshots are effectively compressed
to define an affine subspace. The method was shown to accurately approximate
shock-dominated flows with limited training.

Implicit feature tracking and any nonlinear manifold approaches based on composing
a low-dimensional affine space with domain mappings are fundamentally based on aligning
features in the reduced basis with the corresponding feature in the solution being
approximated.
This is precisely the idea behind implicit shock tracking methods
\cite{corrigan_moving_2019,zahr_implicit_2020}
that aim to align element faces---locations of discontinuities in a piecewise
polynomial basis---with flow discontinuities. In this sense, implicit feature tracking
is the extension of implicit shock tracking to the model reduction setting,
which is further supported by the fact both methods have a similar optimization
formulation.
Unlike implicit shock tracking that requires the underlying solution
representation to support discontinuities between elements (e.g., finite volume or
discontinuous Galerkin methods), implicit feature tracking only requires approximations
of relevant features to exist in the underlying reduced basis.
Although implicit feature tracking is a general model reduction method
that can be applied to any PDE discretization, it is natural to apply to implicit
shock tracking discretizations due to the similarities between the formulations.

In this work, we adapt an empirical quadrature-based hyperreduction
procedure \cite{yano_lp_2019} to the implicit feature tracking
setting to remove a key bottleneck associated with assembling nonlinear
terms in the optimality system. While hyperreduction can be accomplished via
any of the techniques available in the literature
\cite{barrault2004empirical,astrid2008missing,chaturantabut2010nonlinear},
we choose an approach based on empirical quadrature
\cite{farhat2014dimensional,yano_lp_2019}
because they are natural and convenient to implement in a finite element
setting. Additionally, a greedy algorithm \cite{maday2002priori,veroy2003posteriori}
is developed to systematically sample the parameter domain to train the proposed
reduced-order model. For efficiency, we use hyperreduced solutions from previous
greedy iterations to train the empirical quadrature weights at the current iteration,
which keeps training costs low.

The remainder of the paper is organized as follows.
Section~\ref{sec:govern} introduces the governing system of
parametrized partial differential equations, its transformation
to a fixed reference domain, and the corresponding discretization.
Section~\ref{sec:rom} recalls the implicit feature tracking method
of \cite{mirhoseini2023model} and introduces an empirical quadrature
hyperreduction procedure to accelerate assembly of nonlinear terms, a
Levenberg-Marquardt solver for the hyperreduced optimization problems,
and a greedy training algorithm. Section~\ref{sec:numexp} demonstrates
the merits of the proposed hyperreduction framework for two shock-dominated
benchmarks and Section~\ref{sec:concl} concludes the paper.

\section{Governing equations and discretization}
\label{sec:govern}
In this section, we introduce the governing partial differential equations
(Section~\ref{sec:govern:claw}), its transformation to a reference domain
so that a domain mapping appears explicitly in the governing equations
(Section~\ref{sec:govern:tclaw}), and a discretization of the transformed
conservation law to yield the fully discrete governing equations
(Section~\ref{sec:govern:disc}).

\subsection{Parametrized system of conservation laws}
\label{sec:govern:claw}
Consider a $\mubold$-parametrized system of $m$ conservation laws in
$d$-dimensions of the form
\begin{equation} \label{eqn:claw-phys}
 \nabla\cdot f(u,\nabla u;\mubold) = s(u,\nabla u;\mubold)
 \quad\text{in}\quad\Omega,
\end{equation}
where $\mubold\in\Dcal$ is the parameter and $\Dcal$ the parameter domain,
$u(x;\mubold)\in\Rbb^m$ is the solution (assumed unique) at a point
$x\in\Omega$ implicitly defined by (\ref{eqn:claw-phys}),
$\nabla\coloneqq[\partial_{x_1},\cdots,\partial_{x_d}]$
is the gradient operator on the domain $\Omega\subset\Rbb^d$ such that
$\nabla z\coloneqq [\partial_{x_1}z \cdots \partial_{x_d}z]$,
$\func{f}{\Rbb^m\times\Rbb^{m\times d}\times\Dcal}{\Rbb^{m\times d}}$
is the flux function, and
$\func{s}{\Rbb^m\times\Rbb^{m\times d}\times\Dcal}{\Rbb^m}$
is the source term. The boundary of the domain is $\partial\Omega$
with outward unit normal $\func{n}{\partial\Omega}{\Rbb^d}$. The
formulation of the conservation law in (\ref{eqn:claw-phys}) is
sufficiently general to encapsulate steady second-order partial differential
equations (PDEs) in  a $d$-dimensional spatial domain or time-dependent
PDEs in a $(d-1)$-dimensional domain, i.e., a $d$-dimensional space-time
domain.

\subsection{Transformed system of conservation laws on a fixed reference domain}
\label{sec:govern:tclaw}
As the proposed model reduction method is fundamentally based on
deforming the domain, it is convenient to recast the PDE
(\ref{eqn:claw-phys}) on a fixed \textit{reference} domain,
$\Omega_0\coloneqq\bar\Gcal^{-1}(\Omega)\subset\Rbb^d$,
following the approach in \cite{persson2009discontinuous}, where
$\func{\bar\Gcal}{\Rbb^d}{\Rbb^d}$ is a smooth invertible mapping.
Let $\Gbb$ be any collection of bijections from the reference
domain $\Omega_0$ to the physical domain $\Omega$ (Figure~\ref{fig:mapping}).
Then, for any $\Gcal\in\Gbb$, (\ref{eqn:claw-phys}) can be
written as a PDE on the reference domain as
\begin{equation} \label{eqn:claw-ref}
 \nabla_0\cdot F(U,\nabla_0 U;\Gcal,\mubold)=S(U,\nabla_0 U;\Gcal,\mubold)
 \quad\text{in}\quad\Omega_0,
\end{equation}
where $U(X;\Gcal,\mubold)$ is the solution of the transformed PDE
at a point $X\in\Omega_0$,
$\nabla_0\coloneqq[\partial_{X_1},\cdots,\partial_{X_d}]$
is the gradient operator on the domain $\Omega_0\subset\Rbb^d$ such that
$\nabla_0 z\coloneqq [\partial_{X_1}z \cdots \partial_{X_d}z]$,
$\func{F}{\Rbb^m\times\Rbb^{m\times d}\times\Gbb\times\Dcal}{\Rbb^{m \times d}}$ is the transformed flux function, and
$\func{S}{\Rbb^m\times\Rbb^{m\times d}\times\Gbb\times\Dcal}{\Rbb^m}$
is the transformed source term. The boundary of the domain is $\partial\Omega_0$
with outward unit normal $\func{N}{\partial\Omega_0}{\Rbb^d}$.
\begin{figure}
\centering
\input{_py/mapping0.tikz}
\caption{Schematic of domain mapping for $\Gcal\in\Gbb$}
\label{fig:mapping}
\end{figure}
The reference and physical solutions and their gradients are related as
\begin{equation}\label{eqn:p2r-stvc}
 U(X;\Gcal,\mubold) =
 u(\Gcal(X);\mubold), \qquad
 \nabla_0 U(X;\Gcal,\mubold) =
 \pgrad u(\Gcal(X);\mubold)G_\Gcal(X),
\end{equation}
where $G_\Gcal\coloneqq\nabla_0\Gcal$ is the mapping
Jacobian and $g_\Gcal\coloneqq\det G_\Gcal$ is its
determinant. The reference and physical flux functions and source
terms are related as
\begin{equation}\label{eqn:p2r-flux}
 F(w,\nabla_0 w;\Gcal,\mubold) =
 g_\Gcal
 f(w, \nabla_0 w \cdot G_\Gcal^{-1};\mubold)
 G_\Gcal^{-T}, \qquad
 S(w,\nabla_0 w;\Gcal,\mubold) =  g_\Gcal s(w, \nabla_0 w \cdot G_\Gcal^{-1};\mubold)
\end{equation}
where $\func{w}{\Omega_0}{\Rbb^m}$ is any $m$-valued function over the reference
domain; see \cite{persson2009discontinuous} for details of derivation.
The conservation laws in
(\ref{eqn:claw-phys}) and (\ref{eqn:claw-ref}) are equivalent
provided $\Gcal\in\Gbb$ is a diffeomorphism,
that is, if $u$ is the solution of (\ref{eqn:claw-phys}) and $U$
satisfies (\ref{eqn:p2r-stvc}), then $U$ is the solution of
(\ref{eqn:claw-ref}), and vice versa. Owing to this equivalence,
the PDE solution $u(\,\cdot\,;\mubold)$ is independent of
the (diffeomorphic) domain mapping $\Gcal\in\Gbb$.

\subsection{Discretization of the transformed conservation law}
\label{sec:govern:disc}
The transformed PDE (\ref{eqn:claw-ref}) is discretized using
a suitable method and the domain mapping is parametrized in
terms of a finite-dimensional vector of coefficients
$\xbm \in \Rbb^{N_\xbm}$, i.e.,
$\Gcal = \Gcal(\,\cdot\,; \xbm)$,
to yield a nonlinear system of algebraic equations
\begin{equation} \label{eqn:hdm}
 \Rbm(\Ubm^\star;\xbm,\mubold) = \zerobold,
\end{equation}
where $\Ubm^\star(\xbm,\mubold)\in\Rbb^{N_\ubm}$ is the (assumed unique) solution
of the discretized PDE for $\xbm\in\Rbb^{N_\xbm}$, $\mubold\in\Dcal$,
and $\func{\Rbm}{\Rbb^{N_\ubm}\times\Rbb^{N_\xbm}\times\Dcal}{\Rbb^{N_\ubm}}$ is the
nonlinear residual function defined by the spatial discretization of
(\ref{eqn:claw-ref}). We assume the discretization has an assembled
structure, that is, the residual function can be written as
\begin{equation} \label{eqn:hdm_elem}
 \Rbm(\Ubm;\xbm,\mubold) =
 \sum_{e=1}^{N_\mathtt{e}} \Pbm_e\Rbm_e(\Ubm_e, \Ubm_e', \xbm_e, \mubold)
\end{equation}
for any $\Ubm\in\Rbb^{N_\ubm}$, where
$\Rbm_e : \Rbb^{n_\ubm} \times \Rbb^{n_\ubm'} \times \Rbb^{n_\xbm} \times \Dcal \rightarrow \Rbb^{n_\ubm}$
with
$\Rbm_e : (\Ubm_e, \Ubm_e', \xbm_e, \mubold) \mapsto
 \Rbm_e(\Ubm_e,\Ubm_e',\xbm_e,\mubold)$
is the residual contribution of element $e$,
$\Ubm_e \in \Rbb^{n_\ubm}$ are the DoFs associated with element $e$,
$\Ubm_e' \in \Rbb^{n_\ubm'}$ are the DoFs associated with elements that neighbor element $e$
(if applicable),
$\xbm_e \in \Rbb^{n_\xbm}$ are the deformation DoFs associated with element $e$,
$\Pbm_e \in \Rbb^{N_\ubm\times n_\ubm}$ is the assembly operator that maps element DoFs for
element $e$ to global DoFs, and
$\Qbm_e \in \Rbb^{N_\xbm\times n_\xbm}$ is the assembly operator that maps
element deformation DoFs of element $e$ to global deformation DoFs.
Furthermore, the element and global DoFs are related via assembly operators
$\Ubm_e = \Pbm_e^T\Ubm$, $\Ubm_e' = (\Pbm_e')^T\Ubm$,
$\xbm_e = \Qbm_e^T\xbm$, where $\Pbm_e' \in \Rbb^{N_\ubm\times n_\ubm'}$ is an
assembly operator that maps element DoFs from neighbors of element $e$
to the corresponding global DoFs. In (\ref{eqn:hdm_elem}) we localized
the domain mapping to an element as
$\left.\Gcal(\,\cdot\,;\xbm)\right|_{\Omega_{0,e}} = \Gcal_e(\,\cdot\,;\xbm_e)$,
where $\Gcal_e : \Omega_{0,e} \times \Rbb^{n_\xbm} \mapsto \Omega$ is an
elemental mapping.

Finite element and element-centered finite volume discretizations
possess the assembled structure (\ref{eqn:hdm_elem}) \cite{wen2023globally}.
We assume $\Pbm_e^T\Pbm_E = \delta_{eE} \Ibm$, where $\Pbm_E$ is the
assembly operator for another element $E$, which holds for for any discretization
where distinct elements do not share degrees of freedom.
In this work, we use a nodal discontinuous Galerkin method \cite{hesthaven2007nodal}
as the spatial discretization and define the domain mapping $\Gcal$ as a
$C^0(\Omega_0)$ piecewise polynomial of degree $q$ over the mesh used for
spatial discretization, parametrized by the nodal coordinates of the mesh.
We let $\Xbm\in\Rbb^{N_\xbm}$ denote the nodes of the reference mesh,
making $\Gcal(X; \Xbm) = X$ (the identity map).
Details on construction of the domain mapping can be found
in \cite{mirhoseini2023model}.
Throughout the document, (\ref{eqn:hdm}) will
be referred to as the high-dimensional model (HDM).

\section{Implicit feature tracking accelerated by empirical quadrature}
\label{sec:rom}
In this section we provide an overview of the implicit feature tracking model
reduction framework \cite{mirhoseini2023model} (Section~\ref{sec:rom:ift}),
introduce an empirical quadrature procedure \cite{yano_lp_2019} to accelerate
assembly of the nonlinear terms (Section~\ref{sec:rom:eqp}), and
propose a greedy procedure to construct these empirical quadrature-based
hyperreduced models with implicit feature tracking (EQP-IFT)
(Section~\ref{sec:rom:greedy}).

\subsection{Model reduction with implicit feature tracking}
\label{sec:rom:ift}
Central to the implicit feature tracking method introduced in
\cite{mirhoseini2023model} is the following low-dimensional
linear approximation of the reference state at a specific domain mapping
that causes features to align in the reference domain
\begin{equation}
 \Ubm^\star(\hat\xbm_{\Phibold,\Psibold}^\star(\mubold), \mubold) \approx
 \hat\Ubm_{\Phibold,\Psibold}^\star(\mubold) \coloneqq \Phibold\hat\wbm_{\Phibold,\Psibold}^\star(\mubold),
 \qquad
 \hat\xbm_{\Phibold,\Psibold}^\star(\mubold)\coloneqq\Xbm+\Psibold\hat\taubold_{\Phibold,\Psibold}^\star(\mubold),
\end{equation}
where $\hat\xbm_{\Phibold,\Psibold}^\star : \Dcal \rightarrow \Rbb^{N_\xbm}$ is a
low-dimensional parametrization of the domain mapping coefficients ($\xbm$)
with reduced basis $\Psibold\in\Rbb^{N_\ubm\times k_\xbm}$ and coefficients
$\hat\taubold_{\Phibold,\Psibold}^\star : \Dcal \rightarrow \Rbb^{k_\xbm}$
that lead to feature alignment ($k_\xbm \ll N_\xbm$),
$\hat\Ubm_{\Phibold,\Psibold}^\star : \Dcal \rightarrow \Rbb^{N_\ubm}$ is the low-dimensional
approximation of the reference state with reduced basis $\Phibold\in\Rbb^{N_\ubm\times k_\ubm}$
and coefficients $\hat\wbm_{\Phibold,\Psibold}^\star : \Dcal \rightarrow \Rbb^{k_\ubm}$
($k_\ubm \ll N_\ubm$). This linear approximation in the reference domain can be interpreted
as a nonlinear approximation in the physical domain. Construction of the reduced bases
$\Phibold$, $\Psibold$ will be deferred to Section~\ref{sec:rom:greedy}.

The reduced coefficients are determined simultaneously as the solution
of the optimization problem
\begin{equation} \label{eqn:ift-opt}
 (\hat\wbm_{\Phibold,\Psibold}^\star(\mubold), \hat\taubold_{\Phibold,\Psibold}^\star(\mubold)) \coloneqq
 \argmin_{(\wbm,\taubold)\in\Rbb^{k_\ubm}\times\Rbb^{k_\xbm}} \hat{J}_{\Phibold,\Psibold}(\wbm,\taubold;\mubold)
\end{equation}
where the objective function $\hat{J}_{\Phibold,\Psibold} : \Rbb^{k_\ubm} \times \Rbb^{k_\xbm} \times \Dcal$ is defined as the norm of the residual function
$\Fbm : \Rbb^{N_\ubm} \times \Rbb^{N_\xbm} \times \Dcal \rightarrow \Rbb^{N_\ubm+N_\mathtt{e}}$
that incorporates the HDM residual and a mesh distortion term
\begin{equation} \label{eqn:ift-obj}
 \hat{J}_{\Phibold,\Psibold} : (\wbm,\taubold;\mubold) \mapsto
 \frac{1}{2}\norm{\Fbm(\Phibold\wbm, \Xbm + \Psibold\taubold; \mubold)}_2^2,
 \qquad
 \Fbm : (\Ubm, \xbm; \mubold) \mapsto
 \begin{bmatrix}
  \Rbm(\Ubm; \xbm, \mubold) \\ \Nbm(\xbm)
 \end{bmatrix}.
\end{equation}
The residual $\Nbm : \Rbb^{N_\xbm} \rightarrow \Rbb^{N_\mathtt{e}}$ is
a measure of mesh distortion
\begin{equation} \label{eqn:dist_elem0}
 \Nbm : \xbm \mapsto \begin{bmatrix} \Nbm_1(\Qbm_1^T\xbm) \\ \vdots \\ \Nbm_{N_\mathtt{e}}(\Qbm_{N_\mathtt{e}}^T\xbm)\end{bmatrix}, \qquad
 \Nbm_e : \xbm_e \mapsto \kappa(\eta_e(\xbm_e) - \eta_e(\Xbm_e))
\end{equation}
and $\eta_e$ is the elemental distortion
\begin{equation} \label{eqn:dist_elem1}
 \eta_e : \xbm_e \mapsto \int_{\Omega_{0,e}} \left(\frac{\norm{G_e(X;\xbm_e)}_F^2}{\max\{g_e(X;\xbm_e), \epsilon\}^{2/d}}\right)^2 \, dX,
\end{equation}
where $G_e(\,\cdot\,;\xbm_e) \coloneqq \nabla_0 \Gcal_e(\,\cdot\,; \xbm_e)$ is the elemental
mapping deformation gradient, $g_e(\,\cdot\,;\xbm_e) \coloneqq \det G_e(\,\cdot\,;\xbm_e)$
is the elemental mapping Jacobian, $\Xbm_e = \Qbm_e^T\Xbm$, $\epsilon\in\Rbb_{>0}$ is a
tolerance preventing the denominator from reaching zero ($\epsilon=10^{-8}$ in this work),
and $\kappa > 0$ is a regularization parameter determined using the algorithm in
\cite{huang2022robust}.

The first-order optimality condition of (\ref{eqn:ift-opt}) leads to an
interpretation of implicit feature tracking as a projection-based reduced-order
model
\begin{equation}
  \hat\Sbm_{\Phibold,\Psibold}(\hat\wbm_{\Phibold,\Psibold}^\star(\mubold), \hat\taubold_{\Phibold,\Psibold}^\star(\mubold); \mubold) = \zerobold, \qquad
  \hat\Tbm_{\Phibold,\Psibold}(\hat\wbm_{\Phibold,\Psibold}^\star(\mubold), \hat\taubold_{\Phibold,\Psibold}^\star(\mubold); \mubold) = \zerobold,
\end{equation}
where the state $\hat\Sbm_{\Phibold,\Psibold} : \Rbb^{k_\ubm} \times \Rbb^{k_\xbm} \times \Dcal \rightarrow \Rbb^{k_\ubm}$ and deformation $\hat\Tbm_{\Phibold,\Psibold} : \Rbb^{k_\ubm} \times \Rbb^{k_\xbm} \times \Dcal \rightarrow \Rbb^{k_\xbm}$ residuals are defined as
\begin{equation} \label{eqn:ift-res}
\begin{aligned}
 \hat\Sbm_{\Phibold,\Psibold} : (\wbm, \taubold; \mubold) &\mapsto
   \Phibold^T
   \pder{\Rbm}{\Ubm}(\Phibold\wbm; \Xbm + \Psibold\taubold, \mubold)^T
   \Rbm(\Phibold\wbm; \Xbm + \Psibold\taubold, \mubold) \\
 \hat\Tbm_{\Phibold,\Psibold} : (\wbm, \taubold; \mubold) &\mapsto
   \Psibold^T
   \pder{\Rbm}{\xbm}(\Phibold\wbm; \Xbm + \Psibold\taubold; \mubold)^T
   \Rbm(\Phibold\wbm; \Xbm + \Psibold\taubold, \mubold)
  +
   \Psibold^T
   \pder{\Nbm}{\xbm}(\Xbm + \Psibold\taubold)^T
   \Nbm(\Xbm + \Psibold\taubold).
\end{aligned}
\end{equation}
Owing to the assembled form of the residual, the reduced residuals can be written
in assembled form as (derivation in \ref{app:elemental})
\begin{equation} \label{eqn:ift-res-elem1}
 \hat\Sbm_{\Phibold,\Psibold}(\wbm, \taubold; \mubold) = \sum_{e=1}^{N_\mathtt{e}} \hat\Sbm_{\Phibold,\Psibold}^{(e)}(\wbm,\taubold;\mubold), \qquad
 \hat\Tbm_{\Phibold,\Psibold}(\wbm, \taubold; \mubold) = \sum_{e=1}^{N_\mathtt{e}} \hat\Tbm_{\Phibold,\Psibold}^{(e)}(\wbm,\taubold;\mubold),
\end{equation}
where the elemental reduced residuals
$\hat\Sbm_{\Phibold,\Psibold}^{(e)} : \Rbb^{k_\ubm} \times \Rbb^{k_\xbm} \times \Dcal \rightarrow \Rbb^{k_\ubm}$ and
$\hat\Tbm_{\Phibold,\Psibold}^{(e)} : \Rbb^{k_\ubm} \times \Rbb^{k_\xbm} \times \Dcal \rightarrow \Rbb^{k_\xbm}$ are defined as
\begin{equation} \label{eqn:ift-res-elem2}
\begin{aligned}
 \hat\Sbm_{\Phibold,\Psibold}^{(e)} : (\wbm,\taubold;\mubold) &\mapsto
 \left(\Phibold_e^T\pder{\Rbm_e}{\Ubm_e}(\hat\Ubm_e,\hat\Ubm_e',\hat\xbm_e,\mubold)^T + (\Phibold_e')^T\pder{\Rbm_e}{\Ubm_e'}(\hat\Ubm_e,\hat\Ubm_e',\hat\xbm_e,\mubold)^T\right)
 \Rbm_e(\hat\Ubm_e,\hat\Ubm_e',\hat\xbm_e,\mubold) \\
 \hat\Tbm_{\Phibold,\Psibold}^{(e)} : (\wbm,\taubold;\mubold) &\mapsto
 \Psibold_e^T\left(\pder{\Rbm_e}{\xbm_e}(\hat\Ubm_e,\hat\Ubm_e',\hat\xbm_e,\mubold)^T\Rbm_e(\hat\Ubm_e,\hat\Ubm_e',\hat\xbm_e,\mubold) + \pder{\Nbm_e}{\xbm_e}(\hat\xbm_e)^T\Nbm_e(\hat\xbm_e)\right),
\end{aligned}
\end{equation}
where $\hat\Ubm_e=\Phibold_e\wbm$, $\hat\Ubm_e'=\Phibold_e'\wbm$,
$\hat\xbm_e=\Xbm_e+\Psibold_e\taubold$, $\Phibold_e = \Pbm_e^T\Phibold$,
$\Phibold_e' = (\Pbm_e')^T\Phibold$, and $\Psibold_e = \Qbm_e^T\Psibold$.

\begin{remark}
For $\kappa=0$, the optimization problem in (\ref{eqn:ift-opt}) corresponds to residual
minimization over the nonlinear manifold defined as the composition of the linear space
defined by $\Phibold$ with the space of domain mappings defined by $\Psibold$, which can
be directly contrasted to traditional minimum-residual approaches
\cite{legresley2006application} that freeze the domain and optimize only over the linear
subspace. Because the optimization problem is posed over a richer space, we expect to
obtain an approximation with a smaller residual than possible with standard reduced-order
models that use frozen modes. As a result, despite the focus of this work on
convection-dominated PDEs, this approach can be used for a larger class of
parametrized PDEs.
\end{remark}

\begin{remark}
The proposed optimization formulation avoids local loss of bijectivity ($g_e$ approaching
zero) because the objective function blows up for such mappings and would
be rejected by the solver line search (Section~\ref{sec:rom:eqp:solver}).
As a result, the maximum function in (\ref{eqn:dist_elem1}) does not slow the convergence
of gradient-based solvers because the kink in the maximum function only occurs at
sub-optimal configurations rejected by the line search ($g_e \approx \epsilon$).
\end{remark}

\begin{remark}
As described in \cite{mirhoseini2023model}, the domain mappings are parametrized by
coefficients in such a way that the domain geometry will be represented with high-order
accuracy regardless of the value of these coefficients (within some bounds). This means
the proposed approach that uses a reduced model for the domain mapping will not degrade
the accuracy of the boundary representation.
\end{remark}

\subsection{Empirical quadrature procedure for implicit feature tracking}
\label{sec:rom:eqp}
From the elemental form of the reduced residuals in (\ref{eqn:ift-res-elem1}),
it is clear that despite the reduced dimensionality ($k_\ubm \ll N_\ubm$, $k_\xbm \ll N_\xbm$),
assembly of the reduced residuals requires $\Ocal(N_\mathtt{e})$ operations. To circumvent this
bottleneck, we develop an empirical quadrature procedure \cite{yano_lp_2019} that
enables assembly of the nonlinear terms at a cost independent of $N_\mathtt{e}$.

\subsubsection{Formulation}
\label{sec:rom:eqp:form}
To accelerate the assembly of the reduced residual, we replace the reduced
residuals $\hat\Sbm_{\Phibold,\Psibold}$ and $\hat\Tbm_{\Phibold,\Psibold}$
with the \textit{weighted} residuals
$\tilde\Sbm_{\Phibold,\Psibold,\rhobold} : \Rbb^{k_\ubm} \times \Rbb^{k_\xbm} \times \Dcal \rightarrow \Rbb^{k_\ubm}$
and
$\tilde\Tbm_{\Phibold,\Psibold,\rhobold} : \Rbb^{k_\ubm} \times \Rbb^{k_\xbm} \times \Dcal \rightarrow \Rbb^{k_\xbm}$,
defined as
\begin{equation} \label{eqn:ift-eqp-res0}
 \tilde\Sbm_{\Phibold,\Psibold,\rhobold} : (\wbm,\taubold;\mubold) \mapsto
   \sum_{e=1}^{N_\mathtt{e}} \rho_e\hat\Sbm_{\Phibold,\Psibold}^{(e)}(\wbm,\taubold;\mubold),
 \qquad
 \tilde\Tbm_{\Phibold,\Psibold,\rhobold} : (\wbm,\taubold;\mubold) \mapsto
   \sum_{e=1}^{N_\mathtt{e}} \rho_e\hat\Tbm_{\Phibold,\Psibold}^{(e)}(\wbm,\taubold;\mubold).
\end{equation}
If $\rhobold$ is sparse, the weighted residuals can be assembled
efficiently; construction of $\rhobold$ deferred to Section~\ref{sec:rom:eqp:train}.
We define
$\tilde\wbm_{\Phibold,\Psibold,\rhobold}^\star : \Dcal \rightarrow \Rbb^{k_\ubm}$
and
$\tilde\taubold_{\Phibold,\Psibold,\rhobold}^\star : \Dcal \rightarrow \Rbb^{k_\xbm}$,
as the root of the weighted residuals
\begin{equation} \label{eqn:ift-eqp}
  \tilde\Sbm_{\Phibold,\Psibold,\rhobold}(\tilde\wbm_{\Phibold,\Psibold,\rhobold}^\star(\mubold), \tilde\taubold_{\Phibold,\Psibold,\rhobold}^\star(\mubold); \mubold) = \zerobold, \qquad
  \tilde\Tbm_{\Phibold,\Psibold,\rhobold}(\tilde\wbm_{\Phibold,\Psibold,\rhobold}^\star(\mubold), \tilde\taubold_{\Phibold,\Psibold,\rhobold}^\star(\mubold); \mubold) = \zerobold,
\end{equation}
and reconstruct an approximation to the HDM reference state $\Ubm^\star$ as
\begin{equation}
 \Ubm^\star(\tilde\xbm_{\Phibold,\Psibold,\rhobold}^\star(\mubold), \mubold) \approx
 \tilde\Ubm_{\Phibold,\Psibold,\rhobold}^\star(\mubold) \coloneqq \Phibold\tilde\wbm_{\Phibold,\Psibold,\rhobold}^\star(\mubold),
 \qquad
 \tilde\xbm_{\Phibold,\Psibold,\rhobold}^\star(\mubold)\coloneqq\Xbm+\Psibold\tilde\taubold_{\Phibold,\Psibold,\rhobold}^\star(\mubold).
\end{equation}
Owing to the property $\Pbm_e^T\Pbm_E = \delta_{eE} \Ibm$ of the assembly operator,
the EQP model in (\ref{eqn:ift-eqp-res0})-(\ref{eqn:ift-eqp}) is the first-order
optimality system of the following optimization problem (derivation in \ref{app:elemental})
\begin{equation} \label{eqn:ift-opt-eqp0}
 (\tilde\wbm_{\Phibold,\Psibold,\rhobold}^\star(\mubold), \tilde\taubold_{\Phibold,\Psibold,\rhobold}^\star(\mubold)) \coloneqq
 \argmin_{(\wbm,\taubold)\in\Rbb^{k_\ubm}\times\Rbb^{k_\xbm}} \tilde{J}_{\Phibold,\Psibold,\rhobold}(\wbm,\taubold;\mubold),
\end{equation}
where the objective function $\tilde{J}_{\Phibold,\Psibold} : \Rbb^{k_\ubm} \times \Rbb^{k_\xbm} \times \Dcal \rightarrow \Rbb$ mimics the original objective function
$\hat{J}_{\Phibold,\Psibold,\rhobold}$ in (\ref{eqn:ift-obj})
\begin{equation} \label{eqn:ift-opt-eqp1}
 \tilde{J}_{\Phibold,\Psibold,\rhobold} : (\wbm,\taubold;\mubold) \mapsto
 \frac{1}{2}\norm{\tilde\Fbm_{\Phibold,\Psibold,\rhobold}(\wbm,\taubold;\mubold)}_2^2
\end{equation}
with the hyperreduced residual function
$\tilde\Fbm_{\Phibold,\Psibold,\rhobold} : \Rbb^{k_\ubm}\times\Rbb^{k_\xbm}\times\Dcal\rightarrow \Rbb^{N_\ubm+N_\mathtt{e}}$, defined as
\begin{equation} \label{eqn:ift-eqp-res}
 \tilde\Fbm_{\Phibold,\Psibold,\rhobold} : (\wbm,\taubold;\mubold) \mapsto
 \begin{bmatrix}
  \tilde\Rbm_{\Phibold,\Psibold,\rhobold}(\wbm,\taubold;\mubold) \\
  \tilde\Nbm_{\Psibold,\rhobold}(\taubold)
 \end{bmatrix},
\end{equation}
replacing the original residual $\Fbm$, and the weighted residual
$\tilde\Rbm_{\Phibold,\Psibold,\rhobold} : \Rbb^{k_\ubm}\times\Rbb^{k_\xbm}\times\Dcal \rightarrow \Rbb^{N_\ubm}$
and mesh distortion
$\tilde\Nbm_{\Psibold,\rhobold} : \Rbb^{k_\xbm} \rightarrow \Rbb^{N_\mathtt{e}}$ are given by
\begin{equation}
\begin{aligned}
 \tilde\Rbm_{\Phibold,\Psibold,\rhobold} : (\wbm,\taubold;\mubold) &\mapsto
 \sum_{e=1}^{N_\mathtt{e}} \sqrt{\rho_e} \Pbm_e\Rbm_e(\Phibold_e\wbm,\Phibold_e'\wbm,\Xbm_e+\Psibold_e\taubold,\mubold) \\
 \tilde\Nbm_{\Psibold,\rhobold} : \taubold &\mapsto
 \begin{bmatrix}
   \sqrt{\rho_1} \Nbm_1(\Xbm_1 + \Psibold_1\taubold) \\
   \vdots \\
   \sqrt{\rho_{N_\mathtt{e}}} \Nbm_{N_\mathtt{e}}(\Xbm_{N_\mathtt{e}} + \Psibold_{N_\mathtt{e}}\taubold)
 \end{bmatrix}.
\end{aligned}
\end{equation}
Finally, owing to the $\Pbm_e^T\Pbm_E = \delta_{eE} \Ibm$ property of
the assembly operator, the hyperreduced objective function can be written
as a summation of weighted elemental contributions (derivation in \ref{app:elemental})
\begin{equation} \label{eqn:ift-opt-eqp2}
 \tilde{J}_{\Phibold,\Psibold,\rhobold}(\wbm,\taubold;\mubold) =
 \sum_{e=1}^{N_\mathtt{e}}
   \frac{\rho_e}{2}\norm{\Rbm_e(\Phibold_e\wbm,\Phibold_e'\wbm,\Xbm_e+\Psibold_e\taubold,\mubold)}_2^2
   +
 \sum_{e=1}^{N_\mathtt{e}}
   \frac{\rho_e}{2}\left|\Nbm_e(\Xbm_e+\Psibold_e\taubold)\right|^2,
\end{equation}
which clearly shows $\tilde{J}_{\Phibold,\Psibold}$ can be computed
efficiently provided $\rhobold$ is sparse.

\begin{remark} \label{rem:rom-eqp-ones}
The standard and weighted reduced residuals are related as
\begin{equation}
\hat\Sbm_{\Phibold,\Psibold}(\hat\wbm,\hat\taubold;\mubold) =
\tilde\Sbm_{\Phibold,\Psibold,\onebold}(\hat\wbm,\hat\taubold;\mubold),
\qquad
\hat\Tbm_{\Phibold,\Psibold}(\hat\wbm,\hat\taubold;\mubold) =
\tilde\Tbm_{\Phibold,\Psibold,\onebold}(\hat\wbm,\hat\taubold;\mubold),
\end{equation}
for any $\hat\wbm\in\Rbb^{k_\ubm}$, $\hat\taubold\in\Rbb^{k_\xbm}$, and $\mubold\in\Dcal$.
Therefore, the corresponding reduced coordinates are equal
\begin{equation}
 \hat\wbm_{\Phibold,\Psibold}^\star(\mubold) = \tilde\wbm_{\Phibold,\Psibold,\onebold}^\star(\mubold), \qquad
 \hat\taubold_{\Phibold,\Psibold}^\star(\mubold) = \tilde\taubold_{\Phibold,\Psibold,\onebold}^\star(\mubold).
\end{equation}
\end{remark}


\subsubsection{Training}
\label{sec:rom:eqp:train}
The success of EQP, in terms of both accuracy and efficiency, is inherently linked to 
the construction of a sparse weight vector that ensures the weighted residuals
(\ref{eqn:ift-eqp-res0}) are accurate representations of the reduced residuals
(\ref{eqn:ift-res}). Following the standard EQP constructions \cite{yano_lp_2019},
we choose the weights as the solution of an $\ell_1$ minimization problem to
promote sparsity with constraints included to ensure accuracy at training points.
To this end, let $\check\rhobold\in\Rbb^{N_\mathtt{e}}$ be a guess to the
sparse weights in (\ref{eqn:ift-eqp-res0}) ($\check\rhobold = \onebold$ if no guess
available), and define $\rhobold_{\Phibold,\Psibold,\check\rhobold}^\star$ as the
solution of the following linear program
\begin{equation} \label{eqn:eqp-lp}
 \rhobold_{\Phibold,\Psibold,\check\rhobold}^\star \coloneqq
 \argmin_{\rhobold \in \Ccal_{\Phibold,\Psibold,\check\rhobold}}~\sum_{e=1}^{N_\mathtt{e}} \rho_e,
 \qquad
 \Ccal_{\Phibold,\Psibold,\check\rhobold} \coloneqq
 \Ccal^\mathrm{nn} \cap
 \Ccal^\mathrm{dv} \cap
 \Ccal_{\Phibold,\Psibold,\check\rhobold}^\mathrm{sr} \cap
 \Ccal_{\Phibold,\Psibold,\check\rhobold}^\mathrm{dr},
\end{equation}
where the constraints are defined as
\begin{equation} \label{eqn:eqp-con}
\begin{aligned}
 \Ccal^\mathrm{nn} &\coloneqq
 \left\{\rhobold\in\Rbb^{N_\mathtt{e}} \suchthat \rho_e \geq 0,~e=1,\dots,N_\mathtt{e}\right\} \\
 \Ccal^\mathrm{dv} &\coloneqq
 \left\{\rhobold\in\Rbb^{N_\mathtt{e}} \suchthat \left| |\Omega_0| - \sum_{e=1}^{N_\mathtt{e}} \rho_e|\Omega_{0,e}|\right| \leq \delta |\Omega_0|\right\} \\
 \Ccal_{\Phibold,\Psibold,\check\rhobold}^\mathrm{sr} &\coloneqq
 \left\{ \rhobold\in\Rbb^{N_\mathtt{e}} \suchthat \norm{\hat\Sbm_{\Phibold,\Psibold}(\tilde\vbm_{\Phibold,\Psibold,\check\rhobold}^\star(\mubold);\mubold) - \tilde\Sbm_{\Phibold,\Psibold,\rhobold}(\tilde\vbm_{\Phibold,\Psibold,\check\rhobold}^\star(\mubold);\mubold)}_\infty \leq \delta|\Xibold|^{-1},~\mubold\in\Xibold \right\} \\
 \Ccal_{\Phibold,\Psibold,\check\rhobold}^\mathrm{dr} &\coloneqq
 \left\{ \rhobold\in\Rbb^{N_\mathtt{e}} \suchthat \norm{\hat\Tbm_{\Phibold,\Psibold}(\tilde\vbm_{\Phibold,\Psibold,\check\rhobold}^\star(\mubold);\mubold) - \tilde\Tbm_{\Phibold,\Psibold,\rhobold}(\tilde\vbm_{\Phibold,\Psibold,\check\rhobold}^\star(\mubold);\mubold)}_\infty \leq \delta|\Xibold|^{-1},~\mubold\in\Xibold \right\},
\end{aligned}
\end{equation}
and $\tilde\vbm_{\Phibold,\Psibold,\check\rhobold}^\star(\mubold) \coloneqq (\tilde\wbm_{\Phibold,\Psibold,\check\rhobold}(\mubold),\tilde\taubold_{\Phibold,\Psibold,\check\rhobold}(\mubold))$ is introduced for notational brevity, $\delta\in\Rbb_{>0}$ is a user-defined
tolerance, and $\Xibold\subset\Dcal$ is a collection of parameters used to train the
weights. 
The first constraint set ensures the weights are non-negative ($\mathrm{nn}$)
to maintain their interpretation as quadrature weights. The second constraint
set ensures the weights accurately integrate the domain volume ($\mathrm{dv}$).
The last two constraints ensure the weighted state ($\mathrm{sr}$) and domain
($\mathrm{dr}$) residuals accurately approximate the reduced residuals on
the training set $\Xibold$. The first two constraints were directly proposed
in \cite{yano_lp_2019}, while the latter two were adapted from the manifold accuracy
constraint in \cite{yano_lp_2019} to the implicit feature tracking setting.


\begin{remark}
In the case where $\check\rhobold = \onebold$, the manifold accuracy constraints
are evaluated at the ROM-IFT solution
\begin{equation}
\begin{aligned}
 \Ccal_{\Phibold,\Psibold}^\mathrm{sr} &=
 \left\{ \rhobold\in\Rbb^{N_\mathtt{e}} \suchthat \norm{\hat\Sbm_{\Phibold,\Psibold}(\hat\wbm_{\Phibold,\Psibold}^\star(\mubold),\hat\taubold_{\Phibold,\Psibold}^\star(\mubold);\mubold) - \tilde\Sbm_{\Phibold,\Psibold,\rhobold}(\hat\wbm_{\Phibold,\Psibold}^\star(\mubold),\hat\taubold_{\Phibold,\Psibold}^\star(\mubold);\mubold)}_\infty \leq \delta|\Xibold|^{-1},~\mubold\in\Xibold \right\} \\
 \Ccal_{\Phibold,\Psibold}^\mathrm{dr} &=
 \left\{ \rhobold\in\Rbb^{N_\mathtt{e}} \suchthat \norm{\hat\Tbm_{\Phibold,\Psibold}(\hat\wbm_{\Phibold,\Psibold}^\star(\mubold),\hat\taubold_{\Phibold,\Psibold}^\star(\mubold);\mubold) - \tilde\Tbm_{\Phibold,\Psibold,\rhobold}(\hat\wbm_{\Phibold,\Psibold}^\star(\mubold),\hat\taubold_{\Phibold,\Psibold}^\star(\mubold);\mubold)}_\infty \leq \delta|\Xibold|^{-1},~\mubold\in\Xibold \right\}
\end{aligned}
\end{equation}
because of the equivalence of the reduced and weighted solutions
(Remark~\ref{rem:rom-eqp-ones}). This can be further simplified to
\begin{equation}
\begin{aligned}
 \Ccal_{\Phibold,\Psibold}^\mathrm{sr} &=
 \left\{ \rhobold\in\Rbb^{N_\mathtt{e}} \suchthat \norm{\tilde\Sbm_{\Phibold,\Psibold,\rhobold}(\hat\wbm_{\Phibold,\Psibold}^\star(\mubold),\hat\taubold_{\Phibold,\Psibold}^\star(\mubold);\mubold)}_\infty \leq \delta|\Xibold|^{-1},~\mubold\in\Xibold \right\} \\
 \Ccal_{\Phibold,\Psibold}^\mathrm{dr} &=
 \left\{ \rhobold\in\Rbb^{N_\mathtt{e}} \suchthat \norm{\tilde\Tbm_{\Phibold,\Psibold,\rhobold}(\hat\wbm_{\Phibold,\Psibold}^\star(\mubold),\hat\taubold_{\Phibold,\Psibold}^\star(\mubold);\mubold)}_\infty \leq \delta|\Xibold|^{-1},~\mubold\in\Xibold \right\},
\end{aligned}
\end{equation}
from the definition of $\hat\wbm_{\Phibold,\Psibold}^\star$,
$\hat\taubold_{\Phibold,\Psibold}^\star$.
\end{remark}

\begin{remark}
Regardless of $\check\rhobold$, the feasible set of the linear program in
(\ref{eqn:eqp-lp}) is non-empty as $\rhobold = \onebold$ is feasible because
of the equivalence of the reduced and weighted residuals (Remark~\ref{rem:rom-eqp-ones}).
Furthermore, the optimization problem in (\ref{eqn:eqp-lp}) is a linear
program because each weighted term is linear in $\rhobold$ and the infinity-norm
bounds can be recast as a collection of inequality constraints (upper and lower bounds)
on its argument.
\end{remark}

\begin{remark}
As noted in \cite{wen2023globally}, the $\Pbm_e^T\Pbm_E=\delta_{eE}\Ibm$ property
and $\check\rhobold=\onebold$ implies each element residual is zero, which means
$\tilde\Sbm_{\Phibold,\Psibold}^{(e)}(\hat\wbm_{\Phibold,\Psibold}^\star(\mubold),\hat\taubold_{\Phibold,\Psibold}^\star(\mubold);\mubold) = \zerobold$
and
$\hat\Tbm_{\Phibold,\Psibold}^{(e)}(\hat\wbm_{\Phibold,\Psibold}^\star(\mubold),\hat\taubold_{\Phibold,\Psibold}^\star(\mubold);\mubold) = \zerobold$
for any $e = 1,\dots,N_\mathtt{e}$ and therefore the manifold accuracy constraints
are satisfied for any $\rhobold\in\Rbb^{N_\mathtt{e}}$. To ensure the manifold accuracy
constraints are meaningful, we follow the approach in \cite{wen2023globally} and split
the manifold accuracy constraints into volume and face contributions.
\end{remark}

\begin{remark}
Our empirical quadrature approach is the direct generalization of the EQP method in
\cite{yano_lp_2019} to the implicit feature tracking setting. The Galerkin
reduced-order model is replaced by the minimum-residual feature tracking formulation
(\ref{eqn:ift-opt}) and EQP is applied to the first-order optimality system. In the
training, the standard manifold accuracy constraint is replaced with the corresponding
equations from the ROM-IFT optimality system.
\end{remark}

\begin{remark}
There are several approaches to train a hyperreduced model based on empirical
quadrature including the $\ell_1$ minimization approach proposed in \cite{yano_lp_2019}
and a non-negative least squares approach proposed in \cite{farhat2015structure}.
We use $\ell_1$ minimization because multiple accuracy constraints can easily be
incorporated with individual tolerances; however, a non-negative least squares
training formulation could be built as well with potential advantages for
large problems \cite{chapman2017accelerated}.
\end{remark}

\subsubsection{Levenberg-Marquardt solver}
\label{sec:rom:eqp:solver}
To solve the hyperreduced optimization problem in
(\ref{eqn:ift-eqp}), we use the Levenberg-Marquardt solver in \cite{mirhoseini2023model},
which produces a sequence of iterates $\{\wbm_n\}_{n=0}^\infty\subset\Rbb^{k_\ubm}$
and $\{\taubold_n\}\subset\Rbb^{k_\xbm}$ as
\begin{equation}
 \wbm_{n+1} = \wbm_n + \alpha_{n+1} \Delta\wbm_{n+1},
 \quad
 \taubold_{n+1} = \taubold_n + \alpha_{n+1} \Delta\taubold_{n+1},
\end{equation}
such that $\lim_{n\rightarrow \infty} \wbm_n = \tilde\wbm_{\Phibold,\Psibold,\rhobold}^\star(\mubold)$ and $\lim_{n\rightarrow \infty} \taubold_n = \tilde\taubold_{\Phibold,\Psibold,\rhobold}^\star(\mubold)$, where $\wbm_0$ and $\taubold_0$ initialize the iteration
(Remark~\ref{rem:solver}) and $\alpha_{n+1}\in\Rbb_{>0}$ is the step length chosen to
ensure sufficient decrease of $\tilde{J}_{\Phibold,\Psibold,\rhobold}$ based on the Wolfe
conditions \cite{nocedal_numerical_2006}. The steps $\Delta\wbm_n$ and $\Delta\taubold_n$
are defined as the solution of the linear least-squares problem
\begin{equation} \label{eqn:levmarq-step}
 (\Delta\wbm_{n+1},\Delta\taubold_{n+1}) \coloneqq
 \argmin_{(\Delta\wbm,\Delta\taubold)\in\Rbb^{k_\ubm}\times\Rbb^{k_\xbm}}
 \norm{\begin{bmatrix} \tilde\Fbm_n \\ \zerobold \end{bmatrix}
       +
       \begin{bmatrix} \tilde\Jbm_{\wbm,n} & \tilde\Jbm_{\taubold,n} \\ \zerobold & \sqrt{\lambda_n} \Ibm_{k_\xbm}\end{bmatrix}
       \begin{bmatrix}
        \Delta\wbm \\ \Delta\taubold
       \end{bmatrix}}_2,
\end{equation}
where
\begin{equation}
 \tilde\Fbm_n \coloneqq \tilde\Fbm_{\Phibold,\Psibold,\rhobold}(\wbm_n,\taubold_n;\mubold),
 \qquad
 \tilde\Jbm_{\wbm,n} \coloneqq \pder{\tilde\Fbm_{\Phibold,\Psibold,\rhobold}}{\wbm}(\wbm_n,\taubold_n;\mubold),
 \qquad
 \tilde\Jbm_{\taubold,n} \coloneqq \pder{\tilde\Fbm_{\Phibold,\Psibold,\rhobold}}{\taubold}(\wbm_n,\taubold_n;\mubold),
\end{equation}
and $\lambda_n\in\Rbb_{>0}$ is a regularization parameter chosen adaptively based on the
strategy in \cite{huang2022robust} to reduce dependence on user-defined parameters.
The analytical solution of the linear least-squares problem (\ref{eqn:levmarq-step})
leads to the following expression for the step
\begin{equation} \label{eqn:lmstep}
 \tilde\Jbm_n^T\tilde\Jbm_n \Delta\zbm_{n+1} = -\begin{bmatrix} \tilde\Jbm_{\wbm,n}^T\tilde\Fbm_n \\ \tilde\Jbm_{\taubold,n}^T\tilde\Fbm_n \end{bmatrix},
\end{equation}
where
\begin{equation}
 \tilde\Jbm_n = \begin{bmatrix} \tilde\Jbm_{\wbm,n} & \tilde\Jbm_{\taubold,n} \\ \zerobold & \sqrt{\lambda_n} \Ibm_{k_\xbm} \end{bmatrix},
 \qquad
 \Delta\zbm_{n+1} = \begin{bmatrix} \Delta\wbm_{n+1} \\ \Delta\taubold_{n+1} \end{bmatrix},
\end{equation}
which is a (small) $(k_\ubm+k_\xbm) \times (k_\ubm+k_\xbm)$ linear system for
$\Delta\zbm_{n+1}$. The right-hand side of the linear system in
(\ref{eqn:lmstep}) is equivalent to the hyperreduced residuals
(derivation in \ref{app:elemental})
\begin{equation}
 \tilde\Jbm_{\wbm,n}^T\tilde\Fbm_n = \tilde\Sbm_{\Phibold,\Psibold}(\wbm_n,\taubold_n;\mubold,\rhobold), \qquad
 \tilde\Jbm_{\taubold,n}^T\tilde\Fbm_n = \tilde\Tbm_{\Phibold,\Psibold}(\wbm_n,\taubold_n;\mubold,\rhobold).
\end{equation}
Furthermore, the following three terms are required to assemble the
linear system matrix (derivation in \ref{app:elemental})
\begin{equation} \label{eqn:levmarq-elem}
\begin{aligned}
 \tilde\Jbm_{\wbm,n}^T\tilde\Jbm_{\wbm,n} &=
   \sum_{e=1}^{N_\mathtt{e}} \rho_e \left(\Phibold_e^T\pder{\Rbm_e}{\Ubm_e}^T+(\Phibold_e')^T\pder{\Rbm_e}{\Ubm_e'}^T\right) \left(\pder{\Rbm_e}{\Ubm_e}\Phibold_e+\pder{\Rbm_e}{\Ubm_e'}\Phibold_e'\right) \\
 \tilde\Jbm_{\wbm,n}^T\tilde\Jbm_{\taubold,n} &=
   \sum_{e=1}^{N_\mathtt{e}} \rho_e \left(\Phibold_e^T\pder{\Rbm_e}{\Ubm_e}^T+(\Phibold_e')^T\pder{\Rbm_e}{\Ubm_e'}^T\right) \pder{\Rbm_e}{\xbm_e}\Psibold_e \\
 \tilde\Jbm_{\taubold,n}^T\tilde\Jbm_{\taubold,n} &=
   \sum_{e=1}^{N_\mathtt{e}} \rho_e \Psibold_e^T\left(\pder{\Rbm_e}{\xbm_e}^T\pder{\Rbm_e}{\xbm_e}+\pder{\Nbm_e}{\xbm_e}^T\pder{\Nbm_e}{\xbm_e}\right)\Psibold_e,
\end{aligned}
\end{equation}
where all terms involving $\Rbm_e$ and $\Nbm_e$ are evaluated at
$(\Phibold_e\wbm_n,\Phibold_e'\wbm_n,\Xbm_e+\Psibold_e\taubold_n,\mubold)$ and
$\Xbm_e+\Psibold_e\taubold_n$, respectively. From the elemental form of the
hyperreduced residuals (\ref{eqn:ift-eqp-res0}) and Hessian approximation
(\ref{eqn:levmarq-elem}), it is clear the Levenberg-Marquardt linear system can
be assembled efficiently (operations independent of $N_\mathtt{e}$) provided $\rhobold$
is sparse.

\begin{remark} \label{rem:solver}
Following the approach in \cite{mirhoseini2023model}, for a fixed $\mubold\in\Dcal$ and
$\rhobold\in\Rbb^{N_\mathtt{e}}$, the Levenberg-Marquardt iteration is initialized as
\begin{equation}
 \wbm_0 \coloneqq \argmin_{\wbm\in\Rbb^{k_\ubm}} \norm{\tilde\Rbm_{\Phibold,\Psibold,\rhobold}(\wbm,\taubold_0;\mubold)}_2, \qquad
 \taubold_0 = \zerobold,
\end{equation}
and the iterations terminate when the first-order optimality conditions
drop below tolerances. That is, given tolerances $\epsilon_1,\epsilon_2\geq 0$,
$(\wbm_n,\taubold_n)$ is considered a numerical solution of (\ref{eqn:ift-eqp}) if
\begin{equation}
 \norm{\tilde\Jbm_{\wbm,n}^T\tilde\Fbm_n}_2 \leq \epsilon_1, \qquad
 \norm{\tilde\Jbm_{\taubold,n}^T\tilde\Fbm_n}_2 \leq \epsilon_2.
\end{equation}
In this work, we take $\epsilon_1 = \epsilon_2 = 10^{-8}$.
\end{remark}

\subsection{Greedy training}
\label{sec:rom:greedy}
To close this section, we introduce a greedy procedure to construct the
reduced bases, $\Phibold$ and $\Psibold$, as well as the weight vector
$\rhobold$. Following a traditional greedy method \cite{maday2002priori},
we construct a finite sequence of parameters $\{\mubold_j\}_{j=1}^M$
at which to train the EQP-IFT, and in the process construct a sequence
of EQP-IFTs, $\{\Phibold_j\}_{j=1}^M$,
$\{\Psibold_j\}_{j=1}^M$, $\{\rhobold_j\}_{j=1}^M$.
The last entry of the EQP-IFT sequence will be one used for all
online computations, i.e., $\Phibold = \Phibold_M$, $\Psibold = \Psibold_M$,
$\rhobold = \rhobold_M$.

Let $\Dcal_\mathrm{cnd}\subset\Dcal$ be a finite collection of samples in the
parameter domain, called the candidate set, from which the training parameters
will be selected, and define the EQP training set as a subset of $\Dcal_\mathrm{cnd}$,
i.e., $\Xibold \subset \Dcal_\mathrm{cnd}$. Suppose the EQP-IFT at the $j$th entry
in the sequence is available, i.e., $\Phibold_j$, $\Psibold_j$, $\rhobold_j$. Then, the 
next training parameter is the member of $\Dcal_\mathrm{cnd}$ at which the HDM residual
evaluated at the EQP-IFT solution has the largest magnitude
\begin{equation} \label{eqn:greedy0}
 \mubold_{j+1} = \argmax_{\mubold\in\Dcal_\mathrm{cnd}}
  ~\norm{\Rbm(\tilde\Ubm_{\Phibold_j,\Psibold_j,\rhobold_j}^\star(\mubold); \tilde\xbm_{\Phibold_j,\Psibold_j,\rhobold_j}^\star(\mubold),\mubold)}.
\end{equation}
This mimics a traditional greedy search with the EQP-IFT in place of a linear
projection-based ROM and the HDM residual taken as the error estimator. Once
the new training parameter is obtained, the domain mapping that causes the
features in the solution at the new parameter configuration to align in
the reference domain with those in the previous parameters must be determined.
For this, we use the approach introduced in \cite{mirhoseini2023model} and define
\begin{equation} \label{eqn:align0}
 (\tilde\wbm_{j+1}, \xbm_{j+1}^\star) \coloneqq
 \argmin_{(\wbm,\xbm)\in\Rbb^j\times\Rbb^{N_\xbm}}~\frac{1}{2}\norm{\Fbm(\Phibold_j\wbm,\xbm;\mubold_{j+1})}_2^2,
\end{equation}
where $\xbm_{j+1}^\star$ contains the domain mapping coefficients that leads
to feature alignment and $\tilde\wbm_{j+1}$ are unused coefficients of the
state vector. Next, the HDM solution is computed at the new parameter configuration
($\mubold_{j+1}$) and corresponding domain deformation ($\xbm_{j+1}^\star$), and
the new information is appended to the reduced bases
\begin{equation} \label{eqn:greedy1}
 \begin{aligned}
  \Phibold_{j+1} &= \mathtt{GramSchmidt}\left(\begin{bmatrix} \Phibold_j & \Ubm^\star(\xbm_{j+1}^\star,\mubold_{j+1}) \end{bmatrix}\right) \in \Rbb^{N_\ubm\times (j+1)} \\
  \Psibold_{j+1} &= \mathtt{GramSchmidt}\left(\begin{bmatrix} \Psibold_j & \xbm_{j+1}^\star - \Xbm \end{bmatrix}\right) \in \Rbb^{N_\xbm\times (j+1)},
 \end{aligned}
\end{equation}
where $\mathtt{GramSchmidt} : \Rbb^{m\times n} \rightarrow \Rbb^{m\times n}$ is a
operator that applies Gram-Schmidt orthogonalization to the columns of the input
matrix. Finally, the vector of weights is updated by solving the linear program
in (\ref{eqn:eqp-lp}) with the weight vector guess taken to be the weights from
the previous greedy iteration $\check\rhobold = \rhobold_j$
\begin{equation} \label{eqn:greedy2}
 \rhobold_{j+1} = \rhobold_{\Phibold_{j+1},\Psibold_{j+1},\rhobold_j}^\star.
\end{equation}
The greedy algorithm is summarized in Algorithm~\ref{alg:greedy}.

\begin{algorithm}
 \caption{Greedy training for hyperreduced model with implicit feature tracking}
 \label{alg:greedy}
 \begin{algorithmic}[1]
  \REQUIRE Current EQP-IFT $\Phibold_j$, $\Psibold_j$, $\rhobold_j$,
    candidate set $\Dcal_\mathrm{cnd}$, EQP training set $\Xibold$ and tolerance $\delta$
  \ENSURE Updated hyperreduced model $\Phibold_{j+1}$, $\Psibold_{j+1}$, $\rhobold_{j+1}$
  \STATE \textbf{Greedy search:} Solve (\ref{eqn:greedy0}) for $\mubold_{j+1}$ \label{alg:line:greedy}
  \STATE \textbf{Snapshot alignment:} Compute $\xbm_{j+1}^\star$ from (\ref{eqn:align0})
  \STATE \textbf{HDM solution:} Solve HDM on deformed domain $\Ubm^\star(\xbm_{j+1}^\star,\mubold_{j+1})$
  \STATE \textbf{Basis construction:} Compute $\Phibold_{j+1}$, $\Psibold_{j+1}$ from (\ref{eqn:greedy1})
  \STATE \textbf{EQP weights:} Compute $\rhobold_{j+1}$ from (\ref{eqn:greedy2}) \label{alg:line:eqp}
\end{algorithmic}
\end{algorithm}

\begin{remark}
We choose to initialize the greedy method as follows
\begin{equation}
 \mubold_1 =
   \argmin_{\mubold\in\Dcal_\mathrm{cnd}} ~\frac{1}{2}\norm{\mubold-\bar\mubold},
 \quad
 \Phibold_1 = \Ubm^\star(\Xbm,\mubold_1)/\norm{\Ubm^\star(\Xbm,\mubold_1)}_2, \quad
 \Psibold_1 = \Xbm/\norm{\Xbm}_2, \quad
 \rhobold_1 = \rhobold_{\Phibold_1,\Psibold_1,\onebold}^\star,
\end{equation}
where $\bar\mubold$ is the centroid of $\Dcal$. Because a weight vector has not yet
been the initial greedy iteration, we take $\check\rhobold = \onebold$, i.e., the
EQP weights are trained at the ROM-IFT solution. Additionally, at the initial parameter,
the reference and physical domain are identical, which determines the location of
the features as their physical location at $\mubold_1$. The corresponding
features at all subsequent parameters $\{\mubold_2,\dots,\mubold_M\}$ will
be moved to this location using domain deformations.
\end{remark}

\begin{remark}
In this work, we choose the EQP training set to be the entire candidate set,
$\Xibold = \Dcal_\mathrm{cnd}$, which is reasonable when the candidate set
is small-to-moderate in size.
\end{remark}


\begin{remark}
The offline phase of the ROM/EQP-IFT is more expensive than the offline phase usually
associated with the method of snapshots due to the additional snapshot alignment step,
which is required in addition to standard offline computations. However, additional cost
of snapshot alignment is relatively small because it leverages the low-dimensional
manifold approximation and ROM-IFT method. Additionally, as we show in
Section~\ref{sec:numexp}, the ROM/EQP-IFT provides accurate predictions with limited
training, which further reduces the offline cost.
\end{remark}

\section{Numerical experiments}
\label{sec:numexp}
In this section, we investigate the performance of the proposed EQP-IFT
method on two computational fluid dynamics benchmark problems. In
particular, we show that with limited training, the EQP-IFT method
provides accurate predictions across a test set at a fraction of the
cost of directly using the HDM.

For any $\mubold\in\Dcal$, the ROM/EQP-IFT solution
$(\hat\Ubm_{\Phibold_j,\Psibold_j}^\star(\mubold),\hat\xbm_{\Phibold_j,\Psibold_j}^\star(\mubold))$
is defined such that $\hat\Ubm_{\Phibold_j,\Psibold_j}^\star(\mubold)$ approximates the
HDM solution $\Ubm^\star(\hat\xbm_{\Phibold_j,\Psibold_j}^\star(\mubold),\mubold)$; unlike
fixed-domain model reduction, the implicit feature tracking method chooses the
most suitable domain deformation ($\hat\xbm_{\Phibold_j,\Psibold_j}^\star(\mubold)$)
for each parameter $\mubold$. However, this is not necessarily
the domain mapping that leads to the most accurate PDE solution and could
lead to an overly optimistic diagnosis of the ROM/EQP-IFT error. For the
same reason, it is not appropriate to use a fixed-mesh HDM because both
problems exhibit large variation in the shock. To avoid these situations, we take
the HDM to be a high-order implicit shock tracking method \cite{huang2022robust},
which computes an appropriate mesh for each parameter configuration, and compare
the ROM/EQP-IFT directly to this solution in the physical domain
\begin{equation} \label{eqn:errmet}
\begin{aligned}
 \mubold\in\Dcal\mapsto \hat{e}_j(\mubold) &\coloneqq
  \sqrt{\frac{\int_\Omega \left|\hat{u}_j(x;\mubold) - \check{u}(x;\mubold)\right|^2\, dv}{\int_\Omega \left|\check{u}(x;\mubold)\right|^2 \, dv}} \\
 \mubold\in\Dcal\mapsto \tilde{e}_j(\mubold) &\coloneqq
  \sqrt{\frac{\int_\Omega \left|\tilde{u}_j(x;\mubold) - \check{u}(x;\mubold)\right|^2\, dv}{\int_\Omega \left|\check{u}(x;\mubold)\right|^2 \, dv}},
\end{aligned}
\end{equation}
where $(x;\mubold) \in \Omega \times \Dcal \mapsto \check{u}(x;\mubold)$ is the
shock-fitting reference solution (approximation to $u(x;\mubold)$),
$(x;\mubold) \in \Omega \times \Dcal \mapsto \hat{u}_j(x;\mubold)$ is the
functional representation of
$(\hat\Ubm_{\Phibold_j,\Psibold_j}^\star(\mubold),\hat\xbm_{\Phibold_j,\Psibold_j}^\star(\mubold))$
in the physical domain, and
$(x;\mubold) \in \Omega \times \Dcal \mapsto \tilde{u}_j(x;\mubold)$ is the
functional representation of
$(\tilde\Ubm_{\Phibold_j,\Psibold_j,\rhobold_j}^\star(\mubold),\tilde\xbm_{\Phibold_j,\Psibold_j,\rhobold_j}^\star(\mubold))$ in the physical domain.
This choice of error metric (\ref{eqn:errmet}) implies
the ROM/EQP-IFT errors will be nonzero, even at training points. To examine the
performance of the greedy method, we also consider the HDM residual evaluated at
the ROM/EQP-IFT solution
\begin{equation}
 \begin{aligned}
  \mubold \in \Dcal \mapsto \hat{r}_j(\mubold) &\coloneqq \Rbm(\hat\Ubm_{\Phibold_j,\Psibold_j}^\star(\mubold);\hat\xbm_{\Phibold_j,\Psibold_j}^\star(\mubold),\mubold) \\
  \mubold \in \Dcal \mapsto \tilde{r}_j(\mubold) &\coloneqq \Rbm(\tilde\Ubm_{\Phibold_j,\Psibold_j,\rhobold_j}^\star(\mubold);\tilde\xbm_{\Phibold_j,\Psibold_j,\rhobold_j}^\star(\mubold),\mubold).
 \end{aligned}
\end{equation}

To assess the parametric performance of ROM/EQP-IFT, we use the average
$L^2(\Omega)$ error over a finite test set $\Dcal_\mathrm{tst}$,
defined as
\begin{equation}
 \Dcal_\mathrm{tst}\subset\Dcal \mapsto \hat{E}_j(\Dcal_\mathrm{tst}) \coloneqq
 \frac{1}{|\Dcal_\mathrm{tst}|}\sum_{\mu\in\Dcal_\mathrm{tst}}~\hat{e}_j(\mu), \qquad
 \Dcal_\mathrm{tst}\subset\Dcal \mapsto \tilde{E}_j(\Dcal_\mathrm{tst}) \coloneqq
 \frac{1}{|\Dcal_\mathrm{tst}|}\sum_{\mu\in\Dcal_\mathrm{tst}}~\tilde{e}_j(\mu).
\end{equation}
Similarly for the HDM residuals
\begin{equation}
 \Dcal_\mathrm{tst}\subset\Dcal \mapsto \hat{R}_j(\Dcal_\mathrm{tst}) \coloneqq
 \frac{1}{|\Dcal_\mathrm{tst}|}\sum_{\mu\in\Dcal_\mathrm{tst}}~\hat{r}_j(\mu), \qquad
 \Dcal_\mathrm{tst}\subset\Dcal \mapsto \tilde{R}_j(\Dcal_\mathrm{tst}) \coloneqq
 \frac{1}{|\Dcal_\mathrm{tst}|}\sum_{\mu\in\Dcal_\mathrm{tst}}~\tilde{r}_j(\mu).
\end{equation}
The decoration is dropped, i.e.,
$e_j(\mubold)$, $r_j(\mubold)$, $E_j(\Dcal)$, $R_j(\Dcal)$,
to denote the error/residual of either the ROM-IFT or EQP-IFT method.

The computational cost of three components of the algorithm will be measured directly:
1) $T_j^\mathrm{gs}$ is the time to complete the greedy search
(Algorithm~\ref{alg:greedy}, Line~\ref{alg:line:greedy}),
2) $T_j^\rhobold$ is the time to solve the EQP linear program
(Algorithm~\ref{alg:greedy}, Line~\ref{alg:line:eqp}), and
3) $T_j^\mathrm{tst}$ is the average time to compute the ROM/EQP-IFT
solution at a test point. The timings will be presented
as dimensionless where the scaling factor is the average time for
a single HDM simulation, i.e., all timings reported with be in terms
of equivalent number of HDM solves. For both problems, we will also
report the fraction of element weights that are non-zero ($\Ncal$).

\subsection{Inviscid Burgers' equation}
\label{sec:numexp:burg}
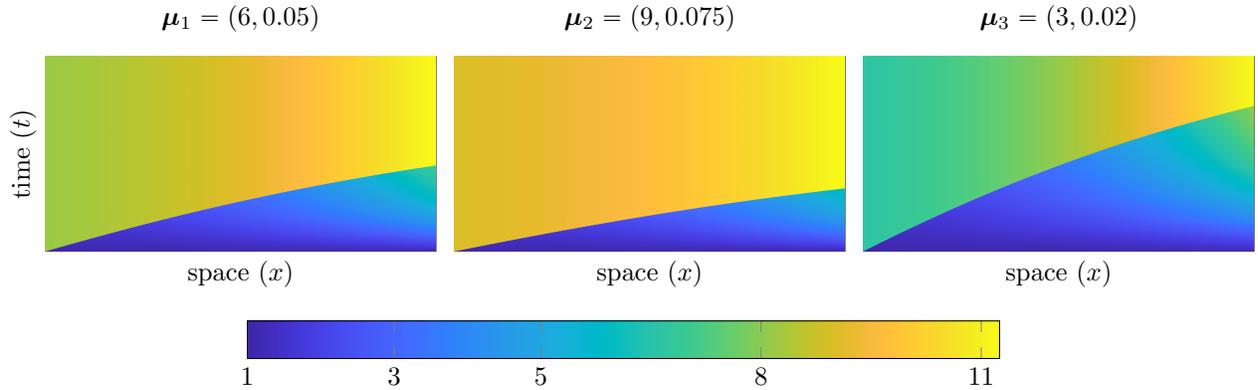
\begin{figure}[H]
\centering
\begin{tikzpicture}
\begin{groupplot}[
  group style={
      group size=3 by 1,
      horizontal sep=0.25cm
  },
  width=0.41\textwidth,
  axis equal image,
  xlabel={space ($x$)},
  xtick = {0.0,  100.0},
  ytick = {0.0,   50.0},
  xmin=0, xmax=100,
  ymin=0, ymax=50
]
\nextgroupplot[title={$\mubold_1=(6,0.05)$},xtick=\empty,ytick=\empty,ylabel={time ($t$)}]
\addplot graphics [xmin=0, xmax=100, ymin=0, ymax=50] {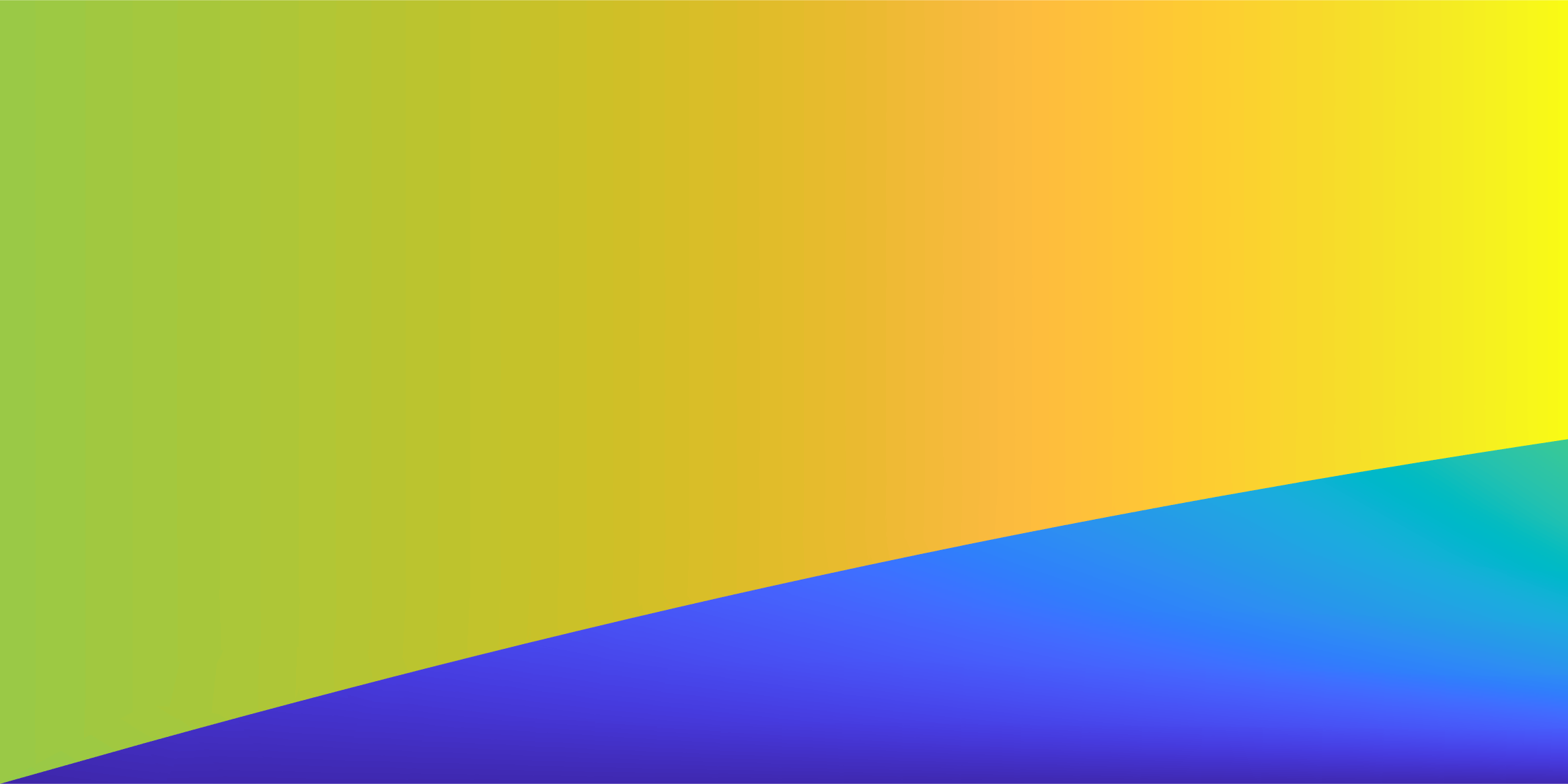};

\nextgroupplot[title={$\mubold_2=(9,0.075)$},xtick=\empty,ytick=\empty ]
\addplot graphics [xmin=0, xmax=100, ymin=0, ymax=50] {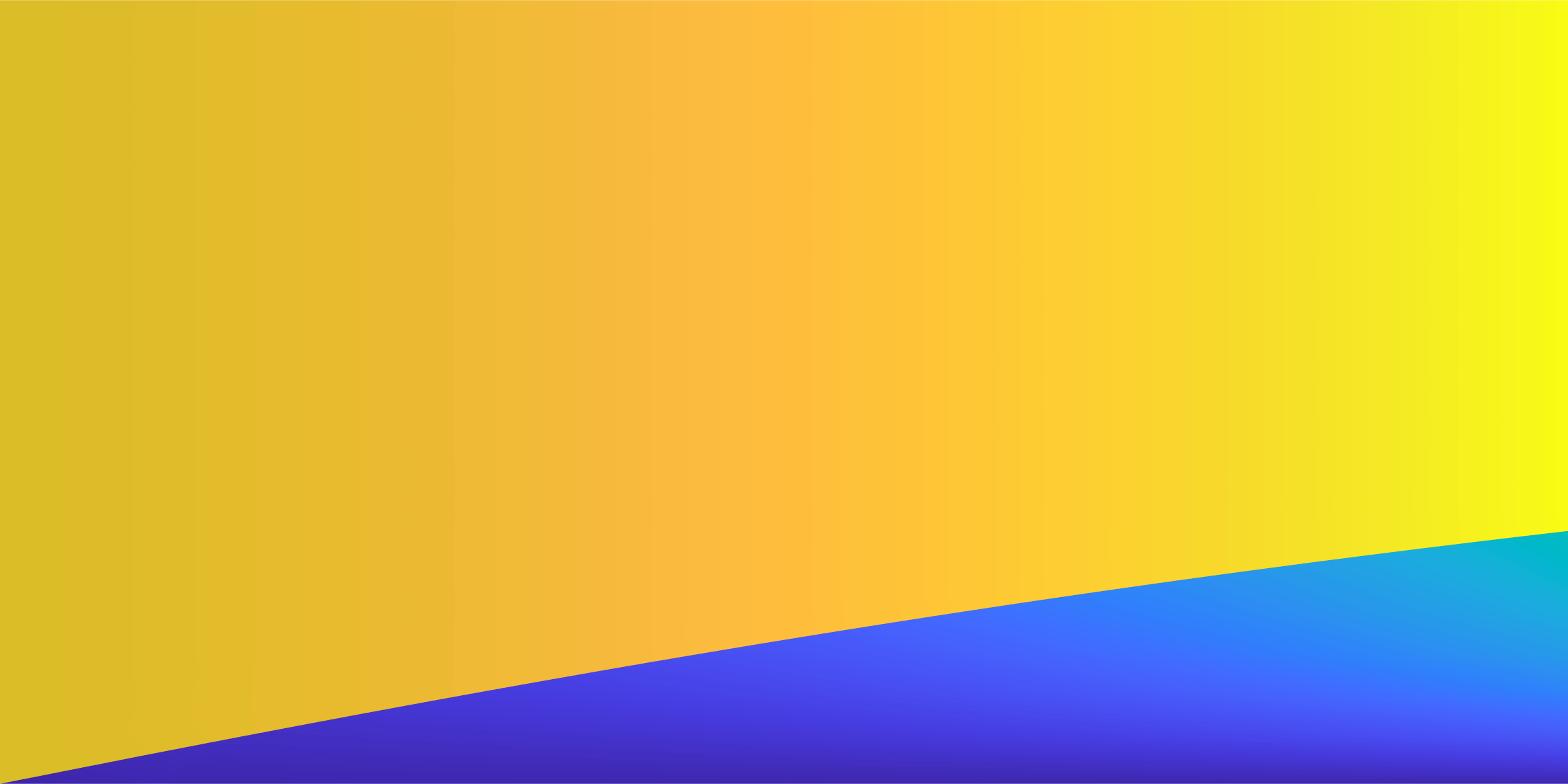};

\nextgroupplot[title={$\mubold_3=(3,0.02)$},xtick=\empty,ytick=\empty ]
\addplot graphics [xmin=0, xmax=100, ymin=0, ymax=50]{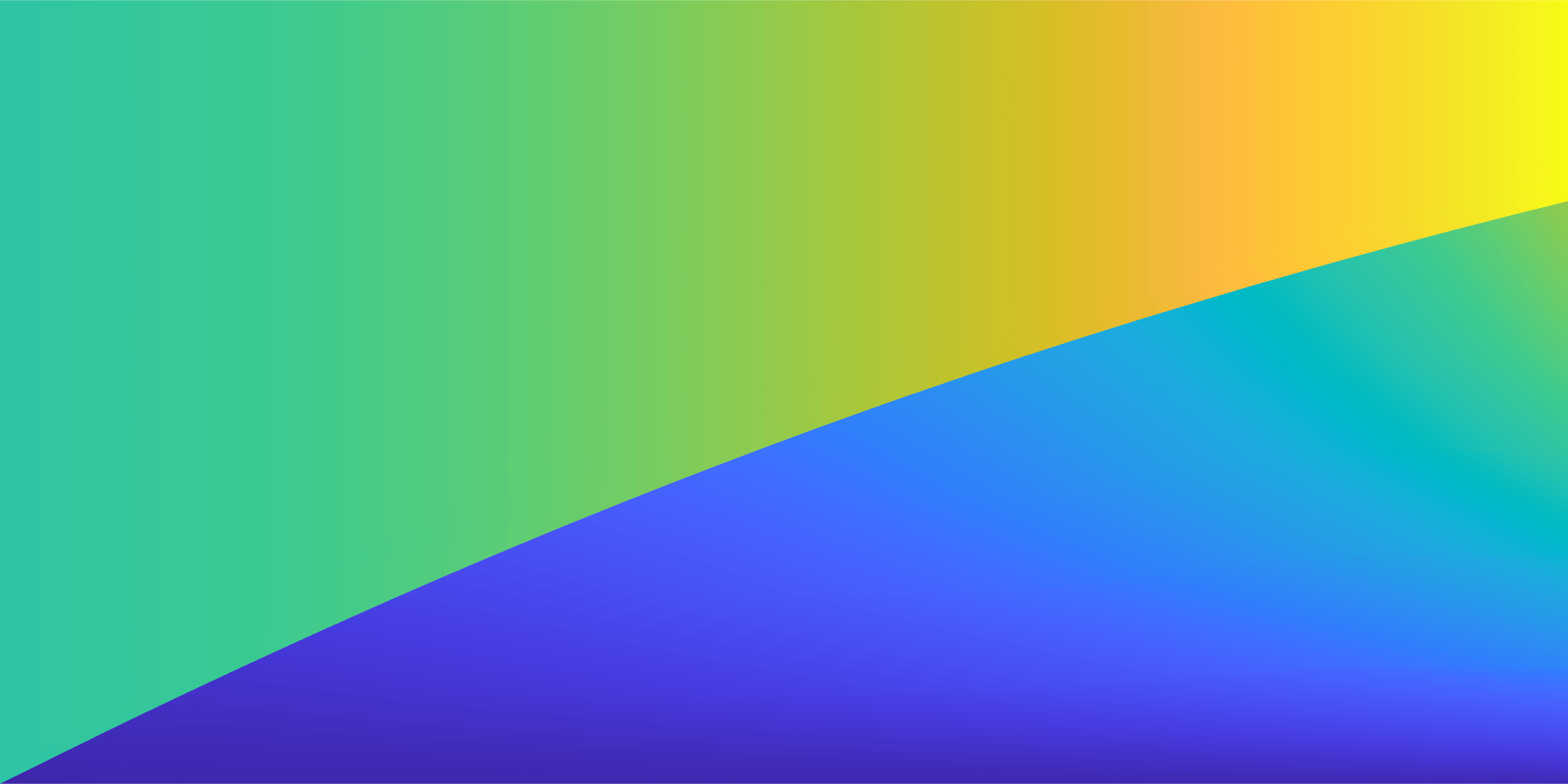};
\end{groupplot}
\end{tikzpicture}
\begin{tikzpicture}
\begin{axis}[
   hide axis, scale only axis,
   height=0pt, width=0pt,
   colormap={parula}{rgb255=(62,38,168) rgb255=(62,39,172) rgb255=(63,40,175) rgb255=(63,41,178) rgb255=(64,42,180) rgb255=(64,43,183) rgb255=(65,44,186) rgb255=(65,45,189) rgb255=(66,46,191) rgb255=(66,47,194) rgb255=(67,48,197) rgb255=(67,49,200) rgb255=(67,50,202) rgb255=(68,51,205) rgb255=(68,52,208) rgb255=(69,53,210) rgb255=(69,55,213) rgb255=(69,56,215) rgb255=(70,57,217) rgb255=(70,58,220) rgb255=(70,59,222) rgb255=(70,61,224) rgb255=(71,62,225) rgb255=(71,63,227) rgb255=(71,65,229) rgb255=(71,66,230) rgb255=(71,68,232) rgb255=(71,69,233) rgb255=(71,70,235) rgb255=(72,72,236) rgb255=(72,73,237) rgb255=(72,75,238) rgb255=(72,76,240) rgb255=(72,78,241) rgb255=(72,79,242) rgb255=(72,80,243) rgb255=(72,82,244) rgb255=(72,83,245) rgb255=(72,84,246) rgb255=(71,86,247) rgb255=(71,87,247) rgb255=(71,89,248) rgb255=(71,90,249) rgb255=(71,91,250) rgb255=(71,93,250) rgb255=(70,94,251) rgb255=(70,96,251) rgb255=(70,97,252) rgb255=(69,98,252) rgb255=(69,100,253) rgb255=(68,101,253) rgb255=(67,103,253) rgb255=(67,104,254) rgb255=(66,106,254) rgb255=(65,107,254) rgb255=(64,109,254) rgb255=(63,110,255) rgb255=(62,112,255) rgb255=(60,113,255) rgb255=(59,115,255) rgb255=(57,116,255) rgb255=(56,118,254) rgb255=(54,119,254) rgb255=(53,121,253) rgb255=(51,122,253) rgb255=(50,124,252) rgb255=(49,125,252) rgb255=(48,127,251) rgb255=(47,128,250) rgb255=(47,130,250) rgb255=(46,131,249) rgb255=(46,132,248) rgb255=(46,134,248) rgb255=(46,135,247) rgb255=(45,136,246) rgb255=(45,138,245) rgb255=(45,139,244) rgb255=(45,140,243) rgb255=(45,142,242) rgb255=(44,143,241) rgb255=(44,144,240) rgb255=(43,145,239) rgb255=(42,147,238) rgb255=(41,148,237) rgb255=(40,149,236) rgb255=(39,151,235) rgb255=(39,152,234) rgb255=(38,153,233) rgb255=(38,154,232) rgb255=(37,155,232) rgb255=(37,156,231) rgb255=(36,158,230) rgb255=(36,159,229) rgb255=(35,160,229) rgb255=(35,161,228) rgb255=(34,162,228) rgb255=(33,163,227) rgb255=(32,165,227) rgb255=(31,166,226) rgb255=(30,167,225) rgb255=(29,168,225) rgb255=(29,169,224) rgb255=(28,170,223) rgb255=(27,171,222) rgb255=(26,172,221) rgb255=(25,173,220) rgb255=(23,174,218) rgb255=(22,175,217) rgb255=(20,176,216) rgb255=(18,177,214) rgb255=(16,178,213) rgb255=(14,179,212) rgb255=(11,179,210) rgb255=(8,180,209) rgb255=(6,181,207) rgb255=(4,182,206) rgb255=(2,183,204) rgb255=(1,183,202) rgb255=(0,184,201) rgb255=(0,185,199) rgb255=(0,186,198) rgb255=(1,186,196) rgb255=(2,187,194) rgb255=(4,187,193) rgb255=(6,188,191) rgb255=(9,189,189) rgb255=(13,189,188) rgb255=(16,190,186) rgb255=(20,190,184) rgb255=(23,191,182) rgb255=(26,192,181) rgb255=(29,192,179) rgb255=(32,193,177) rgb255=(35,193,175) rgb255=(37,194,174) rgb255=(39,194,172) rgb255=(41,195,170) rgb255=(43,195,168) rgb255=(44,196,166) rgb255=(46,196,165) rgb255=(47,197,163) rgb255=(49,197,161) rgb255=(50,198,159) rgb255=(51,199,157) rgb255=(53,199,155) rgb255=(54,200,153) rgb255=(56,200,150) rgb255=(57,201,148) rgb255=(59,201,146) rgb255=(61,202,144) rgb255=(64,202,141) rgb255=(66,202,139) rgb255=(69,203,137) rgb255=(72,203,134) rgb255=(75,203,132) rgb255=(78,204,129) rgb255=(81,204,127) rgb255=(84,204,124) rgb255=(87,204,122) rgb255=(90,204,119) rgb255=(94,205,116) rgb255=(97,205,114) rgb255=(100,205,111) rgb255=(103,205,108) rgb255=(107,205,105) rgb255=(110,205,102) rgb255=(114,205,100) rgb255=(118,204,97) rgb255=(121,204,94) rgb255=(125,204,91) rgb255=(129,204,89) rgb255=(132,204,86) rgb255=(136,203,83) rgb255=(139,203,81) rgb255=(143,203,78) rgb255=(147,202,75) rgb255=(150,202,72) rgb255=(154,201,70) rgb255=(157,201,67) rgb255=(161,200,64) rgb255=(164,200,62) rgb255=(167,199,59) rgb255=(171,199,57) rgb255=(174,198,55) rgb255=(178,198,53) rgb255=(181,197,51) rgb255=(184,196,49) rgb255=(187,196,47) rgb255=(190,195,45) rgb255=(194,195,44) rgb255=(197,194,42) rgb255=(200,193,41) rgb255=(203,193,40) rgb255=(206,192,39) rgb255=(208,191,39) rgb255=(211,191,39) rgb255=(214,190,39) rgb255=(217,190,40) rgb255=(219,189,40) rgb255=(222,188,41) rgb255=(225,188,42) rgb255=(227,188,43) rgb255=(230,187,45) rgb255=(232,187,46) rgb255=(234,186,48) rgb255=(236,186,50) rgb255=(239,186,53) rgb255=(241,186,55) rgb255=(243,186,57) rgb255=(245,186,59) rgb255=(247,186,61) rgb255=(249,186,62) rgb255=(251,187,62) rgb255=(252,188,62) rgb255=(254,189,61) rgb255=(254,190,60) rgb255=(254,192,59) rgb255=(254,193,58) rgb255=(254,194,57) rgb255=(254,196,56) rgb255=(254,197,55) rgb255=(254,199,53) rgb255=(254,200,52) rgb255=(254,202,51) rgb255=(253,203,50) rgb255=(253,205,49) rgb255=(253,206,49) rgb255=(252,208,48) rgb255=(251,210,47) rgb255=(251,211,46) rgb255=(250,213,46) rgb255=(249,214,45) rgb255=(249,216,44) rgb255=(248,217,43) rgb255=(247,219,42) rgb255=(247,221,42) rgb255=(246,222,41) rgb255=(246,224,40) rgb255=(245,225,40) rgb255=(245,227,39) rgb255=(245,229,38) rgb255=(245,230,38) rgb255=(245,232,37) rgb255=(245,233,36) rgb255=(245,235,35) rgb255=(245,236,34) rgb255=(245,238,33) rgb255=(246,239,32) rgb255=(246,241,31) rgb255=(246,242,30) rgb255=(247,244,28) rgb255=(247,245,27) rgb255=(248,247,26) rgb255=(248,248,24) rgb255=(249,249,22) rgb255=(249,251,21) },
   colorbar horizontal,
   point meta min=1.000000e+00, point meta max=11.25,
   colorbar style={width=10cm, xtick={1.000000e+00,3.00000e+00,5.00000e+00,8.000000e+00,11}}
]
\addplot [draw=none] coordinates {(0,0)};
\end{axis}
\end{tikzpicture}
\caption{Reference solutions of the inviscid Burgers' equation problem at selected parameters.}
\label{fig:burg-offline}
\end{figure}

We begin by considering the time-dependent inviscid Burgers' equations in
the space-time domain $\Omega \coloneqq (0, 100) \times (0, 50)$ with the
parametrization introduced in \cite{rewienski2003trajectory}
\begin{equation} \label{eqn:burgers}
\begin{aligned}
 \partial_t u(x,t;\mubold) + u(x,t;\mubold)\partial_{x}u(x,t;\mubold) = s(x;\mubold) \qquad (x,t) \in \Omega,
\end{aligned}
\end{equation}
where $u(x,t;\mubold)\in\Rbb$ is implicitly defined as the solution, the parameter vector
is $\mubold\in \Dcal\coloneqq[3,9]\times[0.02,0.075]$, and the source term
$s: \Omega \times \Dcal \rightarrow \Rbb$ and initial/boundary conditions are
\begin{equation}
s(x; (a,b)) = 0.02 e^{b x}, \qquad
u(0, t; (a,b)) = a, \qquad
u(x, 0; (a,b)) = 1.
\end{equation}
The solution of the parametrized conservation law contains a propagating shock with
the overall magnitude of the solution growing as the shock propagates due to the source
term. The first parameter controls the inflow boundary and therefore the shock speed, and
the second parameter controls the width of the source term (Figure~\ref{fig:burg-offline}).

The conservation law is discretized using a standard nodal discontinuous Galerkin
method using $N_\mathtt{e}=860$ linear simplex elements for a total of $N_\ubm=2580$ degrees of
freedom. We also use the same mesh of linear simplex elements to discretize the domain
mapping $\Gcal$ with Lagrangian basis functions, which means the domain deformation
coefficients ($\xbm$) can be interpreted as the nodal coordinates of the mesh in the
physical domain. The reference domain and mesh are constructed using implicit shock-fitting
\cite{zahr_implicit_2020, huang2022robust} at the parameter used to seed the greedy
method ($\bar\mubold$), which avoids the need for shock capturing and provides an
accurate representation of the conservation law solution.

\begin{figure}[H]
\centering
\begin{tikzpicture}
\begin{groupplot}[
  group style={
      group size=3 by 2,
      horizontal sep=0.25cm,
      vertical sep=0.25cm
  },
  width=0.41\textwidth,
  axis equal image,
  xlabel={$x$},
  ylabel={$t$},
  xtick = {0.0,  100.0},
  ytick = {0.0,   50.0},
  xmin=0, xmax=100,
  ymin=0, ymax=50
]
\nextgroupplot[title={HDM}, xlabel={} ,xtick=\empty,ylabel={} ,ytick=\empty ]
\addplot graphics [xmin=0, xmax=100, ymin=0, ymax=50] {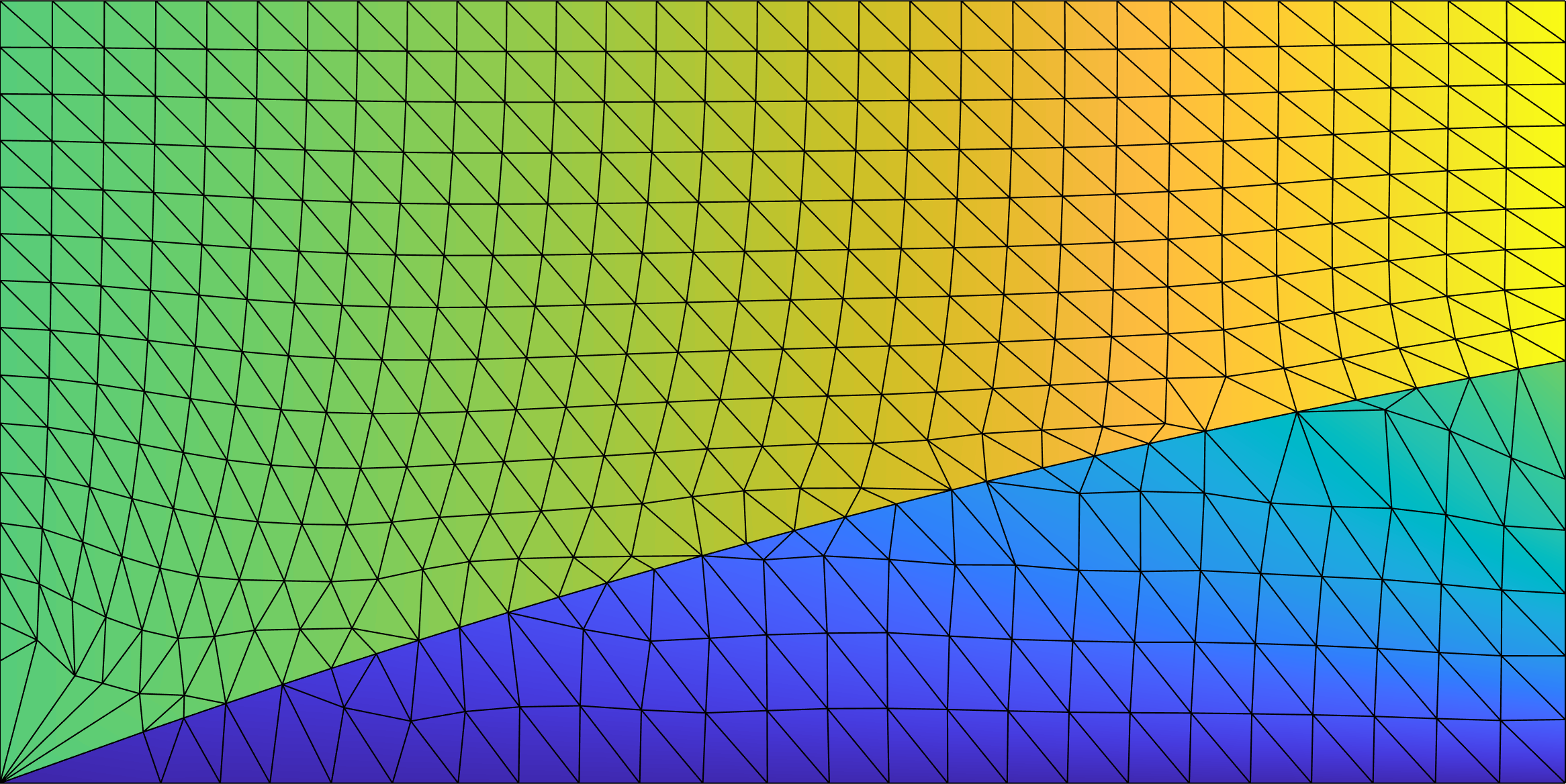};

\nextgroupplot[title={ROM-IFT}, ylabel={} ,ytick=\empty, xlabel={} ,xtick=\empty]
\addplot graphics [xmin=0, xmax=100, ymin=0, ymax=50]{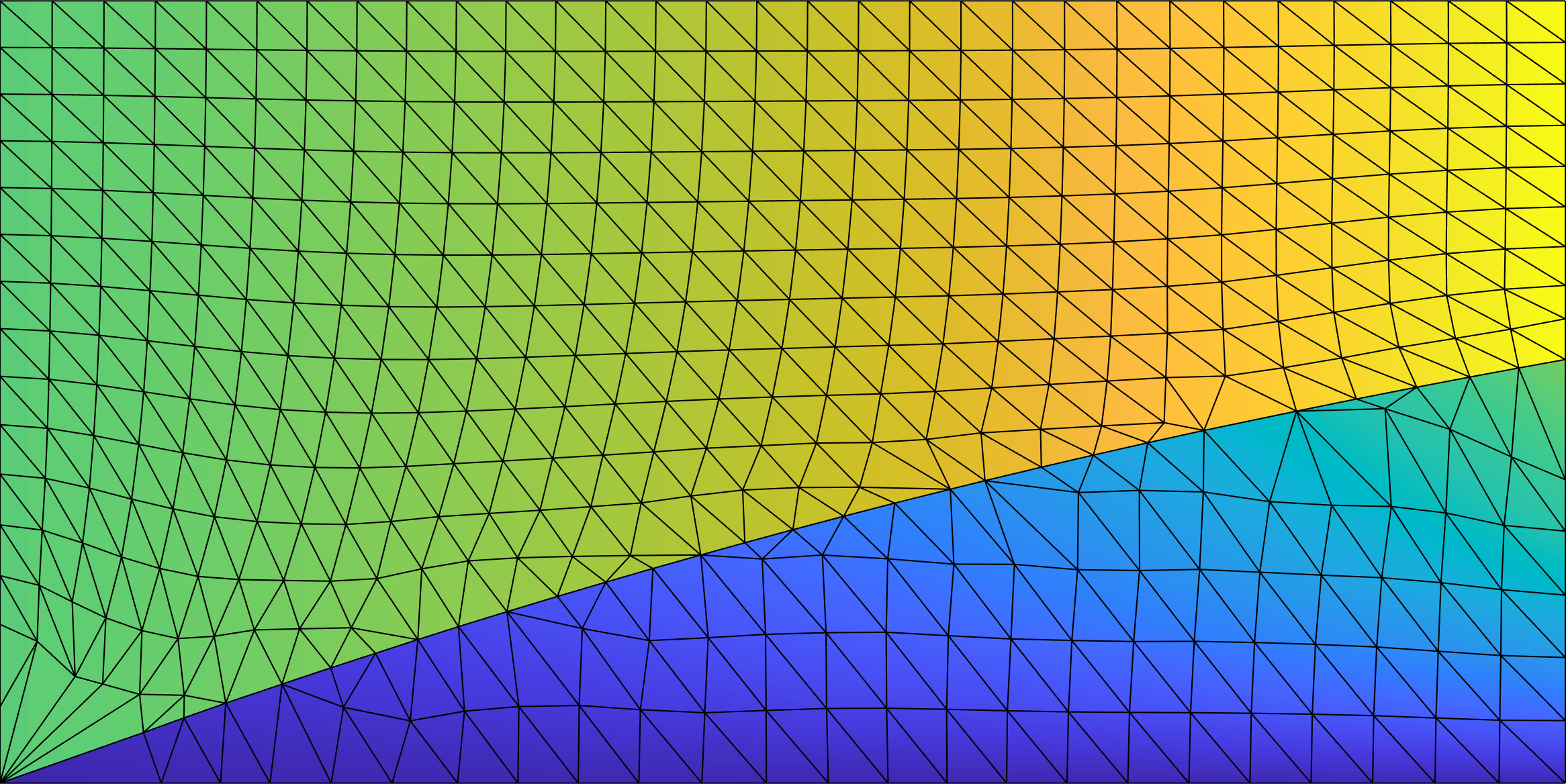};

\nextgroupplot[title={EQP-IFT ($\Ncal=6.05\%$)} , ylabel={},ytick=\empty, xlabel={} ,xtick=\empty]
\addplot graphics [xmin=0, xmax=100, ymin=0, ymax=50] {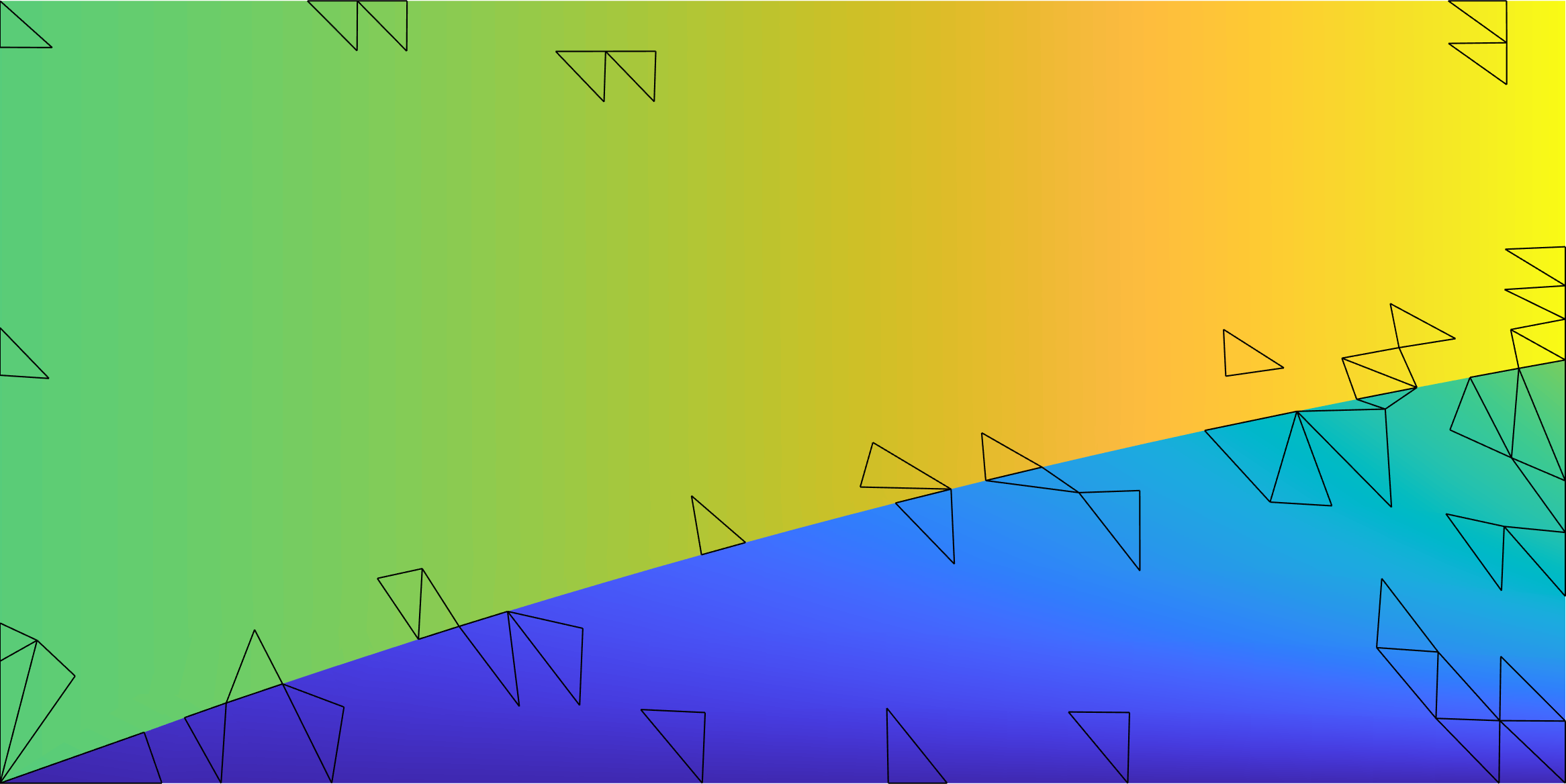};

\nextgroupplot[title={}, xlabel={} ,xtick=\empty,ylabel={} ,ytick=\empty ]
\addplot graphics [xmin=0, xmax=100, ymin=0, ymax=50] {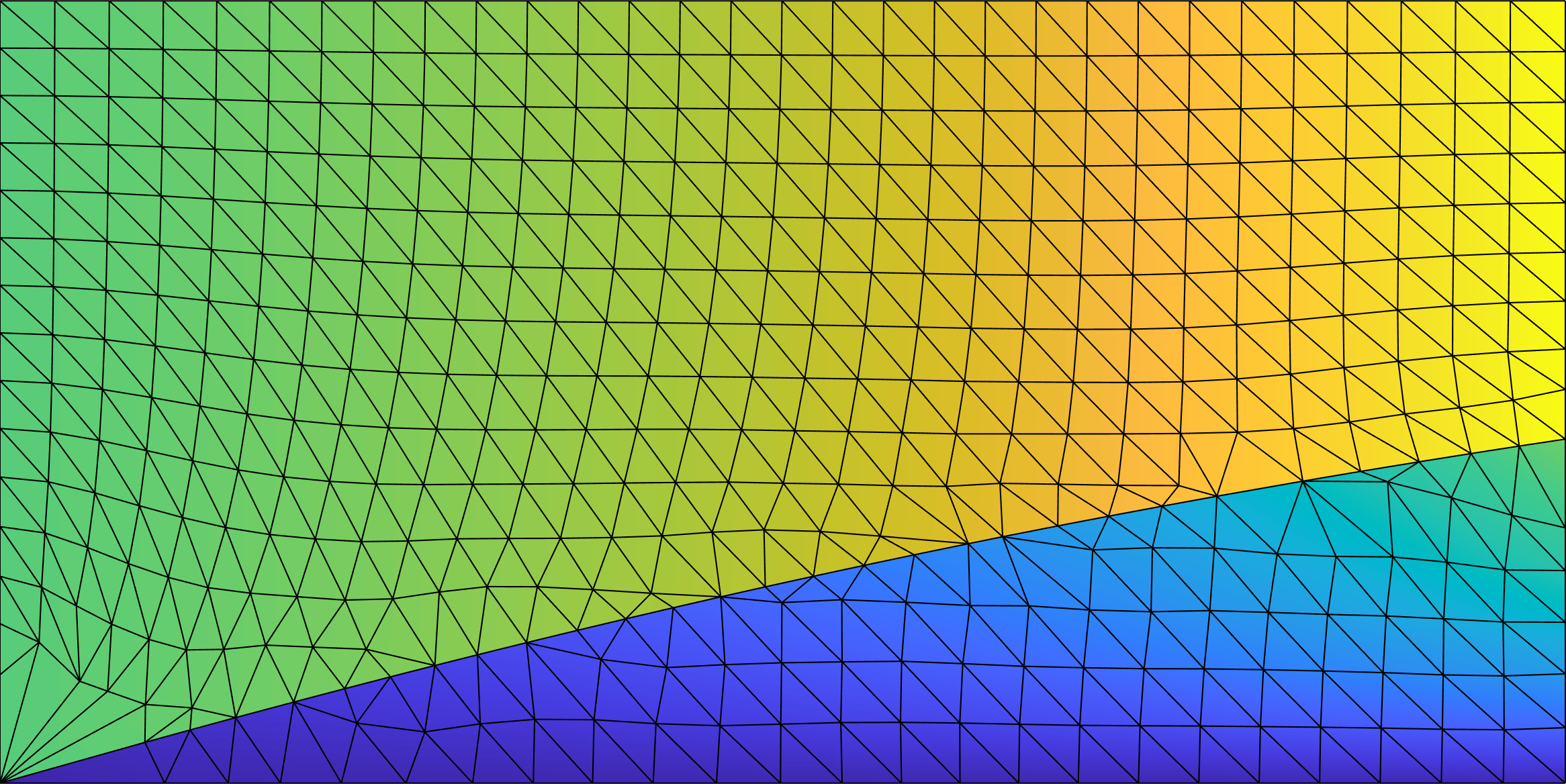};

\nextgroupplot[title={}, xlabel={} ,xtick=\empty,ylabel={} ,ytick=\empty]
\addplot graphics [xmin=0, xmax=100, ymin=0, ymax=50]{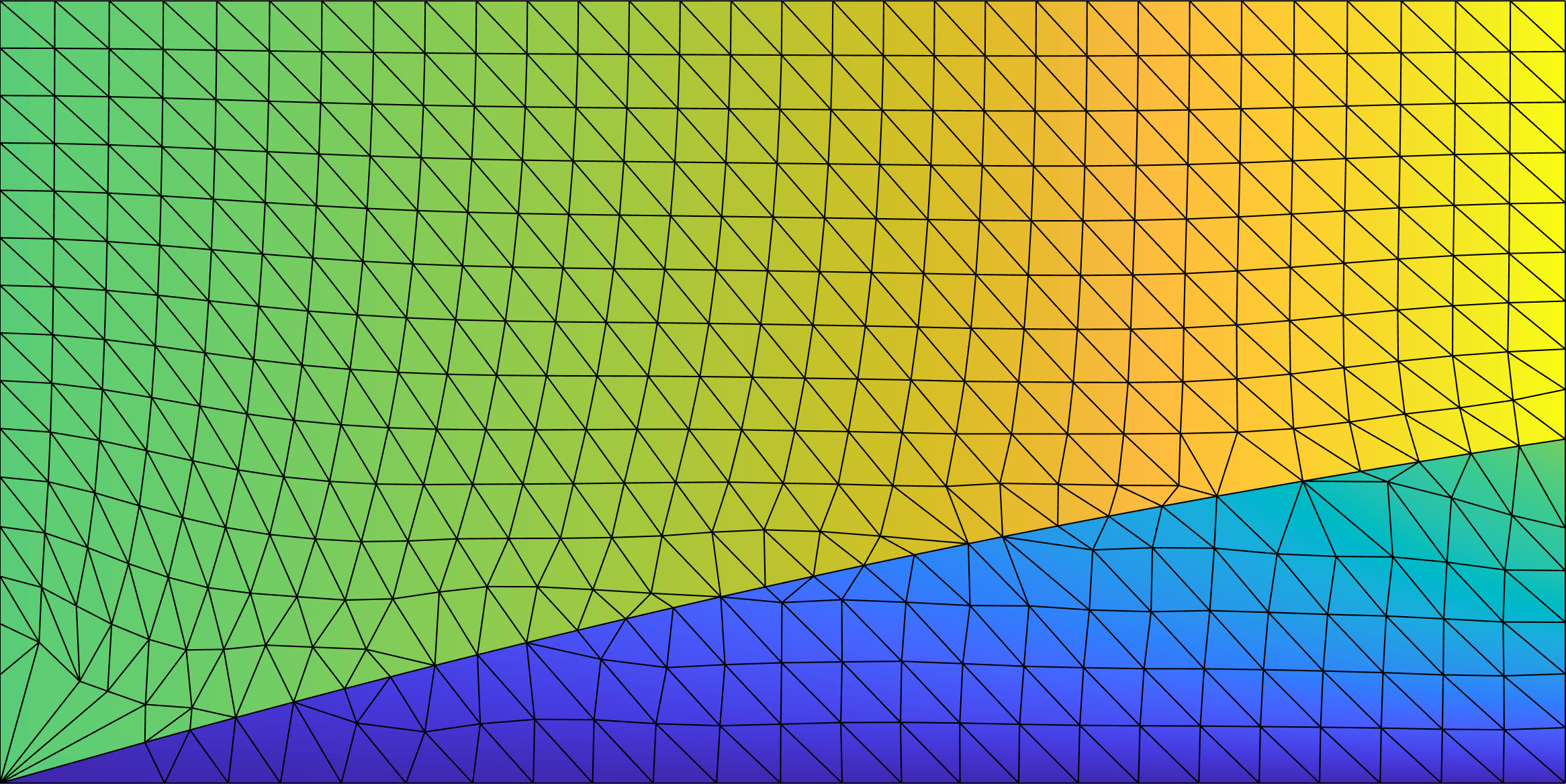};

\nextgroupplot[title={}, xlabel={} ,xtick=\empty,ylabel={} ,ytick=\empty]
\addplot graphics [xmin=0, xmax=100, ymin=0, ymax=50] {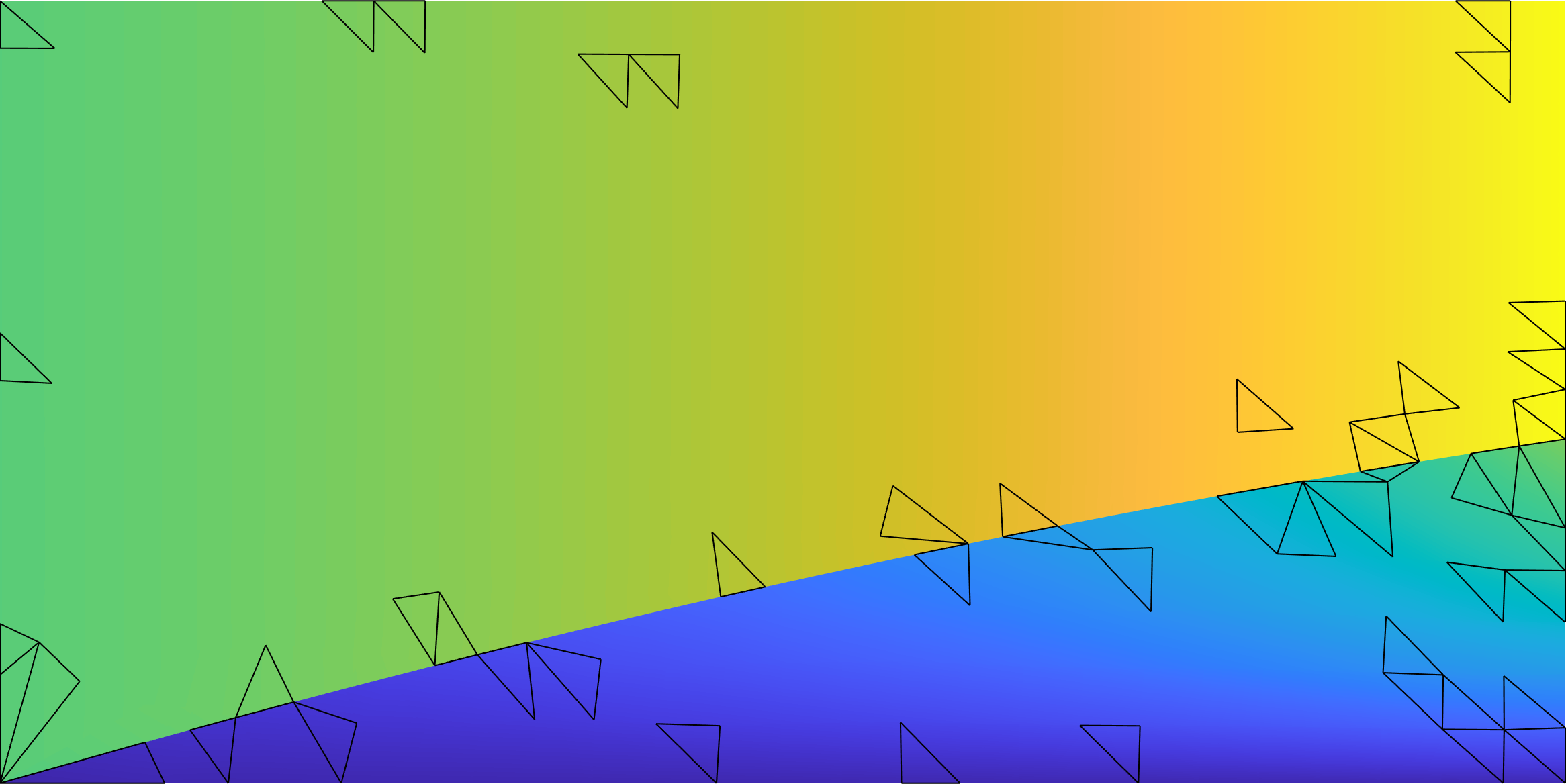};

\end{groupplot}
\end{tikzpicture}
\caption{Online solution of the inviscid Burgers' equation problem at $\mubold=(4.5,0.038)$
 in the physical domain (\textit{top}) and reference domain (\textit{bottom}). The mesh
 edges are included to show the elements over which the nonlinear terms must be assembled.}
\label{fig:burgers-online}
\end{figure}

\begin{figure}[H]
\centering
\input{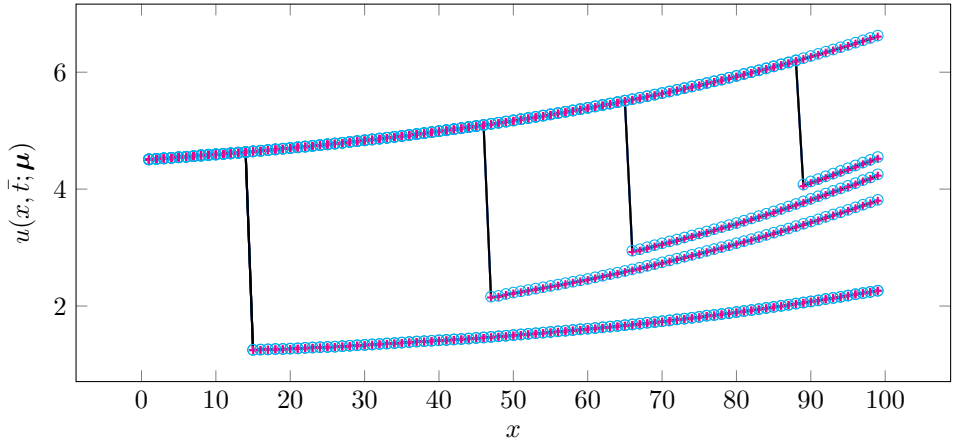}
\caption{Online solution of the inviscid Burgers' equation problem at $\mubold=(4.5,0.038)$
 along a five temporal slices $\{(s,5t) \mid s\in(0,100), t\in\{1,\dots,6\}\}$. Legend:
 HDM (\ref{line:burg_hdm}), ROM-IFT $j=3$ (\ref{line:burg_rom}), EQP-IFT $j=3$ (\ref{line:naca_hyper}).}
\label{fig:burgers-slice}
\end{figure}

We use the greedy algorithm proposed in Section~\ref{sec:rom:greedy} with a
candidate set $\Dcal_\mathrm{cnd}$ with size $|\Dcal_\mathrm{cnd}|=9$ constructed
by uniformly sampling $\Dcal$ with three samples per dimension. We use this procedure
to build up a ROM-IFT and EQP-IFT of sizes $j=1,2,3$. For EQP-IFT, we take the
EQP training set to be the candidate set ($\Xibold = \Dcal_\mathrm{cnd}$) and
the tolerance $\delta = 10^{-10}$. We test the performance of the reduced models
using the error metrics introduced in Section~\ref{sec:numexp} over a test set
$\Dcal_\mathrm{tst}$ with size $|\Dcal_\mathrm{tst}|=49$ constructed by uniformly
sampling $\Dcal$ with seven samples per dimension. The performance of the ROM-IFT
and EQP-IFT is summarized in Table~\ref{tab:burg}, which shows the both the offline
and online cost of EQP-IFT are reduced compared to ROM-IFT and the accuracy is comparable.
Figure~\ref{fig:burg-offline} shows the HDM
snapshot and domain deformation at the training parameters selected by the greedy
algorithm. Even at the parameter in the test set at which the EQP-IFT error is largest
($\mubold=(4.5,0.038)$), the ROM/EQP-IFT provide accurate approximations as can be seen
from the state throughout the domain (Figure~\ref{fig:burgers-online}) and along
a slice through the domain (Figure~\ref{fig:burgers-slice}).

\begin{table}[H]
\centering
\caption{Performance of HDM, ROM-IFT, EQP-IFT on test set $\Dcal_\mathrm{tst}$ for the inviscid Burgers' equation problem. Timings are reported as normalized by the cost to compute a single HDM solution (median CPU time for HDM = 58.2s).}
\label{tab:burg}
\begin{tabular}{r|ccccccc}
& $j$ & $\Ncal$ & ${E}_j$ & ${R}_j$ & $T_j^\mathrm{gs}$ & $T_j^\rhobold$ & $T_j^\mathrm{tst}$ \\
\hline ROM-IFT &1&-&0.107&18.2&-&-&0.24\\
 &2&-&0.0311&4.47&11.3&-&0.33\\
 &3&-&0.0124&2.21&16.1&-&0.41\\
\hline EQP-IFT&1&2.1\%&0.0864&35.4&-&2.4&0.014 \\
&2&4.3\%&0.0221&5.84&1.76&0.95&0.039 \\
&3&6.1\%&0.0120&3.81&2.78&1.4&0.057
\end{tabular}
\end{table}


\begin{figure}[H]
\centering
\input{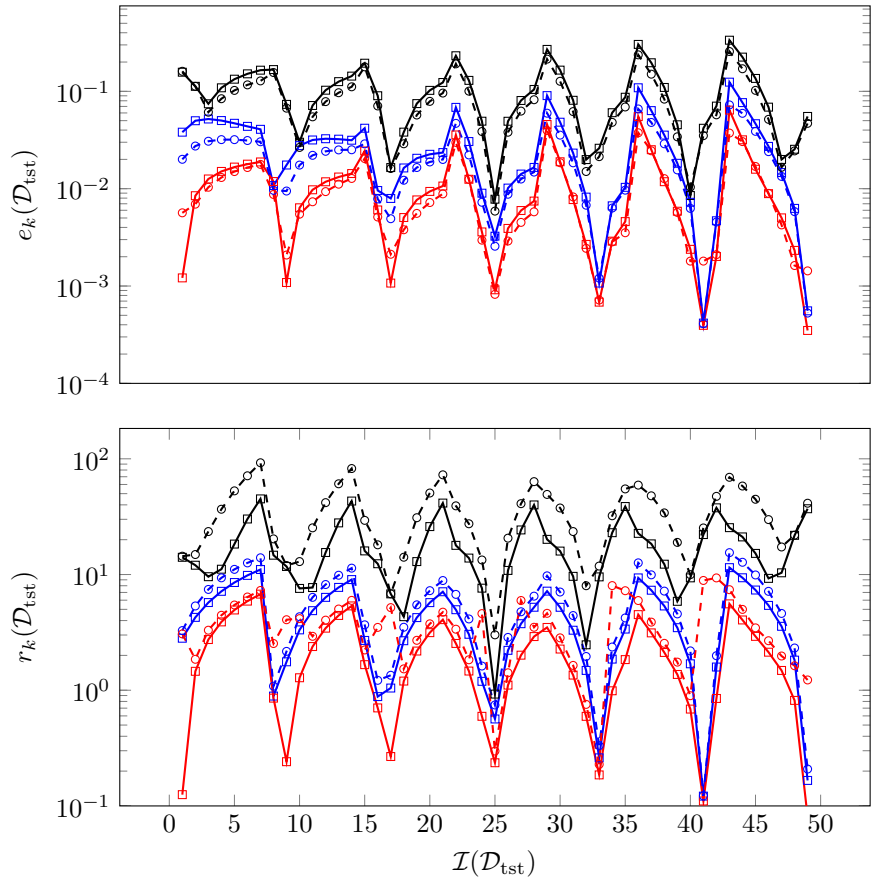}
\caption{Error (\textit{top}) and residual (\textit{bottom}) of the ROM/EQP-IFT
 at each test parameter at each greedy iteration for the ROM-IFT with $j=1$
 (\ref{line:burg_rom1}), $j=2$ (\ref{line:burg_rom2}), $j=3$ (\ref{line:burg_rom3}),
 and EQP-IFT with $j=1$ ($\Ncal=2.1\%$) (\ref{line:burg_eqp1}), $j=2$ ($\Ncal=4.3\%$)
 (\ref{line:burg_eqp2}), $j=3$ ($\Ncal=6.05\%$) (\ref{line:burg_eqp3}).}
\label{fig:burgers-greedy}
\end{figure}

Figure~\ref{fig:burgers-greedy} provides a detailed view into the behavior of the greedy
algorithm. We show variations in the two-dimensional parameter space by enumerating all
49 parameters in $\Dcal_\mathrm{tst}$ into a set denoted $\Ical(\Dcal_\mathrm{tst})$
with the first parameter being the fastest varying.
After the greedy algorithm is seeded with $\mubold_1=(6,0.05)$, the ROM/EQP-IFT
error is relatively small at many testing points, although exceeds 40\% at some testing
parameters. The sawtooth behavior is attributed to the parametrization of the
problem, where the first parameter (shock speed) has the largest impact on the
shock configuration. Thus, all seven test points whose first parameter matches
that of the single training parameter ($\mubold_1$) have small errors and the
error increases rapidly as the first entry of the parameter varies. This behavior
persists as additional training parameters are added. After the second greedy
iteration, there is a substantial drop in error across the parameter domain, which
is most significant near the training parameter selected $\mubold_2=(9,0.075)$.
When the third training parameter $\mubold_3=(3,0.02)$ is added, the reduction
in error is localized near new training parameters (modest global improvement).
This is consistent with observations in \cite{mirhoseini2023model} that showed
the ROM-IFT provides accurate approximations with limited training and additional
training leads to local improvements in the parameter space. Finally, there are only
minor differences between the ROM-IFT and EQP-IFT predictions. In some regions,
the EQP-IFT actually has a smaller error than the ROM-IFT, although it is
fortuitous that the additional error introduced by the weighted residual
makes the solution slightly closer to the reference solution. On the other
hand, the residual norm at the ROM-IFT solution is strictly less than the
that of the EQP-IFT solution, which is expected because the ROM-IFT method
minimizes the entire HDM residual (\ref{eqn:ift-opt}) whereas the EQP-IFT method only
minimizes the weighted HDM residual (\ref{eqn:ift-opt-eqp0}).

\subsection{Transonic flow over NACA0012}
\label{sec:numexp:euler}
\begin{figure}[H]
\centering
\begin{tikzpicture}
\begin{groupplot}[
  group style={
      group size=3 by 1,
      horizontal sep=0.4cm,
  },
  width=0.63\textwidth,
  axis equal image,
  xtick = \empty,
  ytick = \empty,
  xmin=-0.5, xmax=1.5,
  ymin=0, ymax=3
]
\nextgroupplot[title={$\mubold_1$ = 0.8}]
\addplot graphics [xmin=-0.5, xmax=1.5, ymin=0, ymax=3] {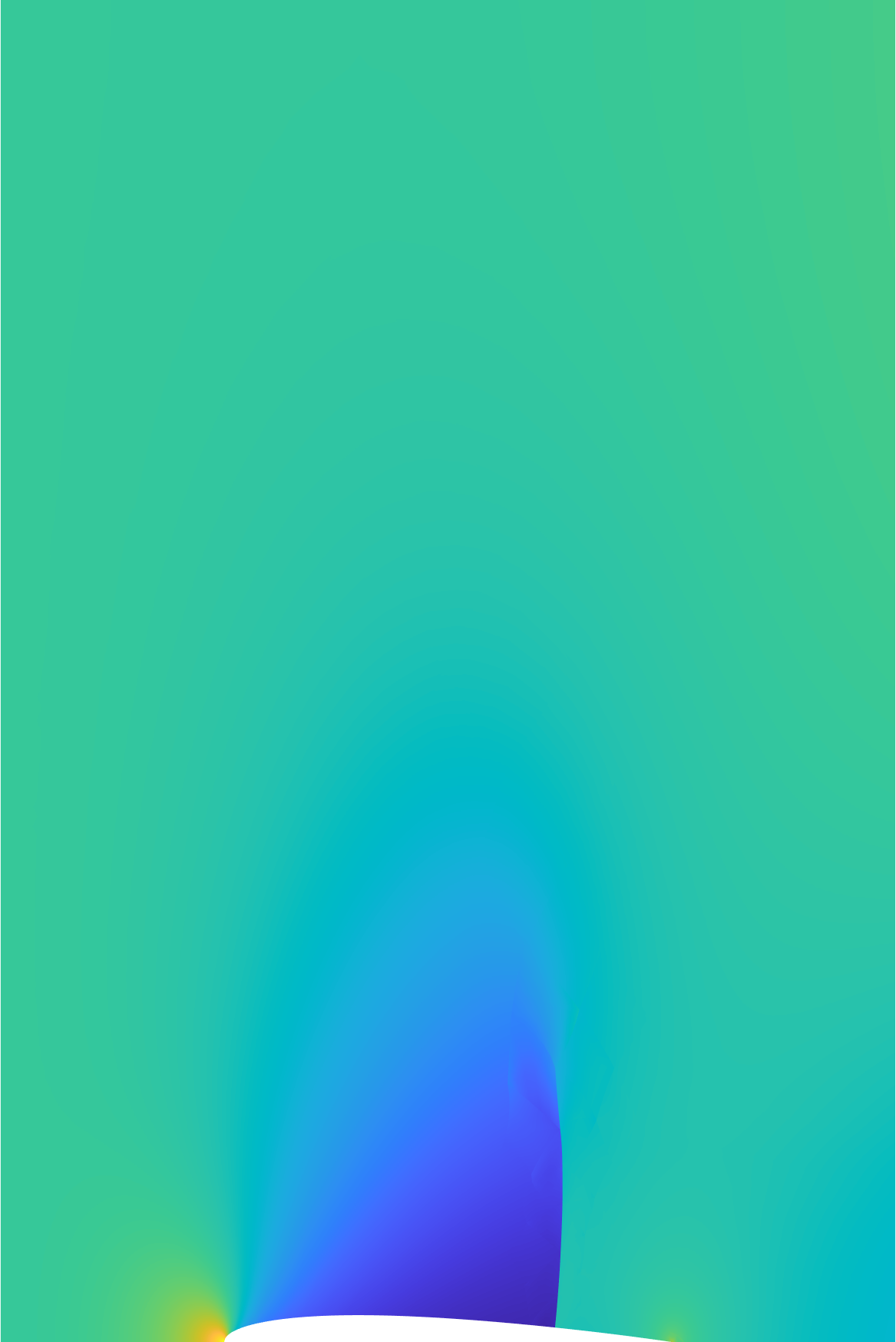};

\nextgroupplot[title={$\mubold_2$ = 0.75},ylabel={},ytick=\empty ]
\addplot graphics [xmin=-0.5, xmax=1.5, ymin=0, ymax=3] {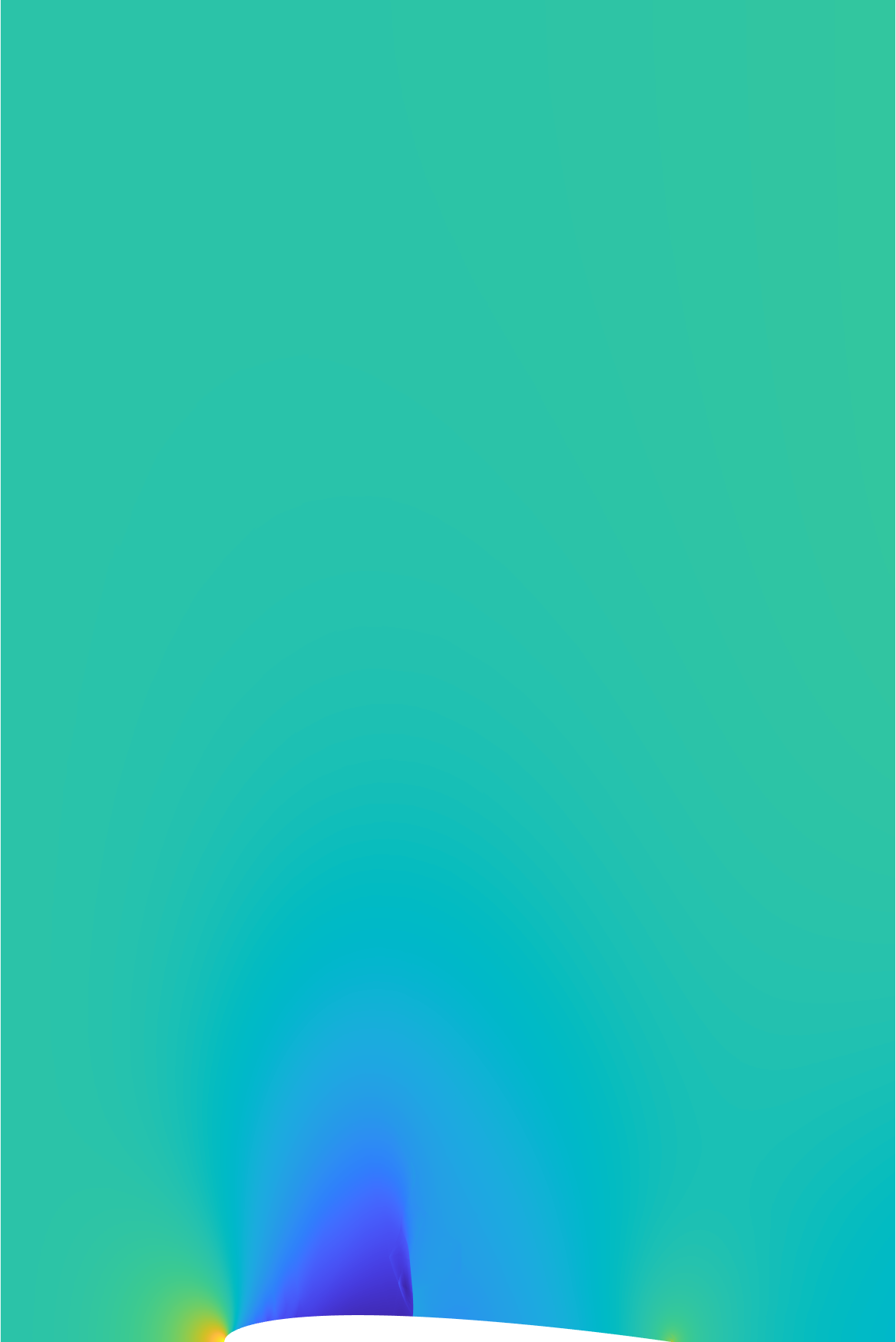};

\nextgroupplot[title={$\mubold_3$ = 0.82}, ylabel={},ytick=\empty]
\addplot graphics [xmin=-0.5, xmax=1.5, ymin=0, ymax=3]{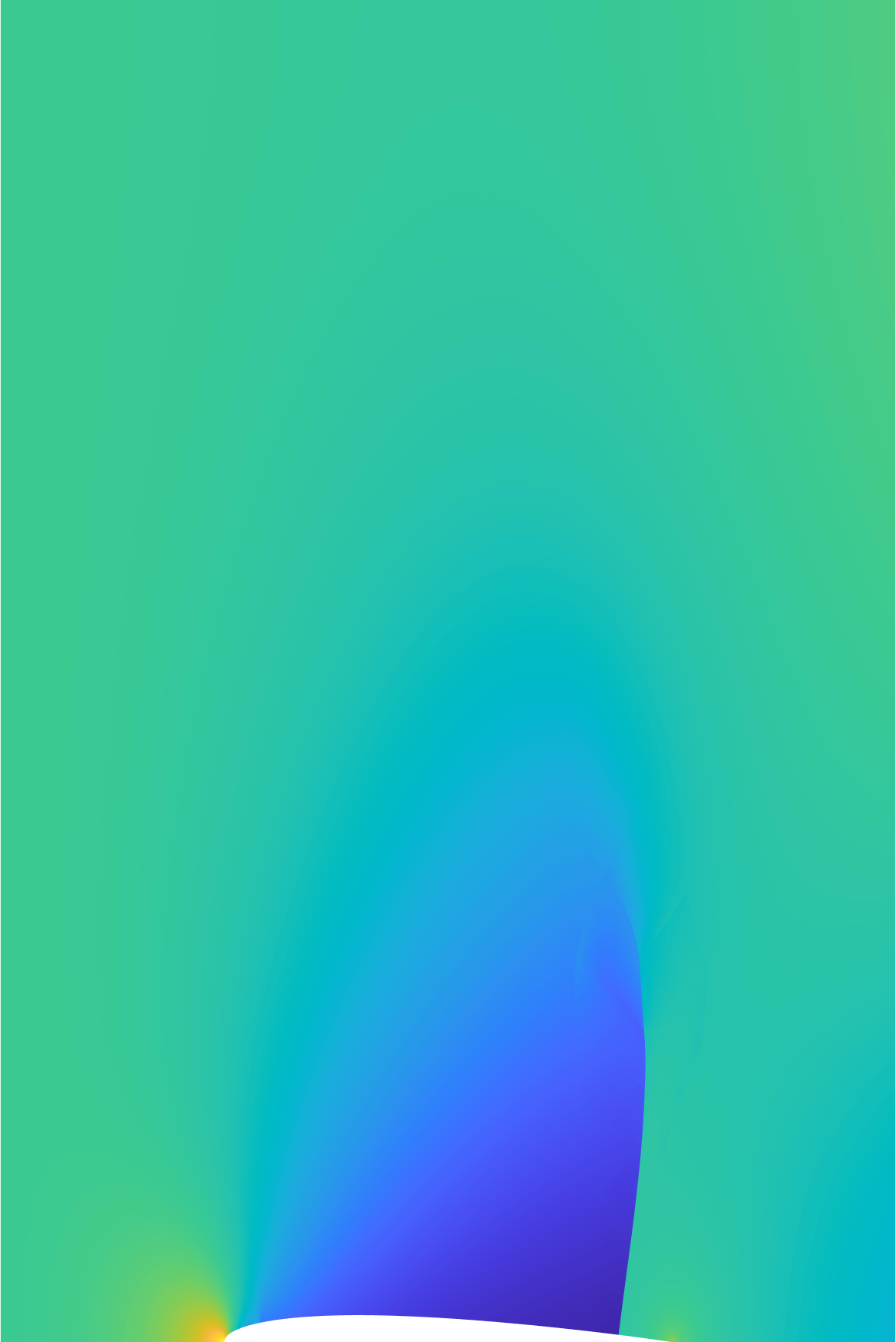};

\end{groupplot}
\end{tikzpicture}
\begin{tikzpicture}
\begin{axis}[
   hide axis, scale only axis,
   height=0pt, width=0pt,
   colormap={parula}{rgb255=(62,38,168) rgb255=(62,39,172) rgb255=(63,40,175) rgb255=(63,41,178) rgb255=(64,42,180) rgb255=(64,43,183) rgb255=(65,44,186) rgb255=(65,45,189) rgb255=(66,46,191) rgb255=(66,47,194) rgb255=(67,48,197) rgb255=(67,49,200) rgb255=(67,50,202) rgb255=(68,51,205) rgb255=(68,52,208) rgb255=(69,53,210) rgb255=(69,55,213) rgb255=(69,56,215) rgb255=(70,57,217) rgb255=(70,58,220) rgb255=(70,59,222) rgb255=(70,61,224) rgb255=(71,62,225) rgb255=(71,63,227) rgb255=(71,65,229) rgb255=(71,66,230) rgb255=(71,68,232) rgb255=(71,69,233) rgb255=(71,70,235) rgb255=(72,72,236) rgb255=(72,73,237) rgb255=(72,75,238) rgb255=(72,76,240) rgb255=(72,78,241) rgb255=(72,79,242) rgb255=(72,80,243) rgb255=(72,82,244) rgb255=(72,83,245) rgb255=(72,84,246) rgb255=(71,86,247) rgb255=(71,87,247) rgb255=(71,89,248) rgb255=(71,90,249) rgb255=(71,91,250) rgb255=(71,93,250) rgb255=(70,94,251) rgb255=(70,96,251) rgb255=(70,97,252) rgb255=(69,98,252) rgb255=(69,100,253) rgb255=(68,101,253) rgb255=(67,103,253) rgb255=(67,104,254) rgb255=(66,106,254) rgb255=(65,107,254) rgb255=(64,109,254) rgb255=(63,110,255) rgb255=(62,112,255) rgb255=(60,113,255) rgb255=(59,115,255) rgb255=(57,116,255) rgb255=(56,118,254) rgb255=(54,119,254) rgb255=(53,121,253) rgb255=(51,122,253) rgb255=(50,124,252) rgb255=(49,125,252) rgb255=(48,127,251) rgb255=(47,128,250) rgb255=(47,130,250) rgb255=(46,131,249) rgb255=(46,132,248) rgb255=(46,134,248) rgb255=(46,135,247) rgb255=(45,136,246) rgb255=(45,138,245) rgb255=(45,139,244) rgb255=(45,140,243) rgb255=(45,142,242) rgb255=(44,143,241) rgb255=(44,144,240) rgb255=(43,145,239) rgb255=(42,147,238) rgb255=(41,148,237) rgb255=(40,149,236) rgb255=(39,151,235) rgb255=(39,152,234) rgb255=(38,153,233) rgb255=(38,154,232) rgb255=(37,155,232) rgb255=(37,156,231) rgb255=(36,158,230) rgb255=(36,159,229) rgb255=(35,160,229) rgb255=(35,161,228) rgb255=(34,162,228) rgb255=(33,163,227) rgb255=(32,165,227) rgb255=(31,166,226) rgb255=(30,167,225) rgb255=(29,168,225) rgb255=(29,169,224) rgb255=(28,170,223) rgb255=(27,171,222) rgb255=(26,172,221) rgb255=(25,173,220) rgb255=(23,174,218) rgb255=(22,175,217) rgb255=(20,176,216) rgb255=(18,177,214) rgb255=(16,178,213) rgb255=(14,179,212) rgb255=(11,179,210) rgb255=(8,180,209) rgb255=(6,181,207) rgb255=(4,182,206) rgb255=(2,183,204) rgb255=(1,183,202) rgb255=(0,184,201) rgb255=(0,185,199) rgb255=(0,186,198) rgb255=(1,186,196) rgb255=(2,187,194) rgb255=(4,187,193) rgb255=(6,188,191) rgb255=(9,189,189) rgb255=(13,189,188) rgb255=(16,190,186) rgb255=(20,190,184) rgb255=(23,191,182) rgb255=(26,192,181) rgb255=(29,192,179) rgb255=(32,193,177) rgb255=(35,193,175) rgb255=(37,194,174) rgb255=(39,194,172) rgb255=(41,195,170) rgb255=(43,195,168) rgb255=(44,196,166) rgb255=(46,196,165) rgb255=(47,197,163) rgb255=(49,197,161) rgb255=(50,198,159) rgb255=(51,199,157) rgb255=(53,199,155) rgb255=(54,200,153) rgb255=(56,200,150) rgb255=(57,201,148) rgb255=(59,201,146) rgb255=(61,202,144) rgb255=(64,202,141) rgb255=(66,202,139) rgb255=(69,203,137) rgb255=(72,203,134) rgb255=(75,203,132) rgb255=(78,204,129) rgb255=(81,204,127) rgb255=(84,204,124) rgb255=(87,204,122) rgb255=(90,204,119) rgb255=(94,205,116) rgb255=(97,205,114) rgb255=(100,205,111) rgb255=(103,205,108) rgb255=(107,205,105) rgb255=(110,205,102) rgb255=(114,205,100) rgb255=(118,204,97) rgb255=(121,204,94) rgb255=(125,204,91) rgb255=(129,204,89) rgb255=(132,204,86) rgb255=(136,203,83) rgb255=(139,203,81) rgb255=(143,203,78) rgb255=(147,202,75) rgb255=(150,202,72) rgb255=(154,201,70) rgb255=(157,201,67) rgb255=(161,200,64) rgb255=(164,200,62) rgb255=(167,199,59) rgb255=(171,199,57) rgb255=(174,198,55) rgb255=(178,198,53) rgb255=(181,197,51) rgb255=(184,196,49) rgb255=(187,196,47) rgb255=(190,195,45) rgb255=(194,195,44) rgb255=(197,194,42) rgb255=(200,193,41) rgb255=(203,193,40) rgb255=(206,192,39) rgb255=(208,191,39) rgb255=(211,191,39) rgb255=(214,190,39) rgb255=(217,190,40) rgb255=(219,189,40) rgb255=(222,188,41) rgb255=(225,188,42) rgb255=(227,188,43) rgb255=(230,187,45) rgb255=(232,187,46) rgb255=(234,186,48) rgb255=(236,186,50) rgb255=(239,186,53) rgb255=(241,186,55) rgb255=(243,186,57) rgb255=(245,186,59) rgb255=(247,186,61) rgb255=(249,186,62) rgb255=(251,187,62) rgb255=(252,188,62) rgb255=(254,189,61) rgb255=(254,190,60) rgb255=(254,192,59) rgb255=(254,193,58) rgb255=(254,194,57) rgb255=(254,196,56) rgb255=(254,197,55) rgb255=(254,199,53) rgb255=(254,200,52) rgb255=(254,202,51) rgb255=(253,203,50) rgb255=(253,205,49) rgb255=(253,206,49) rgb255=(252,208,48) rgb255=(251,210,47) rgb255=(251,211,46) rgb255=(250,213,46) rgb255=(249,214,45) rgb255=(249,216,44) rgb255=(248,217,43) rgb255=(247,219,42) rgb255=(247,221,42) rgb255=(246,222,41) rgb255=(246,224,40) rgb255=(245,225,40) rgb255=(245,227,39) rgb255=(245,229,38) rgb255=(245,230,38) rgb255=(245,232,37) rgb255=(245,233,36) rgb255=(245,235,35) rgb255=(245,236,34) rgb255=(245,238,33) rgb255=(246,239,32) rgb255=(246,241,31) rgb255=(246,242,30) rgb255=(247,244,28) rgb255=(247,245,27) rgb255=(248,247,26) rgb255=(248,248,24) rgb255=(249,249,22) rgb255=(249,251,21) },
   colorbar horizontal,
   point meta min=0.9000000e+00, point meta max=2.04000000e+00,
   colorbar style={width= 9cm, xtick={1.000000e+00,1.2,1.400000e00,1.600000e00,1.8,2.00000e+00}}
]
\addplot [draw=none] coordinates {(0,0)};
\end{axis}
\end{tikzpicture}
\caption{Reference solutions of the transonic flow problem at selected parameters.}
\label{fig:naca-base}
\end{figure}
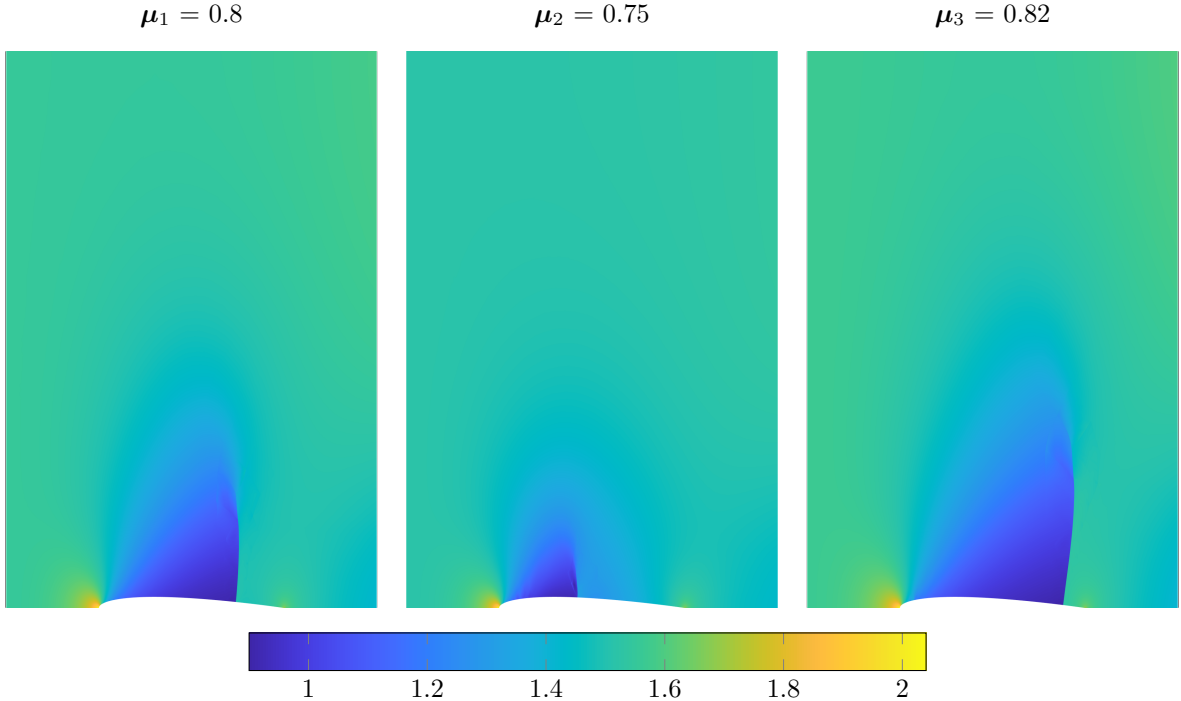

Next, we consider inviscid, transonic flow of an ideal gas over a
NACA0012 airfoil modeled by the compressible Euler equations. Let
$\Omega\subset\Rbb^d$ ($d=2$ in this problem) denote the flow
domain, then the Euler equations can be written as
an inviscid conservation law of the form (\ref{eqn:claw-phys})
with solution $\func{u}{\Omega\times\Dcal}{\Rbb^{d+2}}$
and components $u = (\rho, \rho v, \rho E)$, where
$\func{\rho}{\Omega\times\Dcal}{\Rbb_{>0}}$,
$\func{v}{\Omega\times\Dcal}{\Rbb^d}$,
$\func{E}{\Omega\times\Dcal}{\Rbb_{>0}}$
are the density, velocity, and energy of the fluid, respectively.
The flux function and source term are defined as
\begin{equation} \label{eqn:euler}
 f(u) \coloneqq \begin{bmatrix} \rho v^T \\ \rho v\otimes v+P(u)I_d \\ (\rho E+P(u))v^T \end{bmatrix},
 \qquad
 s(u) \coloneqq 0,
\end{equation}
where $u\mapsto P(u) \coloneqq (\gamma-1)(\rho E - \rho \norm{v}^2/2)$ is the
pressure of the fluid and $\gamma\in\Rbb_{>0}$ is the ratio of specific
heats. In this problem, we take the domain $\Omega$ to be the region around
the airfoil and the freestream density $\rho_\infty=1$, velocity
$v_\infty = (M_\infty c_\infty,0)$, pressure $P_\infty=1$, and
sound speed $c_\infty=\sqrt{\gamma}$.
The freestream Mach number $M_\infty\in[0.75,0.82]$ is used to
parametrize the problem. The solution has a curved shock attached
to the airfoil that strengthens and moves toward the tail as the
freestream Mach number increases (Figure~\ref{fig:naca-base}).

The conservation law is discretized using a standard nodal discontinuous Galerkin
method using $N_\mathtt{e}=814$ quadratic simplex elements for a total of $N_\ubm=19536$ degrees of
freedom. We also use the same mesh of quadratic simplex elements to discretize the domain
mapping $\Gcal$ with Lagrangian basis functions, which means the domain deformation
coefficients ($\xbm$) can be interpreted as the nodal coordinates of the mesh in the
physical domain. The reference domain and mesh are constructed using implicit shock-fitting
\cite{zahr_implicit_2020, huang2022robust} at the parameter used to seed the greedy
method ($\bar\mubold$), which avoids the need for shock capturing and provides an
accurate representation of the conservation law solution.

\begin{figure}
\centering
\begin{tikzpicture}
\begin{groupplot}[
  group style={
      group size=3 by 3,
      horizontal sep=0.1cm,
       vertical sep=0.1cm
  },
  width=0.58\textwidth,
  axis equal image,
  xlabel={$x$},
  ylabel={$\rho$},
  xtick = {-0.5,  1.5},
  ytick = {0.0,   3},
  xmin=-0.5, xmax=1.5,
  ymin=0, ymax=3
]

\nextgroupplot[title={HDM},xtick =\empty, ytick =\empty,
xlabel= {},ylabel= {}]
\addplot graphics [xmin=-0.5, xmax=1.5, ymin=0, ymax=3] {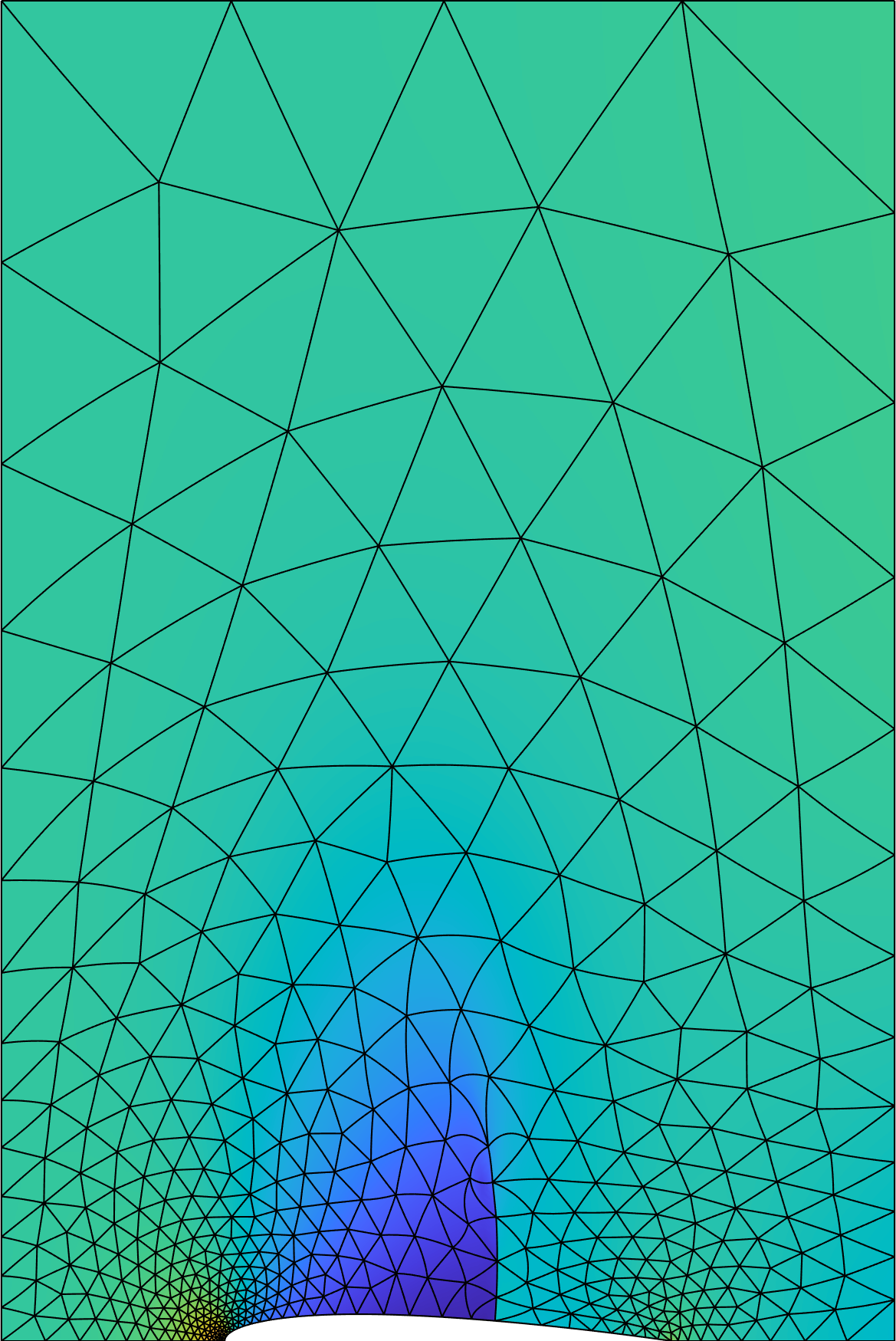};

\nextgroupplot[title={ROM-IFT},xtick =\empty, ytick =\empty,
xlabel= {},ylabel= {}]
\addplot graphics [xmin=-0.5, xmax=1.5, ymin=0, ymax=3]{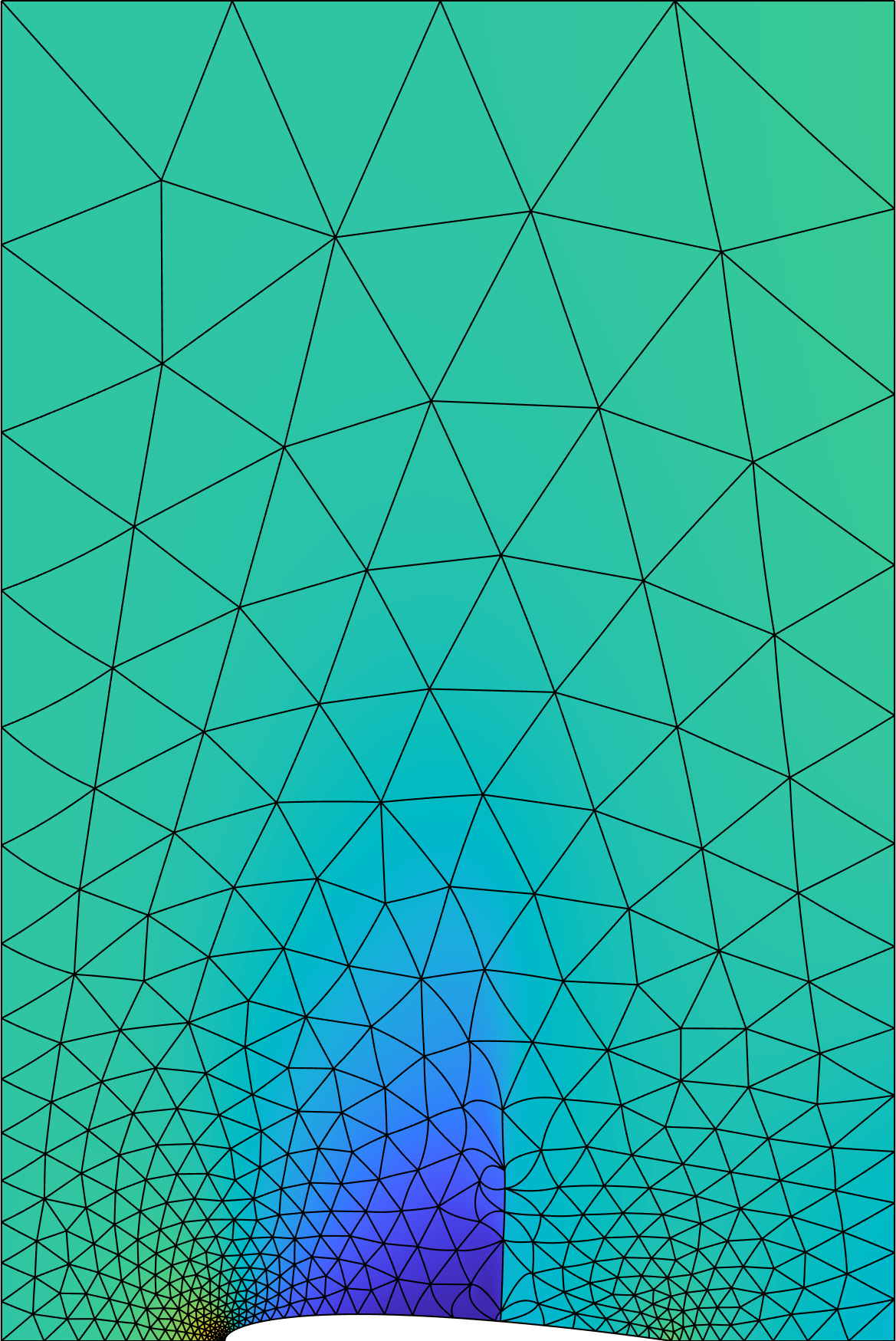};

\nextgroupplot[title={EQP-IFT ($\Ncal = 3.56\%$)} ,xtick =\empty, ytick =\empty,
xlabel= {},ylabel= {}]
\addplot graphics [xmin=-0.5, xmax=1.5, ymin=0, ymax=3] {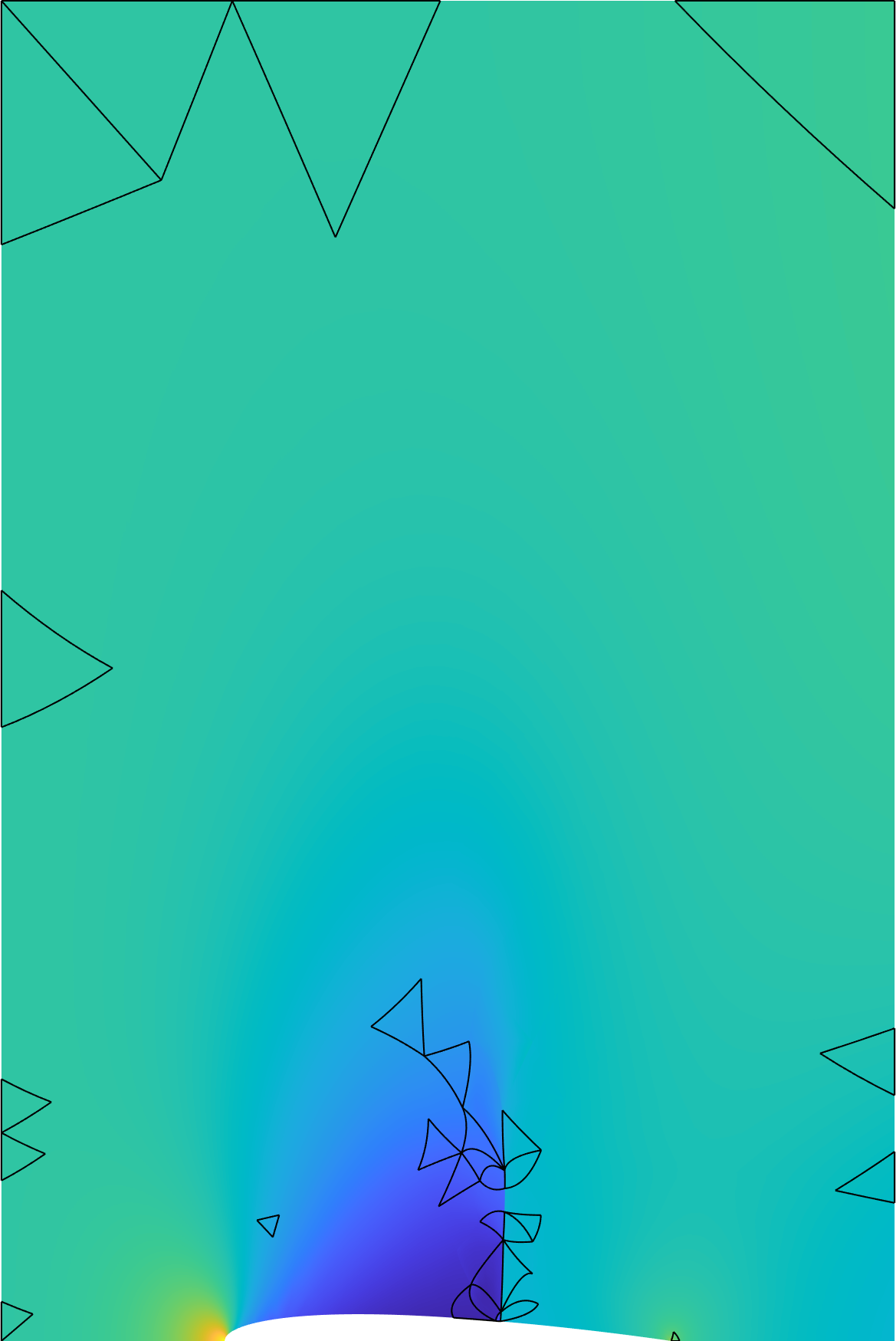};

\nextgroupplot[title={},xtick =\empty, ytick =\empty,
xlabel= {},ylabel= {}]
\addplot graphics [xmin=-0.5, xmax=1.5, ymin=0, ymax=3] {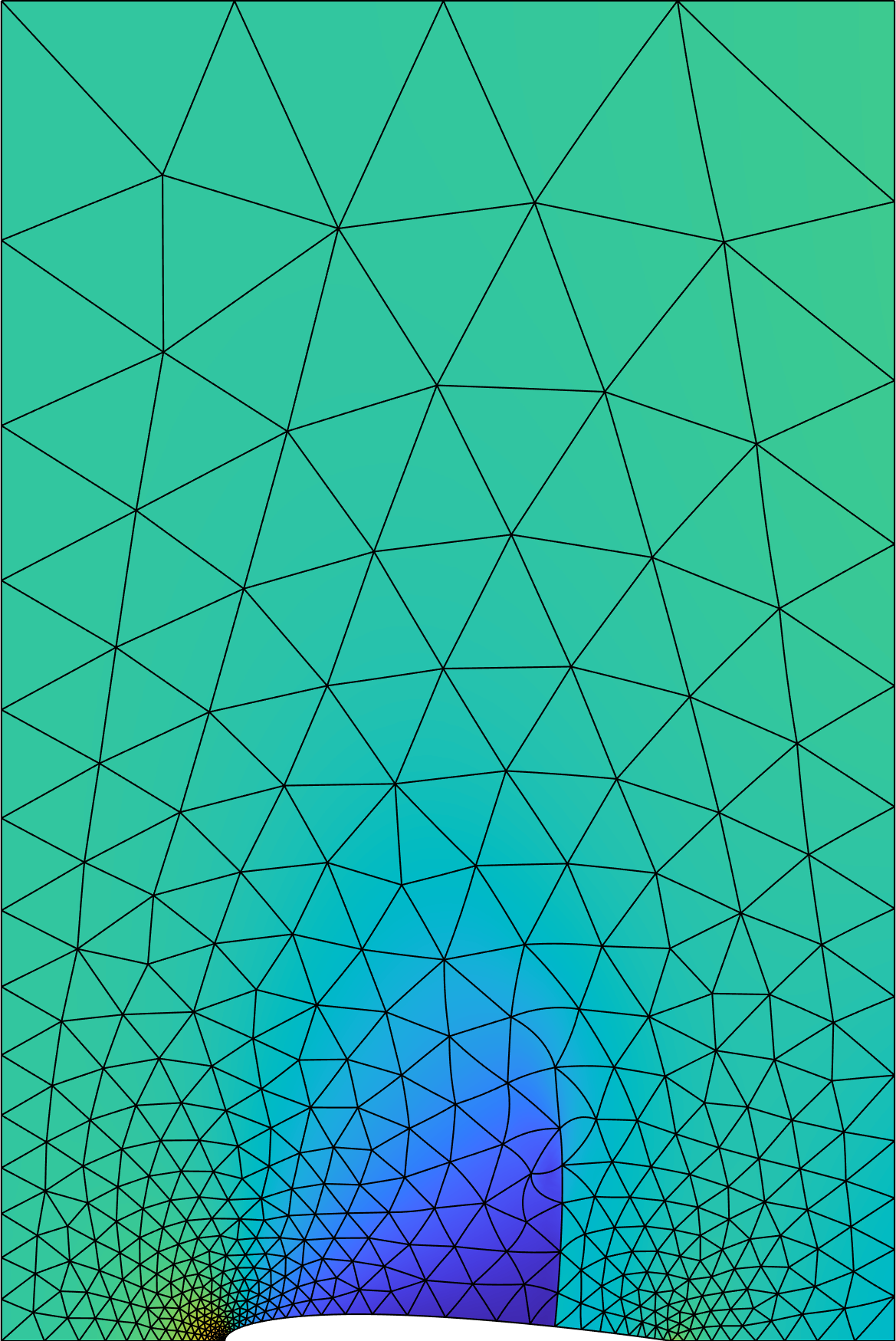};

\nextgroupplot[title={},xtick =\empty, ytick =\empty,
xlabel= {},ylabel= {}]
\addplot graphics [xmin=-0.5, xmax=1.5, ymin=0, ymax=3]{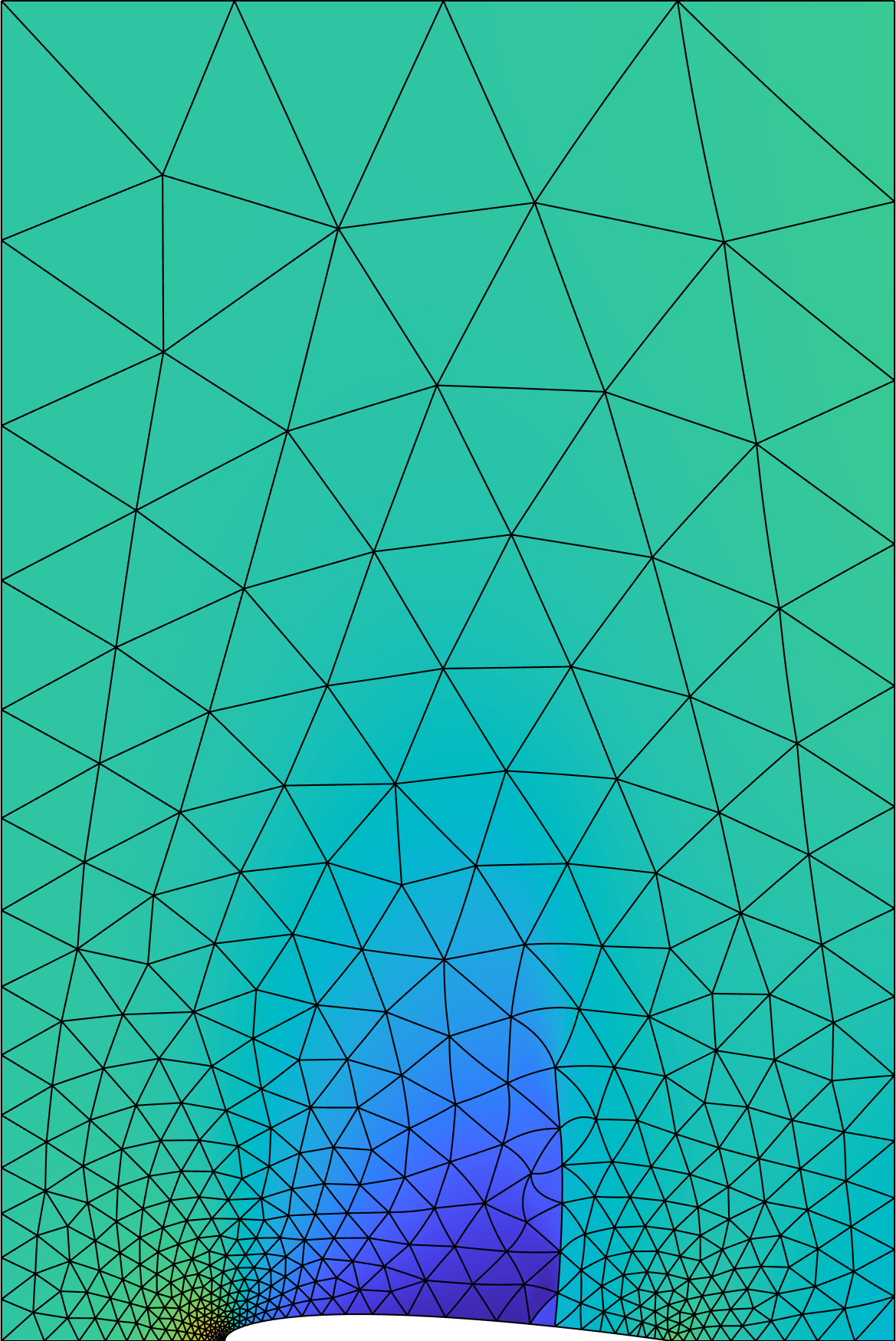};

\nextgroupplot[title={} ,xtick =\empty, ytick =\empty,
xlabel= {},ylabel= {}]
\addplot graphics [xmin=-0.5, xmax=1.5, ymin=0, ymax=3] {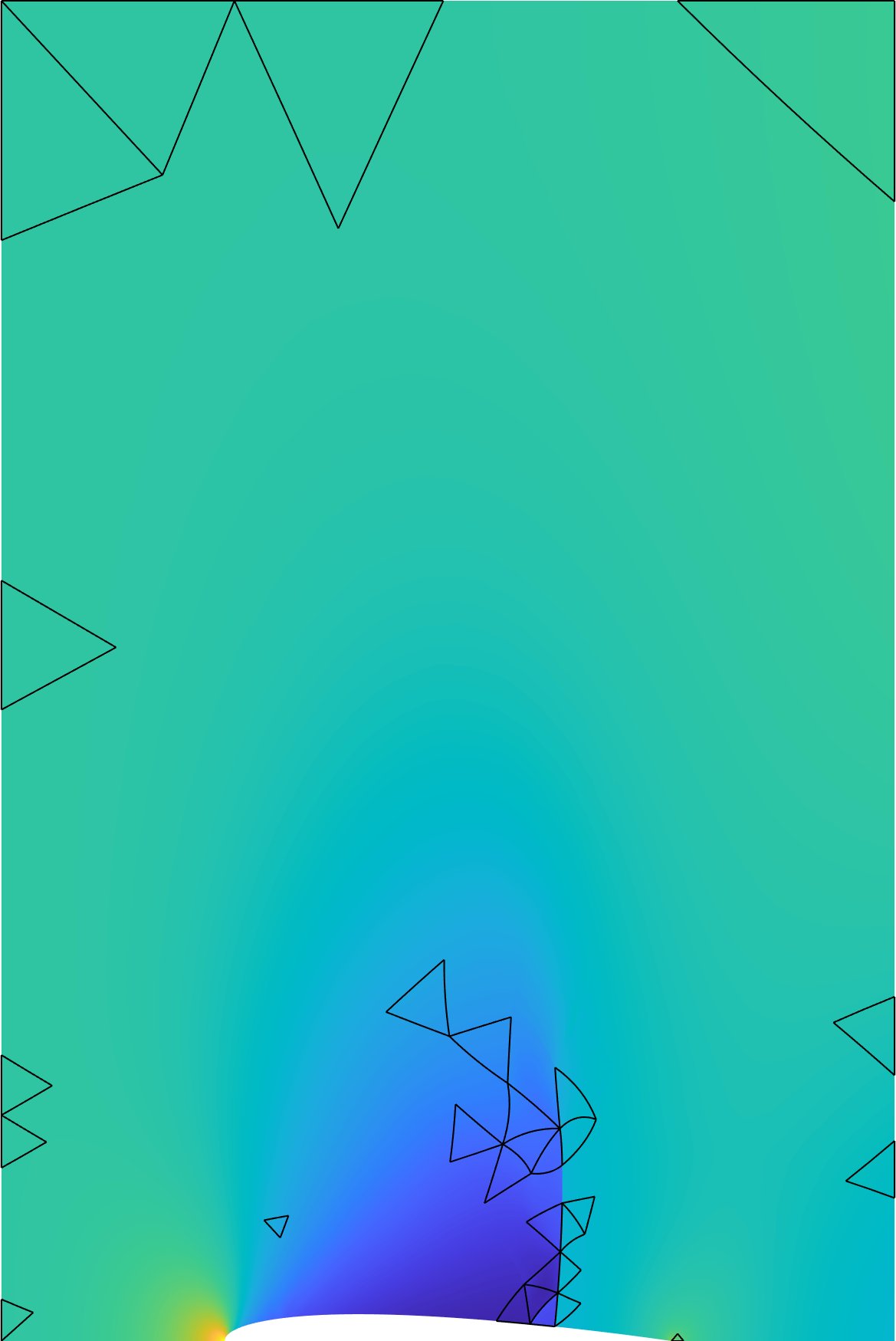};

\end{groupplot}
\end{tikzpicture}
\caption{Online solution of the transonic flow problem at $M_\infty = 0.78$
 in the physical domain (\textit{top}) and reference domain (\textit{bottom}). The mesh
 edges are included to show the elements over which the nonlinear terms must be assembled.}
\label{fig:naca_online}
\end{figure}

\begin{figure}
\centering
\begin{tikzpicture}
\begin{groupplot} [
group style={group size = 1 by 1, horizontal sep = 1cm}]
\nextgroupplot[width=0.6\textwidth, height=0.5\textwidth, ylabel={$\rho$},xlabel={$x$}]
\addplot [mark options={solid, thick}, mark=+, mark size=1.5, color=magenta]
coordinates {
(-5.00000000e-01,          nan)
(-4.59180000e-01,  1.39920000e+00)
(-4.18370000e-01,  1.40100000e+00)
(-3.77550000e-01,  1.40400000e+00)
(-3.36730000e-01,  1.40840000e+00)
(-2.95920000e-01,  1.41410000e+00)
(-2.55100000e-01,  1.42130000e+00)
(-2.14290000e-01,  1.43000000e+00)
(-1.73470000e-01,  1.43990000e+00)
(-1.32650000e-01,  1.44960000e+00)
(-9.18370000e-02,  1.45530000e+00)
(-5.10200000e-02,  1.44780000e+00)
(-1.02040000e-02,  1.41250000e+00)
( 3.06120000e-02,  1.34290000e+00)
( 7.14290000e-02,  1.25690000e+00)
( 1.12240000e-01,  1.17890000e+00)
( 1.53060000e-01,  1.11790000e+00)
( 1.93880000e-01,  1.07430000e+00)
( 2.34690000e-01,  1.03330000e+00)
( 2.75510000e-01,  1.00510000e+00)
( 3.16330000e-01,  9.82030000e-01)
( 3.57140000e-01,  9.63290000e-01)
( 3.97960000e-01,  9.47560000e-01)
( 4.38780000e-01,  9.34630000e-01)
( 4.79590000e-01,  9.24750000e-01)
( 5.20410000e-01,  9.15760000e-01)
( 5.61220000e-01,  9.12810000e-01)
( 6.02040000e-01,  9.03860000e-01)
( 6.42860000e-01,  9.07320000e-01)
( 6.83670000e-01,  1.30550000e+00)
( 7.24490000e-01,  1.30760000e+00)
( 7.65310000e-01,  1.31660000e+00)
( 8.06120000e-01,  1.32720000e+00)
( 8.46940000e-01,  1.33940000e+00)
( 8.87760000e-01,  1.35260000e+00)
( 9.28570000e-01,  1.36600000e+00)
( 9.69390000e-01,  1.37730000e+00)
( 1.01020000e+00,  1.38330000e+00)
( 1.05100000e+00,  1.38350000e+00)
( 1.09180000e+00,  1.37850000e+00)
( 1.13270000e+00,  1.37060000e+00)
( 1.17350000e+00,  1.36200000e+00)
( 1.21430000e+00,  1.35380000e+00)
( 1.25510000e+00,  1.34630000e+00)
( 1.29590000e+00,  1.33930000e+00)
( 1.33670000e+00,  1.33280000e+00)
( 1.37760000e+00,  1.32660000e+00)
( 1.41840000e+00,  1.32070000e+00)
( 1.45920000e+00,  1.31480000e+00)};\label{line:naca_rom}

\addplot [mark options={solid, thin}, mark=o, mark size=2, color=cyan]
coordinates {
(-5.00000000e-01,          nan)
(-4.59180000e-01,  1.39970000e+00)
(-4.18370000e-01,  1.40150000e+00)
(-3.77550000e-01,  1.40450000e+00)
(-3.36730000e-01,  1.40880000e+00)
(-2.95920000e-01,  1.41450000e+00)
(-2.55100000e-01,  1.42180000e+00)
(-2.14290000e-01,  1.43050000e+00)
(-1.73470000e-01,  1.44040000e+00)
(-1.32650000e-01,  1.45010000e+00)
(-9.18370000e-02,  1.45570000e+00)
(-5.10200000e-02,  1.44820000e+00)
(-1.02040000e-02,  1.41290000e+00)
( 3.06120000e-02,  1.34330000e+00)
( 7.14290000e-02,  1.25730000e+00)
( 1.12240000e-01,  1.17930000e+00)
( 1.53060000e-01,  1.11810000e+00)
( 1.93880000e-01,  1.07440000e+00)
( 2.34690000e-01,  1.03350000e+00)
( 2.75510000e-01,  1.00520000e+00)
( 3.16330000e-01,  9.82150000e-01)
( 3.57140000e-01,  9.63420000e-01)
( 3.97960000e-01,  9.47710000e-01)
( 4.38780000e-01,  9.34810000e-01)
( 4.79590000e-01,  9.24950000e-01)
( 5.20410000e-01,  9.15980000e-01)
( 5.61220000e-01,  9.13010000e-01)
( 6.02040000e-01,  9.04160000e-01)
( 6.42860000e-01,  9.07740000e-01)
( 6.83670000e-01,  1.30640000e+00)
( 7.24490000e-01,  1.30820000e+00)
( 7.65310000e-01,  1.31710000e+00)
( 8.06120000e-01,  1.32760000e+00)
( 8.46940000e-01,  1.33970000e+00)
( 8.87760000e-01,  1.35290000e+00)
( 9.28570000e-01,  1.36650000e+00)
( 9.69390000e-01,  1.37790000e+00)
( 1.01020000e+00,  1.38400000e+00)
( 1.05100000e+00,  1.38410000e+00)
( 1.09180000e+00,  1.37900000e+00)
( 1.13270000e+00,  1.37110000e+00)
( 1.17350000e+00,  1.36250000e+00)
( 1.21430000e+00,  1.35440000e+00)
( 1.25510000e+00,  1.34680000e+00)
( 1.29590000e+00,  1.33990000e+00)
( 1.33670000e+00,  1.33340000e+00)
( 1.37760000e+00,  1.32730000e+00)
( 1.41840000e+00,  1.32140000e+00)
( 1.45920000e+00,  1.31550000e+00)};\label{line:naca_hyper}

\addplot [solid, thick, color=black]
coordinates {
(-5.00000000e-01,          nan)
(-4.59180000e-01,  1.40060000e+00)
(-4.18370000e-01,  1.40230000e+00)
(-3.77550000e-01,  1.40520000e+00)
(-3.36730000e-01,  1.40950000e+00)
(-2.95920000e-01,  1.41510000e+00)
(-2.55100000e-01,  1.42230000e+00)
(-2.14290000e-01,  1.43100000e+00)
(-1.73470000e-01,  1.44090000e+00)
(-1.32650000e-01,  1.45080000e+00)
(-9.18370000e-02,  1.45670000e+00)
(-5.10200000e-02,  1.44980000e+00)
(-1.02040000e-02,  1.41530000e+00)
( 3.06120000e-02,  1.34670000e+00)
( 7.14290000e-02,  1.26130000e+00)
( 1.12240000e-01,  1.18240000e+00)
( 1.53060000e-01,  1.12000000e+00)
( 1.93880000e-01,  1.07500000e+00)
( 2.34690000e-01,  1.03630000e+00)
( 2.75510000e-01,  1.00660000e+00)
( 3.16330000e-01,  9.83800000e-01)
( 3.57140000e-01,  9.64820000e-01)
( 3.97960000e-01,  9.48410000e-01)
( 4.38780000e-01,  9.36020000e-01)
( 4.79590000e-01,  9.25220000e-01)
( 5.20410000e-01,  9.15810000e-01)
( 5.61220000e-01,  9.08760000e-01)
( 6.02040000e-01,  9.01830000e-01)
( 6.42860000e-01,  9.07990000e-01)
( 6.83670000e-01,  1.31650000e+00)
( 7.24490000e-01,  1.30980000e+00)
( 7.65310000e-01,  1.31310000e+00)
( 8.06120000e-01,  1.32120000e+00)
( 8.46940000e-01,  1.33320000e+00)
( 8.87760000e-01,  1.34770000e+00)
( 9.28570000e-01,  1.36370000e+00)
( 9.69390000e-01,  1.37860000e+00)
( 1.01020000e+00,  1.38770000e+00)
( 1.05100000e+00,  1.38820000e+00)
( 1.09180000e+00,  1.38210000e+00)
( 1.13270000e+00,  1.37370000e+00)
( 1.17350000e+00,  1.36530000e+00)
( 1.21430000e+00,  1.35740000e+00)
( 1.25510000e+00,  1.35020000e+00)
( 1.29590000e+00,  1.34360000e+00)
( 1.33670000e+00,  1.33750000e+00)
( 1.37760000e+00,  1.33170000e+00)
( 1.41840000e+00,  1.32610000e+00)
( 1.45920000e+00,  1.32050000e+00)};\label{line:naca_hdm}

\end{groupplot}\end{tikzpicture}
\caption{Online solution of the transonic flow problem at $M_\infty = 0.78$ along
 a horizontal slice $\{(s,0.14) \mid s\in(-0.5,1.5)\}$. Legend:
 HDM (\ref{line:naca_hdm}), ROM-IFT $j=3$ (\ref{line:naca_rom}), EQP-IFT $j=3$ (\ref{line:naca_hyper}).}
\label{fig:naca_slice}
\end{figure}
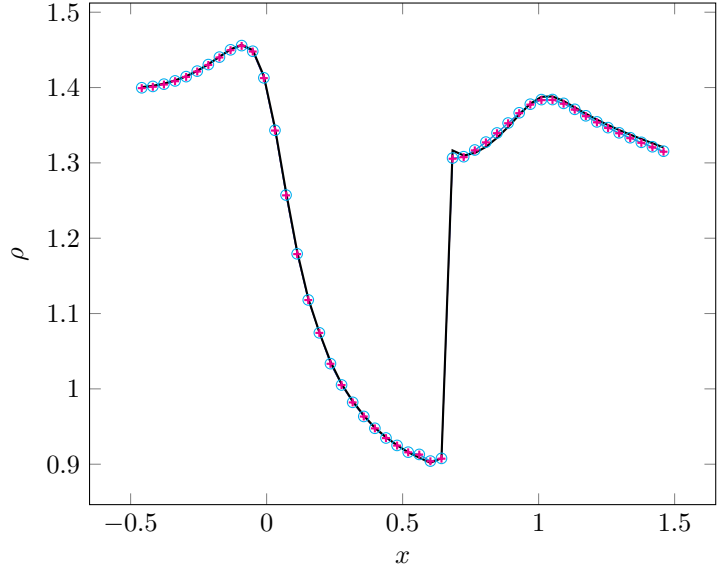

We use the greedy algorithm proposed in Section~\ref{sec:rom:greedy} with a
candidate set $\Dcal_\mathrm{cnd}$ with size $|\Dcal_\mathrm{cnd}|=9$ constructed
by uniformly sampling $\Dcal$ with three samples per dimension. We use this procedure
to build up a ROM-IFT and EQP-IFT of sizes $j=1,2,3$. For the EQP-IFT, we take the
EQP training set to be the candidate set ($\Xibold = \Dcal_\mathrm{cnd}$) and
the tolerance $\delta = 10^{-15}$. We test the performance of the reduced models
using the error metrics introduced in Section~\ref{sec:numexp} over a test set
$\Dcal_\mathrm{tst}$ with size $|\Dcal_\mathrm{tst}|=20$ constructed by uniformly
sampling $\Dcal$. The performance of the ROM-IFT and EQP-IFT methods is summarized in
Table~\ref{tab:naca}, which shows the both the offline and online cost of EQP-IFT
are reduced compared to ROM-IFT and the accuracy is comparable.
Figure~\ref{fig:naca-base} shows the HDM snapshot and domain deformation at the training
parameters selected by the greedy algorithm. Even at the parameter in the test set at
which the EQP-IFT error is largest ($\mubold=0.78$), the ROM/EQP-IFT provide accurate
approximations as can be seen from the state throughout the domain
(Figure~\ref{fig:naca_online}) and the state along a slice through the domain
(Figure~\ref{fig:naca_slice}).

\begin{table}[H]
\centering
\caption{Performance of HDM, ROM-IFT, EQP-IFT on test set $\Dcal_\mathrm{tst}$ for the transonic flow problem. Timings are reported as normalized by the cost to compute a single HDM solution (median CPU time for HDM = 291s).}
\label{tab:naca}
\begin{tabular}{r|ccccccc}
 & $j$ & $\Ncal$ & $E_j$  & $R_j$ & $T_j^\mathrm{gs}$ & $T_j^\rhobold$ & $T_j^\mathrm{tst}$ \\
\hline ROM-IFT&1&-&0.121&0.0537&-&-&0.032\\
&2&-&0.0366&0.0085&2.3&-&0.065\\
&3&-&0.0184&0.0058&3.4&-&0.099\\
\hline EQP-IFT&1&1.1\%&0.118&0.0537&-&0.32&0.0040\\
&2&3.6\%&0.0313&0.0104&0.51 &0.27&0.0059\\
&3&5.3\%&0.0215&0.0067&0.69&0.25&0.012
\end{tabular}
\end{table}

\begin{remark}
To ensure boundary conditions are enforced, elements with non-zero weights are
required on each boundary. In our experiments, the EQP training algorithm naturally
selected weights that satisfy this requirement. While this is the likely outcome
given the training formulation that imposes accuracy requirements on the optimality
criteria, it is not guaranteed. If necessary, the approach introduced in
\cite{washabaugh2016faster} can be used to force elements with non-zero
weights on each boundary or flow region.
\end{remark}

\section{Conclusion}
\label{sec:concl}
This work introduced an empirical quadrature-based hyperreduction procedure
and greedy training algorithm for model reduction with implicit feature
tracking \cite{mirhoseini2023model}. The ROM-IFT approach circumvents
the slowly decaying $n$-width limitation of standard model reduction methods
for convection-dominated problems by using a nonlinear approximation manifold
systematically defined by composing a low-dimensional affine space with a space
of bijections of the underlying domain. The ROM-IFT method is posed as a minimization
problem of the HDM residual to simultaneously determine the reduced state and
deformation coordinates. The hyperreduced version replaces
the elemental form of the ROM-IFT optimality system with a weighted
elemental form with a sparse weight vector determined via $\ell_1$
minimization, which leads to online efficiency as many elements in
the domain can be skipped. Furthermore, a greedy procedure is
developed to systematically train the ROM/EQP-IFT in the
parameter domain. The proposed EQP-IFT method preserves
the accuracy of the original ROM-IFT approach and further reduces
the computational cost by eliminating the bottleneck of assembling
the nonlinear terms.

\appendix

\section{Elemental derivations}
\label{app:elemental}

In this section, we derive the elemental representation of the reduced and
hyperreduced quantities that enable efficient computation in the case where
the element weight vectors is sparse. We begin with the assumption that the
residual $\Rbm$ can be written in the elemental form (\ref{eqn:hdm_elem})
and $\Pbm_e^T\Pbm_E = \delta_{eE}\Ibm$. Direct differentiation of
(\ref{eqn:hdm_elem}) leads to the Jacobians
\begin{equation} \label{eqn:hdm_elem_jac1}
\begin{aligned}
 \pder{\Rbm}{\Ubm}(\Ubm;\xbm,\mubold) &=
  \sum_{e=1}^{N_\mathtt{e}} \Pbm_e\left(\pder{\Rbm_e}{\Ubm_e}(\Ubm_e,\Ubm_e',\xbm,\mubold)\Pbm_e^T
   + \pder{\Rbm_e}{\Ubm_e'}(\Ubm_e,\Ubm_e',\xbm_e,\mubold)(\Pbm_e')^T\right) \\
 \pder{\Rbm}{\xbm}(\Ubm;\xbm,\mubold) &=
  \sum_{e=1}^{N_\mathtt{e}} \Pbm_e\pder{\Rbm_e}{\xbm_e}(\Ubm_e,\Ubm_e',\xbm_e,\mubold)\Qbm_e^T.
\end{aligned}
\end{equation}
Similarly, direct differentiation of the distortion residual leads to
\begin{equation} \label{eqn:eqn:hdm_elem_jac2}
 \pder{\Nbm}{\xbm}(\xbm) =
 \begin{bmatrix}
  \pder{\Nbm_1}{\xbm_1}(\xbm_1)\Qbm_1^T \\
  \vdots \\
  \pder{\Nbm_{N_\mathtt{e}}}{\xbm_{N_\mathtt{e}}}(\xbm_{N_\mathtt{e}})\Qbm_{N_\mathtt{e}}^T \\
 \end{bmatrix}
\end{equation}
In the remainder of this section, we drop arguments for brevity: the elemental
residual $\Rbm_e$ and its partial derivatives are evaluated at
$(\hat\Ubm_e,\hat\Ubm_e',\hat\xbm_e;\mubold)$ and the elemental distortions
$\Nbm_e$ and its partial derivatives are evaluated at $\hat\xbm_e$
($\hat\Ubm_e$, $\hat\Ubm_e'$, $\hat\xbm_e$ are defined at the end of
Section~\ref{sec:rom:ift}).
 
\subsection{Reduced residuals}
We expand the reduced residuals in (\ref{eqn:ift-res}) to their elemental
form in (\ref{eqn:ift-res-elem1})-(\ref{eqn:ift-res-elem2}). The state residual
is expanded as
\begin{equation}
\begin{aligned}
 \hat\Sbm_{\Phibold,\Psibold}(\wbm,\taubold;\mubold)
   &= \Phibold^T\left(\sum_{e=1}^{N_\mathtt{e}}\Pbm_e\left(\pder{\Rbm_e}{\Ubm_e}\Pbm_e^T + \pder{\Rbm_e}{\Ubm_e'}(\Pbm_e')^T\right)\right)^T\left(\sum_{E=1}^{N_\mathtt{e}} \Pbm_E\Rbm_E\right) \\
   &= \Phibold^T\sum_{e,E=1}^{N_\mathtt{e}} \left(\Pbm_e\pder{\Rbm_e}{\Ubm_e}^T+\Pbm_e'\pder{\Rbm_e}{\Ubm_e'}^T\right)(\Pbm_e^T\Pbm_E)\Rbm_E \\
   &= \sum_{e=1}^{N_\mathtt{e}} \left(\Phibold_e^T\pder{\Rbm_e}{\Ubm_e}^T + (\Phibold_e')^T\pder{\Rbm_e}{\Ubm_e'}^T\right)\Rbm_e,
\end{aligned}
\end{equation}
where the elemental form of the residual (\ref{eqn:hdm_elem}) and
Jacobian (\ref{eqn:hdm_elem_jac1}), as well as the
$\Pbm_e^T\Pbm_E = \delta_{eE}\Ibm$ property were used.
An identical procedure leads to the expansion of the deformation residual
\begin{equation}
\begin{aligned}
 \hat\Tbm_{\Phibold,\Psibold}(\wbm,\taubold;\mubold)
   &= \Psibold^T\left(\sum_{e=1}^{N_\mathtt{e}}\Pbm_e\pder{\Rbm_e}{\xbm_e}\Qbm_e^T\right)^T\left(\sum_{E=1}^{N_\mathtt{e}} \Pbm_E\Rbm_E\right) + \Psibold^T\sum_{e=1}^{N_\mathtt{e}} \left(\pder{\Nbm_e}{\xbm_e}\Qbm_e^T\right)^T\Nbm_e\\
   &= \Psibold^T\sum_{e,E=1}^{N_\mathtt{e}} \left(\Qbm_e\pder{\Rbm_e}{\xbm_e}^T\right)(\Pbm_e^T\Pbm_E)\Rbm_E + \Psibold^T\sum_{e=1}^{N_\mathtt{e}} \Qbm_e\pder{\Nbm_e}{\xbm_e}^T\Nbm_e \\
   &= \sum_{e=1}^{N_\mathtt{e}} \Psibold_e^T\left(\pder{\Rbm_e}{\xbm_e}^T\Rbm_e + \pder{\Nbm_e}{\xbm_e}^T\Nbm_e\right).
\end{aligned}
\end{equation}

\subsection{Hyperreduced optimality system}
To derive the hyperreduced residuals, we form the first-order optimality system
of the optimization problem in (\ref{eqn:ift-opt-eqp0}), which leads to
\begin{equation} \label{eqn:ift-opt-eqp00}
 \pder{\tilde\Fbm_{\Phibold,\Psibold,\rhobold}}{\wbm}(\wbm,\taubold;\mubold)^T\tilde\Fbm_{\Phibold,\Psibold,\rhobold}(\wbm,\taubold;\mubold) = \zerobold, \qquad
 \pder{\tilde\Fbm_{\Phibold,\Psibold,\rhobold}}{\taubold}(\wbm,\taubold;\mubold)^T\tilde\Fbm_{\Phibold,\Psibold,\rhobold}(\wbm,\taubold;\mubold) = \zerobold.
\end{equation}
Expanding the first-order optimality system (\ref{eqn:ift-opt-eqp00}) for the
state using the definition of $\tilde\Fbm_{\Phibold,\Psibold}$ in
(\ref{eqn:ift-eqp-res}) leads to
\begin{equation}
\begin{aligned}
 \pder{\tilde\Fbm_{\Phibold,\Psibold,\rhobold}}{\wbm}(\wbm,\taubold;\mubold)^T\tilde\Fbm_{\Phibold,\Psibold,\rhobold}(\wbm,\taubold;\mubold)
 &= \pder{\tilde\Rbm_{\Phibold,\Psibold,\rhobold}}{\wbm}(\wbm,\taubold;\mubold)^T\tilde\Rbm_{\Phibold,\Psibold,\rhobold}(\wbm,\taubold;\mubold) \\
 &= \left(\sum_{e=1}^{N_\mathtt{e}}\sqrt{\rho_e}\Pbm_e\left(\pder{\Rbm_e}{\Ubm_e}\Phibold_e +\pder{\Rbm_e}{\Ubm_e'}\Phibold_e'\right)\right)^T\left(\sum_{E=1}^{N_\mathtt{e}}\sqrt{\rho_E}\Pbm_E\Rbm_E\right) \\
 &= \sum_{e,E=1}^{N_\mathtt{e}} \sqrt{\rho_e\rho_E}\left(\Phibold_e^T\pder{\Rbm_e}{\Ubm_e}^T+(\Phibold_e')^T\pder{\Rbm_e}{\Ubm_e'}^T\right)(\Pbm_e^T\Pbm_E)\Rbm_E \\
 &= \sum_{e=1}^{N_\mathtt{e}} \rho_e \left(\Phibold_e^T\pder{\Rbm_e}{\Ubm_e}^T+(\Phibold_e')^T\pder{\Rbm_e}{\Ubm_e'}^T\right) \Rbm_e \\
 &= \tilde{\Sbm}_{\Phibold,\Psibold,\rhobold}(\wbm,\taubold;\mubold),
\end{aligned}
\end{equation}
where the elemental form of the residual (\ref{eqn:hdm_elem})
and Jacobian (\ref{eqn:hdm_elem_jac1}), the $\Pbm_e^T\Pbm_E = \delta_{eE}\Ibm$
property, and the definition of $\tilde\Sbm_{\Phibold,\Psibold}$ in
(\ref{eqn:ift-eqp-res0}) were used. A similar procedure is used to
derive the elemental form of the hyperreduced deformation residual
\begin{equation}
\begin{aligned}
 \pder{\tilde\Fbm_{\Phibold,\Psibold,\rhobold}}{\taubold}(\wbm,\taubold;\mubold)^T\tilde\Fbm_{\Phibold,\Psibold,\rhobold}(\wbm,\taubold;\mubold)
 &= \pder{\tilde\Rbm_{\Phibold,\Psibold,\rhobold}}{\taubold}(\wbm,\taubold;\mubold)^T\tilde\Rbm_{\Phibold,\Psibold,\rhobold}(\wbm,\taubold;\mubold) + \pder{\tilde\Nbm_{\Psibold,\rhobold}}{\taubold}(\taubold)^T\tilde\Nbm_{\Psibold,\rhobold}(\taubold) \\
 &= \left(\sum_{e=1}^{N_\mathtt{e}}\sqrt{\rho_e}\Pbm_e\pder{\Rbm_e}{\xbm_e}\Psibold_e\right)^T\left(\sum_{E=1}^{N_\mathtt{e}}\sqrt{\rho_E}\Pbm_E\Rbm_E\right) + \sum_{e=1}^{N_\mathtt{e}} \rho_e\Psibold_e^T\pder{\Nbm_e}{\xbm_e}^T \Nbm_e \\
 &= \sum_{e,E=1}^{N_\mathtt{e}} \sqrt{\rho_e\rho_E}\Psibold_e^T\pder{\Rbm_e}{\xbm_e}^T(\Pbm_e^T\Pbm_E)\Rbm_E + \sum_{e=1}^{N_\mathtt{e}} \rho_e\Psibold_e^T\pder{\Nbm_e}{\xbm_e}^T \Nbm_e\\
 &= \sum_{e=1}^{N_\mathtt{e}} \rho_e \Psibold_e^T\left(\pder{\Rbm_e}{\xbm_e}^T\Rbm_e + \pder{\Nbm_e}{\xbm_e}^T \Nbm_e\right) \\
 &= \tilde{\Tbm}_{\Phibold,\Psibold,\rhobold}(\wbm,\taubold;\mubold).
\end{aligned}
\end{equation}
Next, we show the objective function $\tilde{J}_{\Phibold,\Psibold,\rhobold}$ in
(\ref{eqn:ift-opt-eqp1}) can be expanded in the elemental form in
(\ref{eqn:ift-opt-eqp2})
\begin{equation}
\begin{aligned}
 2\tilde{J}_{\Phibold,\Psibold,\rhobold}(\wbm,\taubold;\mubold)
 &= \norm{\tilde\Fbm_{\Phibold,\Psibold,\rhobold}(\wbm,\taubold;\mubold)}_2^2 \\
 &= \norm{\tilde\Rbm_{\Phibold,\Psibold,\rhobold}(\wbm,\taubold;\mubold)}_2^2 +
    \norm{\tilde\Nbm_{\Psibold,\rhobold}(\taubold)}_2^2 \\
 &= \left(\sum_{e=1}^{N_\mathtt{e}}\sqrt{\rho_e}\Pbm_e\Rbm_e\right)^T
    \left(\sum_{E=1}^{N_\mathtt{e}}\sqrt{\rho_E}\Pbm_E\Rbm_E\right) +
    \sum_{e=1}^{N_\mathtt{e}} \rho_e \Nbm_e^2 \\
 &= \sum_{e,E=1}^{N_\mathtt{e}} \sqrt{\rho_e\rho_E}\Rbm_e^T(\Pbm_e^T\Pbm_E)\Rbm_E
    + \sum_{e=1}^{N_\mathtt{e}} \rho_e \Nbm_e^2 \\
 &= \sum_{e=1}^{N_\mathtt{e}} \rho_e\Rbm_e^T\Rbm_e + \sum_{e=1}^{N_\mathtt{e}} \rho_e \Nbm_e^2 \\
 &= \sum_{e=1}^{N_\mathtt{e}} \rho_e\norm{\Rbm_e}_2^2 + \sum_{e=1}^{N_\mathtt{e}} \rho_e |\Nbm_e|^2,
\end{aligned}
\end{equation}
where the first two equalities follows from the definition of
$\tilde{J}_{\Phibold,\Psibold}$ and $\tilde{F}_{\Phibold,\Psibold}$,
the third equality follows from the elemental form of the residual
(\ref{eqn:hdm_elem}) and distortion (\ref{eqn:dist_elem0}), and the
fifth equality follows from $\Pbm_e^T\Pbm_E = \delta_{eE}\Ibm$.

\subsection{Hyperreduced Hessian approximation}
Finally, we derive elemental expressions for hyperreduced Hessian
approximations used for the Levenberg-Marquardt solver in (\ref{eqn:levmarq-elem}).
For brevity, we define individual components of the Hessian approximation
\begin{equation}
\begin{aligned}
 \tilde\Hbm_{\Phibold,\Psibold,\rhobold}^{\wbm\wbm} : (\wbm,\taubold;\mubold) &\mapsto
  \pder{\tilde\Fbm_{\Phibold,\Psibold,\rhobold}}{\wbm}(\wbm,\taubold;\mubold)^T
  \pder{\tilde\Fbm_{\Phibold,\Psibold,\rhobold}}{\wbm}(\wbm,\taubold;\mubold) \\
 \tilde\Hbm_{\Phibold,\Psibold,\rhobold}^{\wbm\taubold} : (\wbm,\taubold;\mubold) &\mapsto
  \pder{\tilde\Fbm_{\Phibold,\Psibold,\rhobold}}{\wbm}(\wbm,\taubold;\mubold)^T
  \pder{\tilde\Fbm_{\Phibold,\Psibold,\rhobold}}{\taubold}(\wbm,\taubold;\mubold) \\
 \tilde\Hbm_{\Phibold,\Psibold,\rhobold}^{\taubold\taubold} : (\wbm,\taubold;\mubold) &\mapsto
  \pder{\tilde\Fbm_{\Phibold,\Psibold,\rhobold}}{\taubold}(\wbm,\taubold;\mubold)^T
  \pder{\tilde\Fbm_{\Phibold,\Psibold,\rhobold}}{\taubold}(\wbm,\taubold;\mubold).
\end{aligned}
\end{equation}
Expanding $\Hbm_{\Phibold,\Psibold,\rhobold}^{\wbm\wbm}$ in its elemental
contribution yields
\begin{equation}
\begin{aligned}
 \tilde\Hbm_{\Phibold,\Psibold,\rhobold}^{\wbm\wbm}(\wbm,\taubold;\mubold)
  &= \pder{\tilde\Rbm_{\Phibold,\Psibold,\rhobold}}{\wbm}(\wbm,\taubold;\mubold)^T
  \pder{\tilde\Rbm_{\Phibold,\Psibold,\rhobold}}{\wbm}(\wbm,\taubold;\mubold) \\
   &= \left(\sum_{e=1}^{N_\mathtt{e}}\sqrt{\rho_e}\Pbm_e\left(\pder{\Rbm_e}{\Ubm_e}\Phibold_e +\pder{\Rbm_e}{\Ubm_e'}\Phibold_e'\right)\right)^T\left(\sum_{E=1}^{N_\mathtt{e}}\sqrt{\rho_E}\Pbm_E\left(\pder{\Rbm_E}{\Ubm_E}\Phibold_E+\pder{\Rbm_E}{\Ubm_E'}\Phibold_E'\right)\right) \\
   &= \sum_{e,E=1}^{N_\mathtt{e}}\sqrt{\rho_e\rho_E}\left(\pder{\Rbm_e}{\Ubm_e}\Phibold_e +\pder{\Rbm_e}{\Ubm_e'}\Phibold_e'\right)^T(\Pbm_e^T\Pbm_E)\left(\pder{\Rbm_E}{\Ubm_E}\Phibold_E+\pder{\Rbm_E}{\Ubm_E'}\Phibold_E'\right) \\
   &= \sum_{e=1}^{N_\mathtt{e}}\rho_e\left(\pder{\Rbm_e}{\Ubm_e}\Phibold_e +\pder{\Rbm_e}{\Ubm_e'}\Phibold_e'\right)^T\left(\pder{\Rbm_e}{\Ubm_e}\Phibold_e+\pder{\Rbm_e}{\Ubm_e'}\Phibold_e'\right).
\end{aligned}
\end{equation}
The steps of the derivation follow the procedure in the previous sections, i.e.,
use the elemental expansion of the Jacobian (\ref{eqn:hdm_elem_jac1}) and
the $\Pbm_e^T\Pbm_E = \delta_{eE}\Ibm$ property.
Expanding $\Hbm_{\Phibold,\Psibold,\rhobold}^{\wbm\taubold}$ in its elemental
contribution yields
\begin{equation}
\begin{aligned}
 \tilde\Hbm_{\Phibold,\Psibold,\rhobold}^{\wbm\taubold}(\wbm,\taubold;\mubold)
  &= \pder{\tilde\Rbm_{\Phibold,\Psibold,\rhobold}}{\wbm}(\wbm,\taubold;\mubold)^T
  \pder{\tilde\Rbm_{\Phibold,\Psibold,\rhobold}}{\taubold}(\wbm,\taubold;\mubold) \\
   &= \left(\sum_{e=1}^{N_\mathtt{e}}\sqrt{\rho_e}\Pbm_e\left(\pder{\Rbm_e}{\Ubm_e}\Phibold_e +\pder{\Rbm_e}{\Ubm_e'}\Phibold_e'\right)\right)^T\left(\sum_{E=1}^{N_\mathtt{e}}\sqrt{\rho_E}\Pbm_E\pder{\Rbm_E}{\xbm_E}\Psibold_E\right) \\
   &= \sum_{e,E=1}^{N_\mathtt{e}}\sqrt{\rho_e\rho_E}\left(\pder{\Rbm_e}{\Ubm_e}\Phibold_e +\pder{\Rbm_e}{\Ubm_e'}\Phibold_e'\right)^T(\Pbm_e^T\Pbm_E)\pder{\Rbm_E}{\xbm_E}\Psibold_E \\
   &= \sum_{e=1}^{N_\mathtt{e}}\rho_e\left(\pder{\Rbm_e}{\Ubm_e}\Phibold_e +\pder{\Rbm_e}{\Ubm_e'}\Phibold_e'\right)^T\pder{\Rbm_e}{\xbm_e}\Psibold_e.
\end{aligned}
\end{equation}
Expanding $\Hbm_{\Phibold,\Psibold,\rhobold}^{\taubold\taubold}$ in its elemental
contribution yields
\begin{equation}
\begin{aligned}
 \tilde\Hbm_{\Phibold,\Psibold,\rhobold}^{\taubold\taubold}(\wbm,\taubold;\mubold)
   &= \pder{\tilde\Rbm_{\Phibold,\Psibold,\rhobold}}{\taubold}(\wbm,\taubold;\mubold)^T
      \pder{\tilde\Rbm_{\Phibold,\Psibold,\rhobold}}{\taubold}(\wbm,\taubold;\mubold) +
      \pder{\tilde\Nbm_{\Psibold,\rhobold}}{\taubold}(\taubold)^T
      \pder{\tilde\Nbm_{\Psibold,\rhobold}}{\taubold}(\taubold) \\
   &= \left(\sum_{e=1}^{N_\mathtt{e}}\sqrt{\rho_e}\Pbm_e\pder{\Rbm_e}{\xbm_e}\Psibold_e\right)^T\left(\sum_{E=1}^{N_\mathtt{e}}\sqrt{\rho_E}\Pbm_E\pder{\Rbm_E}{\xbm_E}\Psibold_E\right) + \sum_{e=1}^{N_\mathtt{e}} \rho_e \Psibold_e^T\pder{\Nbm_e}{\xbm_e}\pder{\Nbm_e}{\xbm_e}\Psibold_e \\
   &= \sum_{e,E=1}^{N_\mathtt{e}}\sqrt{\rho_e\rho_E}\left(\pder{\Rbm_e}{\xbm_e}\Psibold_e\right)^T(\Pbm_e^T\Pbm_E)\pder{\Rbm_E}{\xbm_E}\Psibold_E + \sum_{e=1}^{N_\mathtt{e}} \rho_e \Psibold_e^T\pder{\Nbm_e}{\xbm_e}\pder{\Nbm_e}{\xbm_e}\Psibold_e\\
   &= \sum_{e=1}^{N_\mathtt{e}}\rho_e\Psibold_e^T\left(\pder{\Rbm_e}{\xbm_e}^T\pder{\Rbm_e}{\xbm_e} + \pder{\Nbm_e}{\xbm_e}^T\pder{\Nbm_e}{\xbm_e}\right)\Psibold_e.
\end{aligned}
\end{equation}

\section*{Acknowledgments}
This work is supported by AFOSR award numbers FA9550-20-1-0236, 
FA9550-22-1-0002, FA9550-22-1-0004, and ONR award number 
N00014-22-1-2299. The content of this publication does not necessarily 
reflect the position or policy of any of these supporters, and no official 
endorsement should be inferred.

\bibliographystyle{plain}
\bibliography{biblio}

\begin{thebibliography}{10}

\bibitem{agarwal_trust-region_2013}
Anshul Agarwal and Lorenz~T. Biegler.
\newblock A trust-region framework for constrained optimization using reduced
  order modeling.
\newblock {\em Optimization and Engineering}, 14(1):3--35, March 2013.

\bibitem{amsallem_nonlinear_2012}
David Amsallem, Matthew~J. Zahr, and Charbel Farhat.
\newblock Nonlinear model order reduction based on local reduced-order bases.
\newblock {\em International Journal for Numerical Methods in Engineering},
  92(10):891--916, December 2012.

\bibitem{arian_trust-region_2000}
E.~Arian, M.~Fahl, and E.~W. Sachs.
\newblock Trust-{Region} {Proper} {Orthogonal} {Decomposition} for {Flow}
  {Control}.
\newblock Technical Report ICASE-2000-25, Institue for Computer Applications in
  Science and Engineering, May 2000.

\bibitem{astrid2008missing}
Patricia Astrid, Siep Weiland, Karen Willcox, and Ton Backx.
\newblock Missing point estimation in models described by proper orthogonal
  decomposition.
\newblock {\em IEEE Transactions on Automatic Control}, 53(10):2237--2251,
  2008.

\bibitem{bansal2021model}
Harshit Bansal, Stephan Rave, Laura Iapichino, W~Schilders, and N~Wouw.
\newblock Model order reduction framework for problems with moving
  discontinuities.
\newblock In {\em Numerical Mathematics and Advanced Applications ENUMATH
  2019}, pages 83--91. Springer, 2021.

\bibitem{barrault2004empirical}
Maxime Barrault, Yvon Maday, Ngoc~Cuong Nguyen, and Anthony~T Patera.
\newblock An `empirical interpolation'method: application to efficient
  reduced-basis discretization of partial differential equations.
\newblock {\em Comptes Rendus Mathematique}, 339(9):667--672, 2004.

\bibitem{black_projection-based_2020}
Felix Black, Philipp Schulze, and Benjamin Unger.
\newblock Projection-based model reduction with dynamically transformed modes.
\newblock {\em ESAIM: Mathematical Modelling and Numerical Analysis},
  54(6):2011--2043, 2020.

\bibitem{cagniart_model_2019}
Nicolas Cagniart, Yvon Maday, and Benjamin Stamm.
\newblock Model order reduction for problems with large convection effects.
\newblock In B.~N. Chetverushkin, W.~Fitzgibbon, Y.A. Kuznetsov,
  P.~Neittaanm{\"a}ki, J.~Periaux, and O.~Pironneau, editors, {\em
  Contributions to {Partial} {Differential} {Equations} and {Applications}},
  Computational {Methods} in {Applied} {Sciences}, pages 131--150. Springer
  International Publishing, Cham, 2019.

\bibitem{carlberg_adaptive_2015}
Kevin Carlberg.
\newblock Adaptive h-refinement for reduced-order models.
\newblock {\em International Journal for Numerical Methods in Engineering},
  102(5):1192--1210, 2015.

\bibitem{chapman2017accelerated}
Todd Chapman, Philip Avery, Pat Collins, and Charbel Farhat.
\newblock Accelerated mesh sampling for the hyper reduction of nonlinear
  computational models.
\newblock {\em International Journal for Numerical Methods in Engineering},
  109(12):1623--1654, 2017.

\bibitem{chaturantabut2010nonlinear}
Saifon Chaturantabut and Danny Sorensen.
\newblock Nonlinear model reduction via discrete empirical interpolation.
\newblock {\em SIAM Journal on Scientific Computing}, 32(5):2737--2764, 2010.

\bibitem{constantine_reduced_2012}
Paul Constantine and Gianluca Iaccarino.
\newblock Reduced order models for parameterized hyperbolic conservations laws
  with shock reconstruction.
\newblock Technical report, Center for Turbulence Research Annual Brief, 2012.

\bibitem{corrigan_moving_2019}
Andrew Corrigan, Andrew~D. Kercher, and David~A. Kessler.
\newblock A moving discontinuous {Galerkin} finite element method for flows
  with interfaces: {A} moving discontinuous {Galerkin} finite element method
  for flows with interfaces.
\newblock {\em International Journal for Numerical Methods in Fluids},
  89(9):362--406, March 2019.

\bibitem{dihlmann_model_2011}
Markus Dihlmann, Martin Drohmann, and Bernard Haasdonk.
\newblock Model reduction of parametrized evolution problems using the reduced
  basis method with adaptive time-partitioning.
\newblock In {\em {Proceedings} of {ADMOS}}, page~64, 2011.

\bibitem{farhat2014dimensional}
Charbel Farhat, Philip Avery, Todd Chapman, and Julien Cortial.
\newblock Dimensional reduction of nonlinear finite element dynamic models with
  finite rotations and energy-based mesh sampling and weighting for
  computational efficiency.
\newblock {\em International Journal for Numerical Methods in Engineering},
  98(9):625--662, 2014.

\bibitem{farhat2015structure}
Charbel Farhat, Todd Chapman, and Philip Avery.
\newblock Structure-preserving, stability, and accuracy properties of the
  energy-conserving sampling and weighting method for the hyper reduction of
  nonlinear finite element dynamic models.
\newblock {\em International Journal for Numerical Methods in Engineering},
  102(5):1077--1110, 2015.

\bibitem{ferrero2022registration}
Andrea Ferrero, Tommaso Taddei, and Lei Zhang.
\newblock Registration-based model reduction of parameterized two-dimensional
  conservation laws.
\newblock {\em Journal of Computational Physics}, 457:111068, 2022.

\bibitem{grepl2005posteriori}
Martin~A Grepl and Anthony~T Patera.
\newblock A posteriori error bounds for reduced-basis approximations of
  parametrized parabolic partial differential equations.
\newblock {\em ESAIM: Mathematical Modelling and Numerical Analysis},
  39(1):157--181, 2005.

\bibitem{haasdonk_training_2011}
Bernard Haasdonk, Markus Dihlmann, and Mario Ohlberger.
\newblock A training set and multiple bases generation approach for
  parameterized model reduction based on adaptive grids in parameter space.
\newblock {\em Mathematical and Computer Modelling of Dynamical Systems},
  17(4):423--442, August 2011.

\bibitem{hartman_deep_2017}
David Hartman and Lalit~K. Mestha.
\newblock A deep learning framework for model reduction of dynamical systems.
\newblock In {\em 2017 {IEEE} {Conference} on {Control} {Technology} and
  {Applications} ({CCTA})}, pages 1917--1922, August 2017.

\bibitem{hesthaven2007nodal}
Jan Hesthaven and Tim Warburton.
\newblock {\em Nodal discontinuous Galerkin methods: algorithms, analysis, and
  applications}.
\newblock Springer Science \& Business Media, 2007.

\bibitem{huang2023predictive}
Cheng Huang and Karthik Duraisamy.
\newblock Predictive reduced order modeling of chaotic multi-scale problems
  using adaptively sampled projections.
\newblock {\em arXiv preprint arXiv:2301.09006}, 2023.

\bibitem{huang2022robust}
Tianci Huang and Matthew~J. Zahr.
\newblock A robust, high-order implicit shock tracking method for simulation of
  complex, high-speed flows.
\newblock {\em Journal of Computational Physics}, 454:110981, 2022.

\bibitem{kashima_nonlinear_2016}
Kenji Kashima.
\newblock Nonlinear model reduction by deep autoencoder of noise response data.
\newblock In {\em 2016 {IEEE} 55th {Conference} on {Decision} and {Control}
  ({CDC})}, pages 5750--5755, Las Vegas, NV, USA, December 2016. IEEE.

\bibitem{kim2020efficient}
Youngkyu Kim, Youngsoo Choi, David Widemann, and Tarek Zohdi.
\newblock Efficient nonlinear manifold reduced order model.
\newblock {\em arXiv preprint arXiv:2011.07727}, 2020.

\bibitem{lee_model_2020}
Kookjin Lee and Kevin~T. Carlberg.
\newblock Model reduction of dynamical systems on nonlinear manifolds using
  deep convolutional autoencoders.
\newblock {\em Journal of Computational Physics}, 404:108973, March 2020.

\bibitem{legresley2006application}
Patrick~Allen LeGresley.
\newblock {\em Application of Proper Orthogonal Decomposition ({POD}) to Design
  Decomposition Methods}.
\newblock Stanford University, 2006.

\bibitem{lucia_reduced_2003}
David~J. Lucia, Paul~I. King, and Philip~S. Beran.
\newblock Reduced order modeling of a two-dimensional flow with moving shocks.
\newblock {\em Computers \& Fluids}, 32(7):917--938, August 2003.

\bibitem{maday2002priori}
Yvon Maday, Anthony~T Patera, and Gabriel Turinici.
\newblock A priori convergence theory for reduced-basis approximations of
  single-parameter elliptic partial differential equations.
\newblock {\em Journal of Scientific Computing}, 17:437--446, 2002.

\bibitem{maulik_reduced-order_2020}
Romit Maulik, Bethany Lusch, and Prasanna Balaprakash.
\newblock Reduced-order modeling of advection-dominated systems with recurrent
  neural networks and convolutional autoencoders.
\newblock {\em Physics of Fluids}, 33(3):037106, 2021.

\bibitem{mirhoseini2023model}
Marzieh~Alireza Mirhoseini and Matthew~J. Zahr.
\newblock Model reduction of convection-dominated partial differential
  equations via optimization-based implicit feature tracking.
\newblock {\em Journal of Computational Physics}, 473:111739, 2023.

\bibitem{mojgani_physics-aware_2020}
Rambod Mojgani and Maciej Balajewicz.
\newblock Physics-aware registration based auto-encoder for convection
  dominated {PDEs}.
\newblock {\em arXiv:2006.15655 [cs, math]}, June 2020.
\newblock arXiv: 2006.15655.

\bibitem{nair_transported_2019}
Nirmal~J. Nair and Maciej Balajewicz.
\newblock Transported snapshot model order reduction approach for parametric,
  steady-state fluid flows containing parameter-dependent shocks: {Model} order
  reduction for fluid flows containing shocks.
\newblock {\em International Journal for Numerical Methods in Engineering},
  117(12):1234--1262, March 2019.

\bibitem{nocedal_numerical_2006}
Jorge Nocedal and Stephen~J. Wright.
\newblock {\em Numerical Optimization}.
\newblock Springer Series in Operations Research. Springer, New York, 2006.

\bibitem{ohlberger_nonlinear_2013}
Mario Ohlberger and Stephan Rave.
\newblock Nonlinear reduced basis approximation of parameterized evolution
  equations via the method of freezing.
\newblock {\em Comptes Rendus Mathematique}, 351(23-24):901--906, December
  2013.

\bibitem{ohlberger_reduced_2016}
Mario Ohlberger and Stephan Rave.
\newblock Reduced {Basis} {Methods}: {Success}, {Limitations} and {Future}
  {Challenges}.
\newblock In {\em {Proceedings} of the {Conference} {Algoritmy}}, 2016.

\bibitem{omata_novel_2019}
Noriyasu Omata and Susumu Shirayama.
\newblock A novel method of low-dimensional representation for temporal
  behavior of flow fields using deep autoencoder.
\newblock {\em AIP Advances}, 9(1):015006, January 2019.

\bibitem{peherstorfer_model_2020}
Benjamin Peherstorfer.
\newblock Model reduction for transport-dominated problems via online adaptive
  bases and adaptive sampling.
\newblock {\em SIAM Journal on Scientific Computing}, 42(5):A2803--A2836,
  January 2020.

\bibitem{persson2009discontinuous}
Per-Olof Persson, Javier Bonet, and Jaime Peraire.
\newblock {D}iscontinuous {G}alerkin solution of the {N}avier--{S}tokes
  equations on deformable domains.
\newblock {\em Computer Methods in Applied Mechanics and Engineering},
  198(17-20):1585--1595, 2009.

\bibitem{reiss_shifted_2018}
J.~Reiss, P.~Schulze, J.~Sesterhenn, and V.~Mehrmann.
\newblock The {Shifted} {Proper} {Orthogonal} {Decomposition}: A mode
  decomposition for multiple transport phenomena.
\newblock {\em SIAM Journal on Scientific Computing}, 40(3):A1322--A1344,
  January 2018.

\bibitem{rewienski2003trajectory}
Micha{\l}~Jerzy Rewienski.
\newblock {\em A trajectory piecewise-linear approach to model order reduction
  of nonlinear dynamical systems}.
\newblock PhD thesis, Massachusetts Institute of Technology, 2003.

\bibitem{rim_displacement_2018}
Donsub Rim and Kyle~T. Mandli.
\newblock Displacement interpolation using monotone rearrangement.
\newblock {\em SIAM/ASA Journal on Uncertainty Quantification},
  6(4):1503--1531, January 2018.

\bibitem{rim_transport_2018}
Donsub Rim, Scott Moe, and Randall~J. LeVeque.
\newblock Transport reversal for model reduction of hyperbolic partial
  differential equations.
\newblock {\em SIAM/ASA Journal on Uncertainty Quantification}, 6(1):118--150,
  January 2018.

\bibitem{rim2023manifold}
Donsub Rim, Benjamin Peherstorfer, and Kyle~T Mandli.
\newblock Manifold approximations via transported subspaces: Model reduction
  for transport-dominated problems.
\newblock {\em SIAM Journal on Scientific Computing}, 45(1):A170--A199, 2023.

\bibitem{taddei_registration_2020}
Tommaso Taddei.
\newblock A registration method for model order reduction: Data compression and
  geometry reduction.
\newblock {\em SIAM Journal on Scientific Computing}, 42(2):A997--A1027,
  January 2020.

\bibitem{taddei2021registration}
Tommaso Taddei and Lei Zhang.
\newblock Registration-based model reduction in complex two-dimensional
  geometries.
\newblock {\em Journal of Scientific Computing}, 88:1--25, 2021.

\bibitem{torlo_model_2020}
Davide Torlo.
\newblock Model reduction for advection dominated hyperbolic problems in an
  {ALE} framework: Offline and online phases.
\newblock {\em arXiv:2003.13735 [cs, math]}, March 2020.
\newblock arXiv: 2003.13735.

\bibitem{veroy2003posteriori}
Karen Veroy, Christophe Prud'Homme, Dimitrios Rovas, and Anthony Patera.
\newblock A posteriori error bounds for reduced-basis approximation of
  parametrized noncoercive and nonlinear elliptic partial differential
  equations.
\newblock In {\em 16th AIAA Computational Fluid Dynamics Conference}, page
  3847, 2003.

\bibitem{washabaugh2016faster}
Kyle Washabaugh.
\newblock {\em Faster fidelity for better design: {A} scalable model order
  reduction framework for steady aerodynamic design applications}.
\newblock PhD thesis, Stanford University, 2016.

\bibitem{welper_interpolation_2017}
Garrit Welper.
\newblock Interpolation of functions with parameter dependent jumps by
  transformed snapshots.
\newblock {\em SIAM Journal on Scientific Computing}, 39(4):A1225--A1250,
  January 2017.

\bibitem{welper_transformed_2020}
Garrit Welper.
\newblock Transformed snapshot interpolation with high resolution transforms.
\newblock {\em SIAM Journal on Scientific Computing}, 42(4):A2037--A2061,
  January 2020.

\bibitem{wen2023globally}
Tianshu Wen and Matthew~J. Zahr.
\newblock A globally convergent method to accelerate large-scale optimization
  using on-the-fly model hyperreduction: application to shape optimization.
\newblock {\em Journal of Computational Physics}, 484:112082, 2023.

\bibitem{xu_multi-level_2020}
Jiayang Xu and Karthik Duraisamy.
\newblock Multi-level convolutional autoencoder networks for parametric
  prediction of spatio-temporal dynamics.
\newblock {\em Computer Methods in Applied Mechanics and Engineering},
  372:113379, December 2020.

\bibitem{yano_lp_2019}
Masayuki Yano and Anthony~T. Patera.
\newblock An {LP} empirical quadrature procedure for reduced basis treatment of
  parametrized nonlinear {PDEs}.
\newblock {\em Computer Methods in Applied Mechanics and Engineering},
  344:1104--1123, February 2019.

\bibitem{yue_accelerating_2013}
Yao Yue and Karl Meerbergen.
\newblock Accelerating optimization of parametric linear systems by model order
  reduction.
\newblock {\em SIAM Journal on Optimization}, 23(2):1344--1370, January 2013.

\bibitem{zahr_progressive_2015}
Matthew~J. Zahr and Charbel Farhat.
\newblock Progressive construction of a parametric reduced-order model for
  {PDE}-constrained optimization.
\newblock {\em International Journal for Numerical Methods in Engineering},
  102(5):1111--1135, May 2015.

\bibitem{zahr_implicit_2020}
Matthew~J. Zahr, Andrew Shi, and Per-Olof Persson.
\newblock Implicit shock tracking using an optimization-based high-order
  discontinuous {Galerkin} method.
\newblock {\em Journal of Computational Physics}, 410:109385, June 2020.

\bibitem{zucatti2023adaptive}
Victor Zucatti and Matthew~J. Zahr.
\newblock An adaptive, training-free reduced-order model for
  convection-dominated problems based on hybrid snapshots.
\newblock {\em arXiv preprint arXiv:2301.01718}, 2023.

\end{thebibliography}

\end{document}